\documentclass[12pt]{article}

\usepackage{graphicx, verbatim,mathrsfs, amsmath, amsthm, amssymb, amsfonts,multirow,moreverb}
\newtheorem{defn}{Definition}[section]

\newtheorem{prop}[defn]{Proposition}
\newtheorem{lemma}[defn]{Lemma}
\newtheorem{thm}[defn]{Theorem}
\newtheorem{cor}[defn]{Corollary}

\theoremstyle{remark}
\newtheorem{rem}[defn]{Remark}

\numberwithin{equation}{section}
\newcommand{\com}[2]{\left[\,#1\,,#2\,\right]}
\newcommand{\inprod}[2]{\left(\,#1\,,#2\,\right)}
\newcommand{\bE}{\mathbb{E}}
\newcommand{\bP}{\mathbb{P}}
\newcommand{\bR}{\mathbb{R}}
\newcommand{\bO}{\mathbb{O}}
\newcommand{\bU}{\mathbb{U}}
\newcommand{\bone}{\mathbf{1}}

\newcommand{\footnoteremember}[2]{\footnote{#2}
 \newcounter{#1}
 \setcounter{#1}{\value{footnote}}} \newcommand{\footnoterecall}[1]{\footnotemark[\value{#1}]}

\DeclareMathOperator{\re}{Re}

\DeclareMathOperator{\sgn}{sgn}
\DeclareMathOperator{\pfaffian}{Pf} \DeclareMathOperator{\tr}{tr}
 \DeclareMathOperator{\airy}{Ai}
\newcommand{\rchi}{\raisebox{.4ex}{$\chi$}}
\def\ra{\rightarrow}
\def\iy{\infty}
\def\hf{{1\over 2}}

\def\be{\begin{equation}}
\def\ee{\end{equation}}
\def\ov{\over}
\def\ve{\varepsilon}
\def\la{\lambda}

\begin{document}
\title{Application of Random Matrix Theory to Multivariate Statistics}
\author{\textsc{Momar Dieng}\\ Department of Mathematics \\ University of Arizona \\
Tucson, AZ 85721, USA\\
email: \texttt{momar@math.arizona.edu}
\and
\textsc{Craig A.~Tracy}\\
Department of Mathematics\\
University of California \\
Davis, CA 95616, USA\\
email: \texttt{tracy@math.ucdavis.edu}}
\maketitle
\begin{abstract}
This is an expository account of the edge eigenvalue distributions in random matrix
theory and their application in multivariate statistics.  The emphasis is on the Painlev\'e
representations of these distribution functions.
\end{abstract}
\newpage
\setcounter{tocdepth}{2}
\tableofcontents
\newpage
\section{Multivariate Statistics}
\subsection{Wishart distribution}
The basic problem in statistics is testing the agreement between actual observations and an underlying probability model.
Pearson in 1900 \cite{Pear1} introduced the famous $\rchi^2$ \textit{test} where the sampling distribution  approaches,
as the sample size increases,  
to the $\rchi^2$ distribution.  Recall that if $X_j$ are independent and identically distributed standard normal random variables, $N(0,1)$,
then the distribution of \be \rchi^2_n:=X_1^2+\cdots +X_n^2\label{chiSq}\ee has density
\be
f_{n}(x)=\left \{ \begin{array}{ll} {1\ov 2^{n/2}\, \Gamma(n/2)}\, x^{n/2 -1} e^{-x/2} & \textrm{for}\> x>0,\\
		0 & \textrm{for}\> x\le 0,
		\end{array}\right.\label{chiDensity}
\ee
where $\Gamma(x)$ is the gamma function.

In classical \textit{multivariate  statistics}\footnote{There are many excellent textbooks on multivariate statistics, 
we mention Anderson \cite{Ande1},  Muirhead \cite{Muir1}, and for a shorter
introduction, Bilodeau and Brenner \cite{Bilo1}.}
 it is commonly assumed that the underlying distribution is the multivariate normal distribution.  
 If $X$ is a $p\times 1$-variate normal with $\bE(X)=\mu$ and $p\times p$ covariance matrix
  $\Sigma=\textrm{cov}(X):=\bE\left((X-\mu)\otimes (X-\mu)\right)$,\footnote{If $u$ and $v$ are vectors we denote by $u\otimes v$ the matrix with $(i,j)$ matrix element $u_i v_j$.}
  denoted $N_p(\mu,\Sigma)$, 
 then if $\Sigma>0$  the density function of $X$ is
 \[ f_X(x)=(2\pi)^{-p/2} \left(\det\Sigma\right)^{-1/2}\, \exp\left[-{1\ov2 }\, \left(x-\mu,\Sigma^{-1} (x-\mu)\right)\right ], \>  x\in\bR^p, \]
where $(\cdot,\cdot)$ is the standard inner product on $\bR^p$.

It is convenient to  introduce a matrix notation:  If $X$ is  a $n\times p$ matrix (the \textit{data matrix})  whose rows $X_j$ are independent $N_p(\mu,\Sigma)$ random
variables,
\[ X =\left( \begin{array} {lcl} \longleftarrow &  X_1  &\longrightarrow \\
				             \longleftarrow & X_2    &        \longrightarrow \\
					 & \vdots  & \\
					\longleftarrow & X_n & \longrightarrow
						\end{array}\right), \]
then we say $X$ is $N_p(\bone\otimes \mu, I_n \otimes \Sigma)$ where $\bone=(1,1,\ldots,1)$ and
$I_n$ is the $n\times n$ identity matrix.  We now introduce the multivariate generalization of (\ref{chiSq}).
\begin{defn}
If $A=X^{t} X$, where the $n\times p$ matrix $X$ is $N_p(0,I_n\otimes \Sigma)$, then $A$ is said to have 
\textit{Wishart distribution} with $n$ degrees of freedom and covariance matrix $\Sigma$.  We write $A$
is $W_p(n,\Sigma)$.
\end{defn}
To state the generalization of (\ref{chiDensity}) we first introduce the \textit{multivariate Gamma function}.  If $\mathcal{S}_m^+$ is the space of $p\times p$
positive definite, symmetric matrices, then
\[ \Gamma_p(a):=\int_{{\mathcal{S}_p^+}} e^{-\textrm{tr}(A)}\, \left(\det A\right)^{a-(p+1)/2} \, dA \]
where $\re(a)>(m-1)/2$ and $dA$ is the product Lebesgue measure  of the ${1\ov 2} p(p+1)$ distinct elements of $A$.  By introducing the matrix factorization
$A=T^{t}T$ where $T$ is upper-triangular with positive diagonal elements, one can evaluate this integral in terms of ordinary gamma functions, see, page 62 in \cite{Muir1}.
Note that $\Gamma_1(a)$ is the usual gamma function $\Gamma(a)$.
The basic fact about the Wishart distributions is
\begin{thm}[Wishart \cite{Wish1}]
If $A$ is $W_p(n,\Sigma)$ with $n\ge p$, then the density function of $A$ is
\be
{1\ov 2^{p\, n/2} \,\Gamma_p(n/2) \left(\det \Sigma\right)^{n/2}}\, e^{-{1\ov 2}\, \tr(\Sigma^{-1} A)}\, 
\left(\det A\right)^{(n-p-1)/ 2}.
\label{wishartDensity}
\ee
\end{thm}
For $p=1$ and $\Sigma=1$ (\ref{wishartDensity}) reduces to (\ref{chiDensity}).  The case $p=2$ was obtained by Fisher in 1915 and for general $p$
by Wishart in 1928 using geometrical arguments.   Most modern proofs  follow James \cite{Jame1}.  The importance of the Wishart distribution lies in the
fact that the \textit{sample covariance matrix}, $S$,  is $W_p(n,{1\ov n}\Sigma)$ where
\[ S:= {1\ov n} \sum_{j=1}^N \left(X_i-\overline{X}\right)\otimes \left(X_j - \overline{X}\right),\> N=n+1,\]
and  $X_j$, $j=1,\ldots, N$,  are independent $N_p(\mu,\Sigma)$ random vectors, and $\overline{X}={1\ov N}\sum_j X_j$.

Principle component analysis,\footnote{See, for example, Chap.~9 in \cite{Muir1},  and \cite{John1}
for a discussion of some current issues in principle component analysis.} a multivariate data reduction technique, requires the eigenvalues of the sample covariance matrix; in particular, the largest eigenvalue (the largest principle component variance) is most important.  The next major result
gives the joint density for the eigenvalues of a Wishart matrix.
\begin{thm}[James \cite{Jame2}]
If $A$ is $W_p(n,\Sigma)$ with $n\ge p$ the joint density function of the eigenvalues $\ell_1,\ldots, \ell_p$ of $A$ is
\begin{eqnarray}
{\pi^{p^2/2} 2^{-p n/2} \left(\det\Sigma\right)^{-n/2}\ov \Gamma_p(p/2) \Gamma_p(n/2)}\, \prod_{j=1}^p \ell_j^{(n-p-1)/2}\, \prod_{j<k}(\ell_j-\ell_k) \nonumber \\
\hspace{5ex}\cdot \int_{\bO(p)} e^{-{1\ov 2}\tr(\Sigma^{-1} H L H^{t})}\, dH \label{WishartEigDensity}\end{eqnarray}
where $\bO(p)$ is the orthogonal group of $p\times p$ matrices, $dH$ is normalized  Haar measure and $L$ is the diagonal matrix $\textrm{diag}(\ell_1,\ldots,\ell_p)$.
(We take $\ell_1>\ell_2>\cdots > \ell_p$.)
\end{thm}
\begin{rem}
The difficult part of this density function is the integral over the orthogonal group $\bO(p)$.  There is no known closed formula 
for this integral though  James and
Constantine (see Chap.~7 in \cite{Muir1} for references) developed the theory of \textit{zonal polynomials}  which
 allow one to write infinite series expansions for this integral.
However, these expansions converge slowly;  and zonal polynomials themselves,  lack explicit formulas such
as are available for Schur polynomials.    For complex
Wishart matrices, the group integral is over the unitary group $\bU(p)$;  and this integral can be evaluated
 using the Harish-Chandra-Itzykson-Zuber
integral \cite{Zinn1}.
\end{rem}
There is one important case where the integral can be (trivially) evaluated.
\begin{cor} 
If $\Sigma=I_p$, then the joint density (\ref{WishartEigDensity}) simplifies to
\be
{\pi^{p^2/2} 2^{-p n/2} \left(\det\Sigma\right)^{-n/2}\ov \Gamma_p(p/2) \Gamma_p(n/2)}\, \prod_{j=1}^p \ell_j^{(n-p-1)/2}\, \exp\left(-{1\ov 2}\sum_j \ell_j\right) \,\prod_{j<k}(\ell_j-\ell_k). 
\ee
\end{cor}
\subsection{An example with $\Sigma\not=c I_p$}
This section uses the theory of zonal polynomials as can be found in Chap.~7 of Muirhead \cite{Muir1} or Macdonald \cite{Macd1}.  This section is not used in
the remainder of the chapter.  Let $\la=(\la_1,\ldots, \la_p)$ be a partition into not more than
$p$ parts.  We let $C_\la(Y)$ denote the zonal polynomial of $Y$ corresponding to $\la$.  It is a symmetric, homogeneous polynomial of degree $\vert \la\vert$
in the eigenvalues $y_1,\ldots, y_p$ of $Y$.  The normalization we adopt is defined by
\[ \left(\tr Y\right)^k=(y_1+\cdots+y_p)^k=\sum_{{\la\vdash k \atop \ell(\la)\le p}} C_\la(Y). \]
The fundamental integral formula for zonal polynomials is\footnote{See, for example, Theorem 7.2.5 in \cite{Muir1}.}
\begin{thm} Let $X,Y\in \mathcal{S}_p^{+}$, then
\be \int_{\bO(p)} C_\la (X H Y H^{t})\,  dH = {C_\la(X) C_\la(Y)\ov C_\la(I_p)} \label{zonalIntegral}\ee
where $dH$ is normalized Haar measure.
\end{thm}
By expanding the exponential and using (\ref{zonalIntegral}) it follows that
\be
\int_{\bO(p)} \exp\left(z \tr(XHYH^{t})\right)\, dH =
 \sum_{k\ge 0} {z^k\ov k!}\, \sum_{{\la\vdash k \atop \ell(\la)\le p}} {C_\la(X) C_\la(Y)\ov C_\la(I_p)}\, . \label{zonalExpansion}
\ee

 We examine (\ref{zonalExpansion}) for the special case ($|\rho|<1$)
\[ \Sigma= (1-\rho) I_p + \rho \, \bone\otimes\bone =\left(\begin{array}{lllll}
						1& \rho & \rho & \cdots & \rho \\
						\rho & 1 & \rho & \cdots &\rho \\
						\vdots & && \vdots \\
						\rho & \rho & \rho &\cdots & 1\end{array}\right).\]
We have
\[ \Sigma^{-1} = (1-\rho)^{-1}\, I_p - {\rho\ov (1-\rho)(1+(p-1)\rho)}\, \bone\otimes\bone \]
and
\[ \det\Sigma=(1-\rho)^{p-1}\left(1+(p-1)\rho\right).\]
For this choice of $\Sigma$, let $Y=\alpha \bone\otimes\bone$ where $\alpha=\rho/((2(1-\rho)(1+(p-1)\rho))$, then
\begin{eqnarray*}
\int_{\bO(p)} e^{-{1\ov 2}\tr(\Sigma^{-1} H L H^{t})}\, dH &=& e^{-{1\ov 2(1-\rho)}\sum_j \la_j}\, \int_{\bO(p)}e^{\tr(YH L H^{t})}\, dH\\
&=&  e^{-{1\ov 2(1-\rho)}\sum_j \la_j}\, \sum_{k\ge 0} {C_{(k)}(\alpha \bone\otimes\bone) C_{(k)}(\Lambda)\ov k! C_{(k)}(I_p)}
\end{eqnarray*}
where we used the fact that the only partition $\la\vdash k$ for which $C_\la(Y)$ is nonzero is $\la=(k)$.  And for this
partition,  
$C_{(k)}(Y)= \alpha^k p^k$.  Define the symmetric functions $g_n$\footnote{In the theory of zonal polynomials,
the $g_n$ are the analogue of the   complete symmetric functions $h_n$.} by
\[ \prod_{j\ge 1} (1-x_j y)^{-1/2} = \sum_{n\ge 0} g_n(x) y^n,\]
then it is known that \cite{Macd1}
\[ C_{(k)}(L) = {2^{2k} (k!)^2\ov (2k)!} \, g_k(L).\]
Using the known value of $C_{(k)}(I_p)$ we find
\[ \int_{\bO(p)}e^{-{1\ov 2}\tr(\Sigma^{-1} H L H^{t})} \, dH=e^{-{1\ov 2(1-\rho)}\sum_j \la_j}\,\sum_{k\ge 0} {(\alpha p)^k\ov ({1\ov 2} p)_k}\, g_k(L)\]
where $(a)_k=a(a+1)\cdots (a+k-1)$ is the Pochammer symbol.

\section{Edge Distribution Functions}
\subsection{Summary of Fredholm determinant representations}\label{sec:edgeDistr}
In this section we define three Fredholm determinants  from which the
edge eigenvalue distributions,   for the three symmetry classes orthogonal, unitary
and symplectic, will ensue.  This section follows \cite{Trac3, Trac2, Trac5}; see also, \cite{Ferr1, Forr2}.

In the unitary case ($\beta=2$),  define the  trace class operator $K_{2}$  on $L^2(s,\iy)$
with  \textit{Airy kernel}
\begin{equation}\label{kairy}
  K_{\airy}(x,y):=\frac{\airy(x)\airy^{\prime}(y)-\airy^{\prime}(x)\airy(y)}{x-y}=\int_{0}^\iy \airy(x+z)\airy(y+z)\, dz
\end{equation}
and associated Fredholm determinant, $0\le \la \le 1$, 
\begin{equation}
  D_{2}(s,\lambda)=\det(I-\lambda\,K_{2}).\label{fred2}
\end{equation}
Then we introduce the distribution functions
\be F_2(s)=F_2(s,1) = D_2(s,1), \label{TW2} \ee
and for $m\ge 2$,  the distribution functions $F_2(s,m)$ are defined recursively  below by (\ref{lderiv}).

In the symplectic case ($\beta=4$), we define the trace class operator $K_4$ on $L^2(s,\iy)\oplus L^2(s,\iy)$ with  matrix kernel 
\be
K_4(x,y):={1\ov 2}\, \left(\begin{array}{cc}
	S_4(x,y) & SD_4(x,y)\\
	IS_4(x,y) & S_4(y,x) \end{array}\right)
	\label{K4}
\ee
where 
\begin{eqnarray*}
S_4(x,y)&=&K_{\airy}(x,y)-{1\ov 2} \airy(x,y)\, \int_y^{\iy} \airy(z)\, dz,\\
SD_4(x,y)&=&-\partial_y K_{\airy}(x,y) -{1\ov 2}\, \airy(x) \airy(y),\\
IS_4(x,y)&=&-\int_x^{\iy} K_{\airy}(z,y)\, dz + {1\ov 2}\, \int_x^{\iy} \airy(z)\, dz \cdot  \int_y^{\iy} \airy(z),
\end{eqnarray*}
and the associated Fredholm determinant, $0\le\la\le 1$,
\begin{equation}
  D_{4}(s,\lambda)=\det(I-\lambda\,K_{4}).\label{fred4}
\end{equation}
Then we introduce the distribution functions (note the square root)
\be F_4(s)=F_4(s,1) =\sqrt{D_4(s,1)}\, , \label{TW4} \ee
and for $m\ge 2$,  the distribution functions $F_4(s,m)$ are defined recursively  below by (\ref{betaderiv}).

In the orthogonal case ($\beta=1$), we introduce the matrix kernel
\be
K_1(x,y):=\left(\begin{array}{cc}
	S_1(x,y) & SD_1(x,y)\\
	IS_1(x,y)-\ve(x,y)  & S_1(y,x) \end{array}\right)
	\label{K1}
\ee
where
\begin{eqnarray*}
S_1(x,y)&=&K_{\airy}(x,y)-{1\ov 2} \airy(x)\, \left(1-\int_y^{\iy} \airy(z)\, dz\right),\\
SD_1(x,y)&=&-\partial_y K_{\airy}(x,y) -{1\ov 2}\, \airy(x) \airy(y),\\
IS_1(x,y)&=&-\int_x^{\iy} K_{\airy}(z,y)\, dz + {1\ov 2}\,\left( \int_y^{x} \airy(z)\, dz +
\int_x^{\iy}\airy(z)\, dz\, \cdot\, \int_y^{\iy} \airy(z)\, dz\right),\\
\ve(x-y)&=&{1\ov 2} \, \textrm{sgn}(x-y).
\end{eqnarray*}
The  operator $K_1$ on $L^2(s,\iy)\oplus L^2(s,\iy)$ with  this matrix kernel is \textit{not}
trace class due to the presence of $\ve$.  As discussed in \cite{Trac5}, one must use the weighted
space $L^2(\rho)\oplus L^2(\rho^{-1})$, $\rho^{-1}\in L^1$.   Now the determinant is the 2-determinant,
 \begin{equation}
  D_{1}(s,\lambda)={\det}_2(I-\lambda\,K_{1}\rchi_J)\label{fred1}
\end{equation}
where $\rchi_J$ is the characteristic function of the interval $(s,\iy)$.
 We introduce the distribution functions (again note the square root)
\be F_1(s)=F_1(s,1) =\sqrt{D_1(s,1)}\, , \label{TW1} \ee
and for $m\ge 2$,  the distribution functions $F_1(s,m)$ are defined recursively  below by (\ref{betaderiv}).
This is the first indication that the determinant $D_1(s,\la)$ might be more subtle than either
$D_2(s,\la)$ or $D_4(s,\la)$.

\subsection{Universality theorems}
Suppose  $A$  is $W_{p}(n,I_p)$ with eigenvalues $\ell_{1}>\cdots >\ell_{p}$. We define scaling constants
\begin{equation*}
\mu_{_{np}}=\left(\sqrt{n-1}+\sqrt{p}\right)^{2}\, ,\quad \sigma_{_{np}}=\left(\sqrt{n-1}+\sqrt{p}\right)\left(\frac{1}{\sqrt{n-1}}+\frac{1}{\sqrt{p}}\right)^{1/3} \, . 
\end{equation*}
The following  theorem establishes,  under the null hypothesis $\Sigma=I_p$,  that the  largest principal component
variance, $\ell_1$,   converges in law to $F_1$.

\begin{thm}[Johnstone, \cite{John1}]
\label{sec:wishart-ensembles-3}
If  $n, p\ra\iy$ such that $n/p\to\gamma, 0<\gamma<\infty$, then
\begin{equation*}
  \frac{\ell_{1}-\mu_{_{np}}}{\sigma_{_{np}}}\xrightarrow{\mathscr{D}}F_{1}(s,1).   
\end{equation*}
\end{thm}
\noindent Johnstone's theorem generalizes to  the $m^{th}$ largest eigenvalue.
\begin{thm}[Soshnikov, \cite{Sosh2}]
\label{sec:wishart-ensembles-4}
If $n, p\ra\iy$ such that $n/p\to\gamma, 0<\gamma<\infty$, then
\begin{equation*}
  \frac{\ell_{m}-\mu_{_{np}}}{\sigma_{_{np}}}\xrightarrow{\mathscr{D}}F_{1}(s,m),\> m=1,2,\ldots.
\end{equation*}
\end{thm}
\noindent Soshnikov proved his result under the additional assumption $n-p=\textrm{O}(p^{1/3})$.
We remark  that a straightforward generalization of Johnstone's proof \cite{John1} together with
results of Dieng \cite{Dien1} show this restriction can be removed.
Subsequently, El Karoui \cite{Elka1} extended  Theorem \ref{sec:wishart-ensembles-4} to  $0<\gamma\leq \infty$. 
The extension to $\gamma=\iy$ is important for modern statistics where $p\gg n$ arises in applications.

Going further, Soshnikov lifted  the Gaussian assumption, again establishing  a $F_1$ universality theorem. In order to state the generalization precisely,
 let us redefine the $n\times p$ matrices $X=\left\{x_{i,j}\right\}$ such that $A=X^{t}X$ to satisfy
\begin{enumerate}
\item $\bE(x_{ij})=0$, $\bE(x_{ij}^{2})=1$.
\item The random variables $x_{ij}$ have symmetric laws of distribution.
\item All even moments of $x_{ij}$ are finite, and they 
 decay at least as fast as  a Gaussian at infinity: $\bE(x_{ij}^{2m})\leq (\textrm{const}\,m)^{m}$.
\item $n-p=\textrm{O}(p^{1/3})$.
\end{enumerate}
With these assumptions, 
\begin{thm}[Soshnikov, \cite{Sosh2}]
\begin{equation*}
  \frac{\ell_{m}-\mu_{_{np}}}{\sigma_{_{np}}}\xrightarrow{\mathscr{D}}F_{1}(s,m),\>  m=1,2,\ldots.
\end{equation*} \end{thm}
\noindent It is an important open problem to remove the restriction $n-p=\textrm{O}(p^{1/3})$.

For real symmetric matrices, 
Deift and Gioev \cite{Deif3}, building on the work of Widom \cite{Wido1},  proved   $F_1$ universality   when the Gaussian weight function $\exp(-x^{2})$ is replaced by  $\exp(-V(x))$ 
where  $V$ is an even degree polynomial with positive leading coefficient.

Table~\ref{table} in Section \ref{sec:numerics}  displays a comparison of  the percentiles of the $F_{1}$ distribution   with percentiles of empirical Wishart distributions. Here $\ell_{j}$ denotes the $j^{th}$ largest eigenvalue in the Wishart Ensemble. The percentiles in the $\ell_{j}$ columns were obtained by finding the ordinates corresponding to the $F_{1}$--percentiles listed in the first column, and computing the proportion of eigenvalues lying to the left of that ordinate in the empirical distributions for the $\ell_{j}$. The bold entries correspond to the levels of confidence commonly used in statistical applications. The reader should compare Table~\ref{table} to  similar ones in \cite{Elka1, John1}.

\section{Painlev\'e Representations: A Summary}
\begin{figure}[htbp] \label{TWdensities}
\centering
  \includegraphics{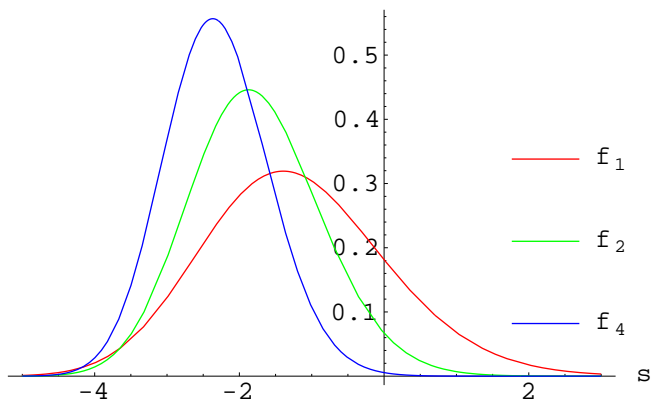}
  \caption{Tracy--Widom Density Functions}
\end{figure}

The Gaussian $\beta$--ensembles are probability spaces on $N$-tuples of random variables $\{\ell_{1},\ldots, \ell_{N}\}$, with joint density functions $P_{\beta}$ 
given by\footnoteremember{norm}{In many places in the random matrix theory literature, the parameter $\beta$ (times $1/2$) appears in front of the summation inside the exponential factor
(\ref{jointdensity}), in addition to being the power of the Vandermonde determinant. That convention originated in \cite{Meht1}, and was justified by the alternative physical and very useful interpretation of \eqref{jointdensity} as a one--dimensional Coulomb gas model. In that language the potential $W=\frac{1}{2}\sum_{i} l_{i}^{2}-\sum_{i<j}\ln|l_{i}-l_{j}|$ and $P_{\beta}^{(N)}(\vec{\ell}\,\,)=C\exp(-W/kT)=C\exp(-\beta\,W)$, so that $\beta=(k\,T)^{-1}$ plays the role of inverse temperature. However, by an appropriate choice of specialization in Selberg's integral, it is possible to remove the $\beta$ in the exponential weight, at the cost of redefining the normalization constant $C_{\beta}^{(N)}$. We choose the latter convention in this work since we will not need the Coulomb gas analogy. Moreover, with computer simulations and statistical applications in mind, this will in our opinion make later choices of standard deviations, renormalizations, and scalings  more transparent. It also allows us to dispose of the   $\sqrt{2}$ that is often present in  $F_4$. }
\begin{equation}\label{jointdensity}
  P_{\beta}(\ell_{1},\ldots,\ell_{N})= P_{\beta}^{(N)}(\vec{\ell}\,\,) = C_{\beta}^{(N)}\,\exp\left[-\sum_{j=1}^{N}\ell_{j}^{\,2}\right]\prod_{j<k}|\ell_{j}-\ell_{k}|^{\beta}.
\end{equation}
The $C_{\beta}^{(N)}$ are normalization constants, given by
\begin{equation}
  \label{eq:34}
C_{\beta}^{(N)} = \pi^{-N/2}\,2^{-N-\beta\,N(N-1)/4}\cdot\prod_{j=1}^{N}\frac{\Gamma(1 + \gamma)\,\Gamma(1 + \frac{\beta}{2})}{\Gamma(1 + \frac{\beta}{2}\,j)}
\end{equation}
By setting $\beta=1, 2, 4$ we recover the (finite $N$) \emph{Gaussian Orthogonal Ensemble} ($\textrm{GOE}_{N}$), \emph{Gaussian Unitary Ensemble} ($\textrm{GUE}_{N}$), and
 \emph{Gaussian Symplectic Ensemble} ($\textrm{GSE}_{N}$), respectively. For the remainder of the chapter we restrict to these
three cases, and refer the reader to \cite{Dumi1} for recent results on the general $\beta$ case. Originally the $\ell_{j}$ are eigenvalues of randomly chosen matrices from corresponding matrix ensembles,
 so we will henceforth refer to them as eigenvalues. With the eigenvalues ordered so that $\ell_{j}\geq\ell_{j+1}$, define
\begin{equation}
  \hat{\ell}_{m}^{(N)}=\frac{\ell_{m}-\sqrt{2\,N}}{2^{-1/2}\,N^{-1/6}},
\end{equation}
to be the rescaled $m^{th}$ eigenvalue measured from edge of spectrum. For the  largest eigenvalue in the $\beta$--ensembles (proved only in the $\beta=1,2,4$ cases) we have
\begin{equation}
  \hat{\ell}_{1}^{\,(N)}\xrightarrow{\mathscr{D}}\hat{\ell}_{1},
\end{equation}
whose law is given by the Tracy--Widom distributions.\footnoterecall{norm}
\begin{thm}[Tracy, Widom \cite{Trac3,Trac2}]
  \begin{eqnarray}
    F_{2}(s)&=&\bP_{2}(\hat{\ell}_{1}\leq s)=\exp\left[-\int_{s}^{\infty}(x-s)\,q^{2}(x)d\,x\right],\label{guemax}\\
    F_{1}(s)&=&\bP_{1}(\hat{\ell}_{1}\leq s)=\left(F_{2}(s)\right)^{1/2} \, \exp\left[-\hf \int_{s}^{\infty}q(x)d\,x\right], \label{goemax}\\
    F_{4}(s)&=&\bP_{4}(\hat{\ell}_{1}\leq s)=\left(F_{2}(s)\right)^{1/2} \, \cosh\left[-\frac{1}{2}\int_{s}^{\infty}q(x)d\,x\right].
    \label{gsemax}
  \end{eqnarray}
\end{thm}
The function $q$ is the unique (see \cite{Clar1,Hasti1}) solution to the Painlev\'e II equation 
\begin{equation}\label{pII}
  q'' = x\,q + 2\,q^{3},
\end{equation}
such that $q(x)\sim \airy(x)$ as $x\to\infty$, where $\airy(x)$ is the solution to the Airy equation which decays like $\frac{1}{2}\,\pi^{-1/2}\,x^{-1/4}\,\exp\left(-\frac{2}{3}\,x^{3/2}\right)$ at $+\infty$. The density functions $f_{\beta}$ corresponding to the $F_\beta$ are graphed in Figure \ref{TWdensities}.\footnote{Actually,
for $\beta=4$,
the density of  $F_{4}(\sqrt{2}\, s)$ is graphed  to agree with Mehta's original normalization \cite{Meht1}
as well as with \cite{Trac2}.}
Let $F_{2}(s,m)$ denote the distribution for the $m^{th}$ largest eigenvalue in GUE. Tracy and Widom showed \cite{Trac3} that if we define $F_{2}(s,0)\equiv 0$, then
\begin{equation}
  F_{2}(s,m+1) -  F_{2}(s,m) = \frac{(-1)^{m}}{m\,!}\frac{d^{m}}{d\,\lambda^{m}}\,D_{2}(s,\lambda)\big{\vert}_{\lambda=1}\,,\quad m\geq 0, \label{lderiv}
\end{equation}
where (\ref{fred2}) has the Painlev\'e representation
\begin{equation}\label{D2}
  D_{2}(s,\lambda)=\exp\left[-\int_{s}^{\infty}(x-s)\,q^{2}(x,\lambda)d\,x \right],
\end{equation}
and $q(x,\lambda)$ is the solution to \eqref{pII} such that $q(x,\lambda)\sim~\sqrt{\lambda}\,\airy(x)$ as $x\to\infty$. 
The same combinatorial argument used to obtain the recurrence \eqref{lderiv} in the $\beta=2$ case also works for the $\beta=1,4$ cases, leading to
\begin{equation}\label{betaderiv}
  F_{\beta}(s,m+1) -  F_{\beta}(s,m) = \frac{(-1)^{m}}{m\,!}\frac{d^{m}}{d\,\lambda^{m}}\,D_{\beta}^{1/2}(s,\lambda)\big{\vert}_{\lambda=1}\,,\quad m\geq 0,\, \beta=1,4,
\end{equation}
where $F_{\beta}(s,0)\equiv 0$. Given the similarity in the arguments up to this point and comparing \eqref{D2} to \eqref{guemax}, it is natural to conjecture that $D_{\beta}(s,\lambda), \beta=1,4,$ can be obtained simply by replacing $q(x)$ by $q(x,\lambda)$ in \eqref{goemax} and \eqref{gsemax}.

That this is not the case for $\beta=1$  was shown by Dieng \cite{Dien1, Dien2}.  A hint that $\beta=1$ is different comes from the following interlacing theorem.
\begin{thm}[Baik, Rains \cite{Baik2}]\label{baikthm}
  In the appropriate  scaling limit, the distribution of the largest eigenvalue in GSE corresponds to that of the second largest in GOE. More generally, the joint distribution of every second eigenvalue in the GOE coincides with the joint distribution of all the eigenvalues in the GSE, with an appropriate number of eigenvalues.
\end{thm}
This  interlacing property  between GOE and GSE had long been in the literature, and had in fact been noticed by Mehta and Dyson \cite{Meht3}. In this context, Forrester and Rains \cite{Forr1} classified all weight functions for which alternate eigenvalues taken from an orthogonal ensemble form a corresponding symplectic ensemble, and similarly those for which alternate eigenvalues taken from a union of two orthogonal ensembles form an unitary ensemble. 
The following theorem gives explicit formulas for $D_1(s,\la)$ and $D_4(s,\la)$; and hence, from (\ref{betaderiv}), a recursive procedure to determine $F_{1}(\cdot,m)$ and $F_{4}(\cdot,m)$ for $m\ge 2$.
\begin{thm}[Dieng \cite{Dien1, Dien2}]\label{mainthm}
  In the edge scaling limit, the distributions for the $m^{th}$ largest eigenvalues in the \textrm{GOE} and \textrm{GSE} satisfy the recurrence \eqref{betaderiv} with\footnoterecall{norm}
 
  \begin{equation}\label{goedet}
    D_{1}(s,\lambda)=D_{2}(s,\tilde{\lambda})\,\frac{\lambda - 1 - \cosh{\mu(s,\tilde{\lambda})} + \sqrt{\tilde{\lambda}}\,\sinh{\mu(s,\tilde{\lambda})}}{\lambda - 2},
  \end{equation}
  \begin{equation}\label{gsedet}
    D_{4}(s,\lambda)=D_{2}(s,\lambda)\,\cosh^{2}\left(\frac{\mu(s,\lambda)}{2}\right),
  \end{equation}
  where
  \begin{equation}
    \mu(s,\lambda):=\int_{s}^{\infty}q(x,\lambda)d\,x, \qquad \tilde{\lambda}:=2\,\lambda-\lambda^{2},
  \end{equation}
  and $q(x,\lambda)$ is the solution to \eqref{pII} such that $q(x,\lambda)\sim~\sqrt{\lambda}\,\airy(x)$ as $x\to\infty$.
\end{thm}
\noindent Note the appearance of $\tilde\lambda$ in the arguments
on the right hand side of (\ref{goedet}).  In Fig.~2 we compare the densities $f_1(s,m)$, $m=1,\ldots, 4$, with finite $N$ GOE simulations.
This last theorem also provides a new proof of the Baik-Rains interlacing theorem.
\begin{cor}[Dieng \cite{Dien1, Dien2}]
  \begin{equation}\label{interlacingcor}
    F_{4}(s,m)= F_{1}(s,2\,m),\quad m\geq 1.
  \end{equation}
\end{cor}
The proofs of these  theorems occupy
the bulk of the remaining part of the chapter.  In the last section, we present an efficient numerical scheme to compute $F_{\beta}(s,m)$ and the associated density functions $f_{\beta}(s,m)$. We implemented this scheme using MATLAB\texttrademark\,,\footnote{MATLAB\texttrademark\, is a registered trademark of The MathWorks, Inc., 3 Apple Hill Drive, Natick, MA 01760-2098; Phone: 508-647-7000; Fax: 508-647-7001. Copies
of the code are available by contacting the first author.} and compared the results to simulated Wishart distributions.
  \begin{figure}\label{nextFigure}
  \centering
    \includegraphics[height=65mm]{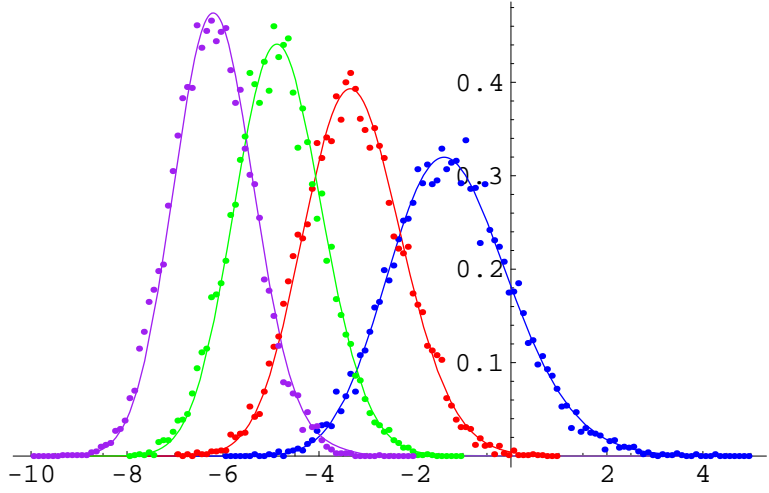}
    \caption{$10^{4}$ realizations of $10^{3}\times 10^{3}$ GOE matrices; the solid curves are, from right to left, the theoretical limiting densities for the first through fourth largest eigenvalue.}
  \end{figure}

\section{Preliminaries}

\subsection{Determinant matters}

We gather in this short section more or less classical results for further reference.
\begin{thm}\label{vandthm}
  \begin{equation*}
    \prod_{0\leq j<k\leq N}(x_{j}-x_{k})^{4}=\det\left(x_{k}^{j}\quad j\,x_{k}^{j-1}\right)_{\substack{j=0,\dots,2\,N-1 \\ k=1,\ldots,N}}
  \end{equation*}
\end{thm}

\begin{thm}\label{detthm} 
  If $A, B$ are Hilbert--Schmidt operators on a general\hspace{.2ex}\footnote{See \cite{Gohb1} for proof.} Hilbert space $\mathcal H$, then
  \begin{equation*}
    \det(I + A\,B)=\det(I + B\,A).
  \end{equation*}
\end{thm}
\begin{thm}[de Bruijn, \cite{Debr1}]\label{debruijnThm}
  \begin{eqnarray}\label{debruijn1} 
    \begin{aligned} 
      \int\cdots\int\det(\varphi_{j}(x_{k}))_{_{1\leq j,k\leq N}}\cdot &\det(\psi_{j}(x_{k}))_{_{1\leq j,k\leq N}}\,d\,\mu(x_{1})\cdots d\,\mu(x_{N}) \\
      & = N!\,\det\left(\int\varphi_{j}(x)\,\psi_{k}(x)\,d\,\mu(x)\right)_{_{1\leq j,k\leq N}},
    \end{aligned}
  \end{eqnarray}
  \begin{eqnarray}\label{debruijn2} 
    \begin{aligned} 
      \underset{x_{1}\leq\ldots\leq x_{N}}{\int\cdots\int} & \det(\varphi_{j}(x_{k}))_{_{1\leq j,k\leq N}}  d\,\mu(x_{1})\cdots d\,\mu(x_{N}) \\
      & = \pfaffian\left(\int\int\sgn(x-y)\varphi_{j}(x)\,\varphi_{k}(x)\,d\,\mu(x)\,d\,\mu(y)\right)_{_{1\leq j,k\leq N}},
    \end{aligned}
  \end{eqnarray}
  \begin{eqnarray}\label{debruijn3}
    \begin{aligned} 
      \int\cdots\int & \det(\varphi_{j}(x_{k})\quad \psi_{j}(x_{k}))_{_{\substack{1\leq j\leq 2\,N\\ 
            1\leq  k \leq N}}} d\,\mu(x_{1})\cdots d\,\mu(x_{N}) \\
      & = (2\,N)!\,\pfaffian\left(\int \varphi_{j}(x)\psi_{k}(x)-\varphi_{k}(x)\psi_{j}(x)\,d\,\mu(x)\right)_{_{1\leq j,k\leq 2\,N}},
    \end{aligned}
  \end{eqnarray}
\end{thm}
\noindent where $\pfaffian$ denotes the Pfaffian. The last two integral identities were discovered by de Bruijn \cite{Debr1} in an attempt to generalize the first one. The first and last are valid in general measure spaces. In the second identity, the space needs to be ordered. In the last identity, the left hand side determinant is a $2\,N\times 2\,N$ determinant whose columns are alternating columns of the $\varphi_{j}$ and $\psi_{j}$ (i.e. the first four columns are $\{\varphi_{j}(x_{1})\}$, $\{\psi_{j}(x_{1})\}$, $\{\varphi_{j}(x_{2})\}$, $\{\psi_{j}(x_{2})\}$, respectively for $j=1,\ldots,2\,N$), hence the notation, and asymmetry in indexing.

A large portion of the foundational theory of random matrices, in the case of invariant measures, can be developed
from Theorems~\ref{detthm} and \ref{debruijnThm} as was demonstrated in \cite{Trac1, Wido1}.

\subsection{Recursion formula for the eigenvalue distributions} 

With the joint density function defined as in \eqref{jointdensity}, let $J$ denote the interval $(t,\infty)$, and $\rchi=\rchi_{_{J}}(x)$ its characteristic function.\footnote{Much of what is said here is still valid if $J$ is taken to be a finite union of open intervals in $\mathbb R$ (see \cite{Trac7}). However, since we will only be interested in edge eigenvalues we restrict ourselves to $(t,\infty)$ from here on.} We denote by $\tilde{\rchi}=1-\rchi$ the characteristic function of the complement of $J$, and define $\tilde{\rchi}_{_{\lambda}}=1-\lambda\,\rchi$. Furthermore, let $E_{\beta,N}(t,m)$ equal the probability that exactly the $m$ largest eigenvalues of a matrix chosen at random from a (finite $N$) $\beta$--ensemble lie in $J$. We also define
\begin{equation}
\label{eq:23}
  G_{\beta,N}(t,\lambda)= \underset{x_{i}\in \mathbb R}{\int\cdots\int} \tilde{\rchi}_{_{\lambda}}(x_{1})\cdots \tilde{\rchi}_{_{\lambda}}(x_{N})\,P_{\beta}(x_{1},\ldots,x_{N})\,d\,x_{1}\cdots d\,x_{N}.
\end{equation}
For $\lambda=1$ this is just $E_{\beta,N}(t,0)$, the probability that no eigenvalues lie in $(t,\infty)$, or equivalently the probability that the largest eigenvalue is less than $t$. In fact we will see in the following propositions that $G_{\beta,N}(t,\lambda)$ is in some sense a generating function for $E_{\beta,N}(t,m)$.
\begin{prop}
  \begin{equation}
    G_{\beta,N}(t,\lambda)=\sum_{k=0}^{N}(-\lambda)^{k}\binom{N}{k}\underset{x_{i}\in J}{\int\cdots\int} P_{\beta}(x_{1},\ldots,x_{N})\,d\,x_{1}\cdots d\,x_{N}.
  \end{equation}
\end{prop}
\begin{proof}
  Using the definition of the $\tilde{\rchi}_{_{\lambda}}(x_{1})$ and multiplying out the integrand of \eqref{eq:23} gives 
  \begin{equation*}
    G_{\beta,N}(t,\lambda) =  \sum_{k=0}^{N}(-\lambda)^{k}\underset{x_{i}\in \mathbb R}{\int\cdots\int} e_{k}(\rchi(x_{1}),\ldots, \rchi(x_{N}))P_{\beta}(x_{1},\ldots,x_{N})\,d\,x_{1}\cdots d\,x_{N},    
  \end{equation*}
  where, in the notation of \cite{Stan2}, $e_{k}=m_{1^{k}}$ is the $k^{th}$ elementary symmetric function. Indeed each term in the summation arises from picking $k$ of the $\lambda\,\rchi$-terms, each of which comes with a negative sign, and $N-k$ of the $1$'s. This explains the coefficient $(-\lambda)^{k}$. Moreover, it follows that $e_{k}$ contains $\binom{N}{k}$ terms. Now the integrand is symmetric under permutations of the $x_{i}$. Also if $x_{i}\not\in J$, all corresponding terms in the symmetric function are $0$, and they are $1$ otherwise. Therefore we can restrict the integration to $x_{i}\in J$, remove the characteristic functions (hence the symmetric function), and introduce the binomial coefficient to account for the identical terms up to permutation.
\end{proof}
\begin{prop}
  \begin{equation}
    \left. E_{\beta,N}(t,m)=\frac{(-1)^{m}}{m\,!}\,\frac{d^{m}}{d\,\lambda^{m}}\,G_{\beta,N}(t,\lambda)\right|_{\lambda=1}\,,\quad m\geq 0.
  \end{equation}
\end{prop}
\begin{proof}
  This is proved by induction. As noted above, $E_{\beta,N}(t,0)=G_{\beta,N}(t,1)$ so it holds for the degenerate case $m=0$. When $m=1$ we have
  \begin{equation*}
    \begin{aligned}
      \left.-\frac{d}{d\,\lambda}\,G_{\beta,N}(t,\lambda)\right|_{\lambda=1} &=\left.-\frac{d}{d\,\lambda}\, \int\cdots\int \tilde{\rchi}_{_{\lambda}}(x_{1})\cdots \tilde{\rchi}_{_{\lambda}}(x_{n})\,P_{\beta}^{(N)}(\vec{x})\,d\,x_{1}\cdots d\,x_{n}\right|_{\lambda=1} \\
      &= -\sum_{j=1}^{N}- \int\cdots\int \tilde{\rchi}(x_{1})\cdots\tilde{\rchi}(x_{j-1})\,\rchi(x_{j})\,\tilde{\rchi}(x_{j+1})\cdots\\
      &   \qquad\qquad \cdots\tilde{\rchi}(x_{N})\, P_{\beta}^{(N)}(\vec{x})\,d\,x_{1}\cdots d\,x_{N}.
    \end{aligned}    
  \end{equation*}
  The integrand is symmetric under permutations so we can make all terms look the same. There are $N=\binom{N}{1}$ of them so we get
  \begin{equation*}
    \begin{aligned}
      \left.-\frac{d}{d\,\lambda}\,G_{\beta,N}(t,\lambda)\right|_{\lambda=1} &=\left.\binom{N}{1}\,\int\cdots\int \rchi(x_{1})\,\tilde{\rchi}(x_{2})\cdots \tilde{\rchi}(x_{N})\,P_{\beta}^{(N)}(\vec{x})\,d\,x_{1}\cdots d\,x_{N}\right|_{\lambda=1} \\
      &=\binom{N}{1}\,\int\cdots\int \rchi(x_{1})\,\rchi(x_{2})\cdots \rchi(x_{N})\,P_{\beta}^{(N)}(\vec{x})\,d\,x_{1}\cdots d\,x_{N}\\
      &= E_{\beta,N}(t,1).
    \end{aligned}  
  \end{equation*}
  When $m=2$ then
  \begin{equation*}
    \begin{aligned}
      \frac{1}{2}&\left(-\frac{d}{d\,\lambda}\right)^{2} \, \left. G_{\beta,N}(t,\lambda)\right|_{\lambda=1} \\
      & = \frac{N}{2}\,\sum_{j=2}^{N} \int\cdots\int \rchi(x_{1})\tilde{\rchi}(x_{2})\cdots\tilde{\rchi}(x_{j-1})\,\rchi(x_{j})\,\tilde{\rchi}(x_{j+1})\cdots \\
      &  \qquad\qquad \cdots\left.\tilde{\rchi}(x_{N})\, P_{\beta}^{(N)}(\vec{x})\,d\,x_{1}\cdots d\,x_{N}\right|_{\lambda=1}  \\
      &=\left.\frac{N\,(N-1)}{2}\,\int\cdots\int \rchi(x_{1})\,\rchi(x_{2})\,\tilde{\rchi}(x_{3})\cdots \tilde{\rchi}(x_{N})\,P_{\beta}^{(N)}(\vec{x})\,d\,x_{1}\cdots d\,x_{N}\right|_{\lambda=1} \\
      & =\left. \binom{N}{2}\int\cdots\int \rchi(x_{1})\,\rchi(x_{2})\,\tilde{\rchi}(x_{3})\cdots \tilde{\rchi}(x_{N})\,P_{\beta}^{(N)}(\vec{x})\,d\,x_{1}\cdots d\,x_{N}\right|_{\lambda=1} \\
      & =\binom{N}{2}\int\cdots\int \rchi(x_{1})\,\rchi(x_{2})\,\rchi(x_{3})\cdots \rchi(x_{N})\,P_{\beta}^{(N)}(\vec{x})\,d\,x_{1}\cdots d\,x_{N} \\
      & = E_{\beta,N}(t,2), 
    \end{aligned}  
  \end{equation*}
  where we used the previous case to get the first equality, and again the invariance of the integrand under symmetry to get the second equality. By induction then,
  \begin{eqnarray*} 
    \lefteqn{\frac{1}{m\,!}\left(-\frac{d}{d\,\lambda}\right)^{m}\, \left. G_{\beta,N}(t,\lambda)\right|_{\lambda=1} =}\\
    & & = \frac{N\,(N-1)\cdots(N-m+2)}{m\,!}\,\sum_{j=m}^{N}\int\cdots\int \rchi(x_{1})\,\tilde{\rchi}(x_{2})\cdots \\
    & & \qquad \cdots\tilde{\rchi}(x_{j-1})\,\rchi(x_{j})\,\tilde{\rchi}(x_{j+1})\cdots\tilde{\rchi}(x_{N})\,P_{\beta}^{(N)}(\vec{x})\,d\,x_{1}\cdots d\,x_{N}  \\
    & & =\frac{N\,(N-1)\cdots(N-m+1)}{m\,!}\,\int\cdots\int\rchi(x_{1})\cdots \\
    & & \qquad \cdots\rchi(x_{m})\,\tilde{\rchi}(x_{m+1})\cdots\tilde{\rchi}(x_{N})\,P_{\beta}^{(N)}(\vec{x})\,d\,x_{1}\cdots d\,x_{N} \\
    & & = \binom{N}{m}\,\int\cdots\int\rchi(x_{1})\cdots\rchi(x_{m})\,\tilde{\rchi}(x_{m+1})\cdots\tilde{\rchi}(x_{N})\,P_{\beta}^{(N)}(\vec{x})\,d\,x_{1}\cdots d\,x_{N} \\
    & & = E_{\beta,N}(t,m). \end{eqnarray*}
\end{proof}
If we define $F_{\beta,N}(t,m)$ to be the distribution of the $m^{th}$ largest eigenvalue in the (finite $N$) $\beta$--ensemble, then the following probabilistic result is immediate from our definition of $E_{\beta,N}(t,m)$.
\begin{cor}
\begin{equation}F_{\beta,N}(t,m+1) -  F_{\beta,N}(t,m) = E_{\beta,N}(t,m)\end{equation}
\end{cor}

\section{The Distribution of the $m^{th}$ Largest Eigenvalue in the GUE}

\sectionmark{The $m^{th}$ largest eigenvalue in the GUE}

\subsection{The distribution function as a Fredholm determinant} 

We  follow \cite{Trac1} for the derivations that follow. The GUE case corresponds to the specialization $\beta=2$ in \eqref{jointdensity} so that
\begin{equation}
  \label{eq:1}
  G_{2,N}(t,\lambda)=C_{2}^{(N)}\underset{x_{i}\in\mathbb{R}}{\int\cdots\int} \prod_{j<k}\left(x_{j}-x_{k}\right)^{2}\,\prod_{j}^{N}w(x_{j})\,\prod_{j}^{N}\left(1+f(x_{j})\right)\,d\,x_{1}\cdots d\,x_{N}
\end{equation}
where $w(x)=\exp\left(-x^{2}\right)$, $f(x)=-\lambda\,\rchi_{J}(x)$, and $C_{2}^{(N)}$ depends only on $N$. In the steps that follow, additional constants depending solely on $N$ (such as $N!$) which appear will be lumped into $C_{2}^{(N)}$. A probability argument will show that the resulting constant at the end of all calculations simply equals $1$. Expressing the Vandermonde as a determinant
\begin{equation}
  \prod_{1\leq j<k\leq N}(x_{j}-x_{k})=\det\left(x_{k}^{j}\right)_{\substack{j=0,\ldots,N\\k=1,\ldots,N}}
\end{equation}
and using \eqref{debruijn1} with $\varphi_{j}(x)=\psi_{j}(x)=x^{j}$ and $d\,\mu(x)=w(x)(1+f(x))$ yields
\begin{equation}
  G_{2,N}(t,\lambda)=C_{2}^{(N)}\det\left(\int_{\mathbb{R}} x^{j+k}\,w(x)\left(1+f(x)\right)\,d\,x)\right)_{j,k=0,\ldots,N-1}.
\end{equation}
Let $\left\{\varphi_{j}(x)\right\}$ be the sequence obtained by orthonormalizing the sequence $\left\{x^{j}\,w^{1/2}(x)\right\}$. It follows that
\begin{eqnarray}
  G_{2,N}(t,\lambda) & = & C_{2}^{(N)}\det\left(\int_{\mathbb{R}} \varphi_{j}(x)\,\varphi_{k}(x)\,\left(1+f(x)\right)\,d\,x)\right)_{j,k=0,\ldots,N-1}\\
  & = & C_{2}^{(N)}\det\left(\delta_{j,k}+\int_{\mathbb{R}} \varphi_{j}(x)\,\varphi_{k}(x)\,f(x)\,d\,x)\right)_{j,k=0,\ldots,N-1}.
\end{eqnarray}
The last expression is of the form $\det(I+AB)$ for $A:L^{2}(\mathbb{R})\to \mathbb{C}^{N}$ with kernel $A(j,x)=~\varphi_{j}(x)f(x)$ whereas $B:\mathbb{C}^{N}\to L^{2}(\mathbb{R})$ with kernel $B(x,j)=\varphi_{j}(x)$. 
Note that $AB:\mathbb{C}^{N}\to \mathbb{C}^{N}$ has kernel 
\begin{equation}
AB(j,k)=\int_{\mathbb{R}} \varphi_{j}(x)\,\varphi_{k}(x)\,f(x)\,d\,x  
\end{equation}
whereas $BA:L^{2}(\mathbb{R})\to L^{2}(\mathbb{R})$ has kernel
\begin{equation}
BA(x,y)=\sum_{k=0}^{N-1}\varphi_{k}(x)\,\varphi_{k}(y) :=K_{2,N}(x,y).
\end{equation}
From Theorem \ref{detthm} it follows that
\begin{equation}
  G_{2,N}(t,\lambda)=C_{2}^{(N)}\det\left(I+K_{2,N}\,f\right),
\end{equation}
where $K_{2,N}$ has kernel $K_{2,N}(x,y)$ and $K_{2,N}\,f$ acts on a function by first multiplying it by $f$ and acting on the product with $K_{2,N}$. From \eqref{eq:1} we see that setting $f=0$ in the last identity yields $C_{2}^{(N)}=1$. Thus the above simplifies to
\begin{equation}
  \label{eq:22}
  G_{2,N}(t,\lambda)=\det\left(I+K_{2,N}\,f\right).
\end{equation}

\subsection{Edge scaling and differential equations}

We specialize $w(x)=\exp\left(-x^{2}\right)$, $f(x)=-\lambda\,\rchi_{J}(x)$, so that the $\left\{\varphi_{j}(x)\right\}$ are in fact the Hermite polynomials times the square root of the weight. Using the 
Plancherel-Rotach asymptotics of Hermite polynomials, it follows that in the \emph{edge scaling limit},
\begin{equation}
  \label{eq:edgescaling}
\lim_{N\to\infty}\frac{1}{2^{1/2}N^{1/6}}\,K_{N,2}\left(\sqrt{2N}+\frac{x}{2^{1/2}N^{1/6}},\sqrt{2N}+\frac{y}{2^{1/2}N^{1/6}}\right)\,\rchi_{J}\left(\sqrt{2N}+\frac{
y}{2^{1/2}N^{1/6}}\right)
\end{equation}
is $K_{\airy}(x,y)$ as defined in \eqref{kairy}.  As operators,
the convergence is in trace class norm  to $K_{2}$.  (A proof of this last fact
can be found in \cite{Trac5}.) For notational convenience, we denote the corresponding operator $K_{2}$ by $K$ in the rest of this subsection.
It is convenient to view $K$ as the integral operator on $\bR$  with kernel
\begin{equation}
  \label{eq:15}
  K(x,y)=\frac{\varphi(x)\psi(y)-\psi(x)\varphi(y)}{x-y}\,\rchi_{J}(y),
\end{equation}
where $\varphi(x)=\sqrt{\lambda}\airy(x)$, $\psi(x)=\sqrt{\lambda}\airy^{\prime}(x)$ and $J$ is $\left(s,\infty\right)$ with 
\begin{equation}
  t=\sqrt{2N}+ \frac{s}{\sqrt{2}N^{1/6}}.
\end{equation}
Note that although $K(x,y)$, $\varphi$ and $\psi$ are functions of $\lambda$ as well, this dependence will not affect our calculations in what follows. Thus we omit it to avoid cumbersome notation. The Airy equation implies that $\varphi$ and $\psi$ satisfy the relations 
\begin{eqnarray}
  \frac{d}{d\,x}\,\varphi &= &\psi,\nonumber\\
  \frac{d}{d\,x}\,\psi &=& x\,\varphi.
\end{eqnarray}
We define $D_{2,N}(s,\lambda)$ to be the Fredholm determinant $\det\left(I-K_{N,2}\right)$. Thus in the edge scaling limit
\begin{equation*}
  \lim_{N\ra\iy}D_{2,N}(s,\lambda)=D_{2}(s,\lambda).
\end{equation*}
We define the operator
\begin{equation}
  \label{eq:rdef}
  R=(I-K)^{-1}K,
\end{equation}
whose kernel we denote $R(x,y)$. Incidentally, we shall use the notation $\doteq$ in reference to an operator to mean ``has kernel''. For example $R\doteq R(x,y)$. We also let $M$ stand for the operator whose action is multiplication by $x$. It is well known that
\begin{equation}
  \label{rderiv}
  \frac{d}{ds}\log\det\left(I-K\right)=-R(s,s).
\end{equation}
For functions $f$ and $g$, we write $f\otimes g$ to denote the operator specified by
\begin{equation}
  f\otimes g \doteq f(x)g(y),
\end{equation}
and define
\begin{eqnarray}
  Q(x,s)&=& Q(x) = \left(\left(I-K\right)^{-1}\varphi\right)(x),\\
  P(x,s) &=& P(x) = \left(\left(I-K\right)^{-1}\psi\right)(x).
\end{eqnarray}
Then straightforward computation yields the following facts
\begin{eqnarray}
  \com{M}{K} &=& \varphi\otimes\psi - \psi\otimes\varphi,\nonumber\\
  \com{M}{\left(I-K\right)^{-1}} &=& \left(I-K\right)^{-1}\com{M}{K}\left(I-K\right)^{-1} \nonumber\\
  & = & Q\otimes P-P\otimes Q.
\end{eqnarray}
On the other hand if $\left(I-K\right)^{-1}\doteq\rho(x,y)$, then 
\begin{equation}
  \rho(x,y)=\delta(x-y)+R(x,y),
\end{equation}
and it follows that 
\begin{equation}
\com{M}{\left(I-K\right)^{-1}} \doteq (x-y)\rho(x,y)=(x-y)R(x,y).
\end{equation}
Equating the two representation for the kernel of $\com{M}{\left(I-K\right)^{-1}}$ yields
\begin{equation}
  R(x,y)=\frac{Q(x)P(y)-P(x)Q(y)}{x-y}.
\end{equation}
Taking the limit $y\to x$ and defining $q(s)=Q(s,s)$, $p(s)=P(s,s)$, we obtain
\begin{equation}
  \label{RDiag}
  R(s,s)=Q^{\prime}(s,s)\,p(s)-P^{\prime}(s,s)\,q(s).
\end{equation}
Let us now derive expressions for $Q^{\prime}(x)$ and $P^\prime(x)$. If we let the operator $D$ stand for differentiation with respect to $x$,
\begin{eqnarray}
  Q^\prime(x,s)&=& D \left(I-K\right)^{-1} \varphi \nonumber\\
  &=& \left(I-K\right)^{-1} D\varphi +
  \left[D,\left(I-K\right)^{-1}\right]\varphi\nonumber\\
  &=& \left(I-K\right)^{-1} \psi +
  \left[D,\left(I-K\right)^{-1}\right]\varphi\nonumber\\
  &=& P(x) + \left[D,\left(I-K\right)^{-1}\right]\varphi. \label{Qderiv1}
\end{eqnarray}
We need the commutator
\begin{equation} 
  \left[D,\left(I-K\right)^{-1}\right]=\left(I-K\right)^{-1} \left[D,K\right] \left(I-K\right)^{-1}. 
\end{equation}
Integration by parts shows
\begin{equation}
  \left[D,K\right] \doteq \left( \frac{\partial K}{\partial x} + \frac{\partial K}{\partial y}\right) + K(x,s) \delta(y-s). 
\end{equation} 
The $\delta$ function comes from differentiating the characteristic function $\chi$. Moreover,
\begin{equation}
  \left( \frac{\partial K}{\partial x} + \frac{\partial K}{\partial y}\right) = \varphi(x) \varphi(y).
\end{equation} 
Thus
\begin{equation} 
 \label{DComm}
\left[D,\left(I-K\right)^{-1}\right]\doteq - Q(x) Q(y) + R(x,s) \rho(s,y).
\end{equation}
(Recall $(I-K)^{-1}\doteq \rho(x,y)$.)  We now use this in (\ref{Qderiv1}) to obtain
\begin{eqnarray*}
  Q^\prime(x,s)&=&P(x) - Q(x) \left(Q,\varphi\right) + R(x,s) q(s) \\
  &=& P(x) - Q(x) u(s) + R(x,s) q(s),
\end{eqnarray*}
where the inner product $\left(Q,\varphi\right)$ is denoted by $u(s)$. Evaluating  at $x=s$  gives
\begin{equation}
  \label{Qderiv2}
  Q^\prime(s,s) = p(s) - q(s) u(s) +R(s,s) q(s).
\end{equation}
We now apply the same procedure to compute $P^\prime$. 
\begin{eqnarray*}
  P^\prime(x,s)&=& D \left(I-K\right)^{-1} \psi \\
  &=& \left(I-K\right)^{-1} D\psi + \left[D,\left(I-K\right)^{-1}\right]\psi\\
  &=& M \left(I-K\right)^{-1} \varphi +  \left[\left(I-K\right)^{-1},M\right]\varphi+   \left[D,\left(I-K\right)^{-1}\right]\psi\\
  &=& x Q(x) +\left(P\otimes Q-Q\otimes P\right)\varphi +(-Q\otimes Q)\psi +  R(x,s) p(s)\\
  &=& x Q(x) + P(x)\left(Q,\varphi\right) -  Q(x) \left(P,\varphi\right)  - Q(x) \left(Q,\psi\right)+R(x,s)p(s)\\
  &=& x Q(x) - 2 Q(x) v(s) + P(x) u(s) + R(x,s) p(s).
\end{eqnarray*}
Here $v=\left(P,\varphi\right)=\left(\psi,Q\right)$. Setting $x=s$ we obtain
\begin{equation}
  \label{Pderiv}
  P^{\prime}(s,s) = s q(s) + 2 q(s) v(s) +p(s) u(s) +R(s,s) p(s).
\end{equation} 
Using this and the expression for $Q^\prime(s,s)$ in (\ref{RDiag}) gives
\begin{equation}
  \label{RDiag2}
  R(s,s)= p^2-s q^2 + 2 q^2 v - 2 p q u.
\end{equation}
Using the chain rule, we have
\begin{equation} 
  \label{qDeriv}
  \frac{d\,q}{d\,s} = \left( \frac{\partial}{\partial x} + \frac{\partial}{\partial s}\right) Q(x,s)\left\vert_{x=s}. \right.
\end{equation}
The first term is known. The partial with respect to $s$ is
\begin{eqnarray*}
  \frac{\partial Q(x,s)}{\partial s}&=& \left(I-K\right)^{-1} \frac{\partial K}{\partial s} \left(I-K\right)^{-1} \varphi\\
  &=& - R(x,s) q(s),
\end{eqnarray*}
where we used the fact that
\begin{equation}
  \frac{\partial K}{\partial s}\doteq -K(x,s)\delta(y-s).
\end{equation}
Adding the two partial derivatives  and evaluating at $x=s$ gives
\begin{equation}
  \label{qEqn}
  \frac{d\,q}{d\,s} = p - q u.
\end{equation}
A similar calculation gives
\begin{equation}
  \label{pEqn}
  \frac{d\,p}{d\,s}= s q - 2 q v + p u.
\end{equation}
We derive first order differential equations for  $u$ and $v$  by differentiating the inner products. Recall that
\begin{equation*}
  u(s) = \int_s^\infty \varphi(x) Q(x,s)\, d\,x.
\end{equation*}
Thus
\begin{eqnarray*}
  \frac{d\,u}{ds}&=& -\varphi(s) q(s) + \int_s^\infty \varphi(x) \frac{\partial Q(x,s)}{\partial s}\, d\,x \\
  &=& -\left(\varphi(s)+\int_s^\infty R(s,x) \varphi(x)\,d\,x\right) q(s)\\
  &=& -\left(I-K\right)^{-1} \varphi(s) \, q(s)\\
  &=& - q^2.
\end{eqnarray*}
Similarly,
\begin{equation}
  \frac{d\,v}{d\,s} = - p q.
\end{equation}
From the first order differential equations for $q$, $u$ and $v$ it follows immediately  that the derivative of   $ u^2-2v-q^2 $ is zero.  Examining the behavior near $s=\infty$ to check that the constant of integration is zero then gives
\begin{equation}
  u^2-2v=q^2.
\end{equation}
We now differentiate (\ref{qEqn}) with respect to $s$, use the first order differential equations for $p$ and $u$, and then the  first integral to deduce that $q$ satisfies the Painlev\'e II equation (\ref{pII}).
Checking the asymptotics of the Fredholm determinant $\det(I-K)$ for large $s$ shows we want the solution  with boundary condition
\begin{equation}
  \label{bc}
  q(s,\lambda)\sim \sqrt{\lambda}\airy(s) \qquad \textrm{as} \qquad s\rightarrow\infty.
\end{equation}
That a solution $q$ exists and is unique follows from the representation of the Fredholm determinant in terms of it.  Independent proofs of this, as well as the asymptotics as $s\ra\-\iy$ were given by \cite{Hasti1}, \cite{Clar1}, \cite{Deif2}. Since $\com{D}{(I-K)^{-1}}\doteq (\partial/\partial x~+~\partial/\partial y)~R(x,y)$, (\ref{DComm}) says
\begin{equation}
  \left(\frac{\partial}{\partial x}+\frac{\partial}{\partial y}\right)R(x,y)=-Q(x)Q(y)+R(x,s)\rho(s,y).
\end{equation}
In computing $\partial Q(x,s)/\partial s$ we showed that
\begin{equation}
  \frac{\partial}{\partial s} \left(I-K\right)^{-1}\doteq \frac{\partial}{\partial s}R(x,y) = -R(x,s)\rho(s,y).
\end{equation}
Adding these two expressions,
\begin{equation}
  \left(\frac{\partial}{\partial x} + \frac{\partial}{\partial y}+ \frac{\partial}{\partial s}\right)R(x,y) = -Q(x)\,Q(y),
\end{equation}
and then evaluating at $x=y=s$ gives
\begin{equation}
  \frac{d}{d\,s}R(s,s)=-q^2. \label{Rderiv}
\end{equation}
Integration  (and recalling (\ref{rderiv})) gives,
\begin{equation}
  \frac{d}{d\,s}\log\det\left(I-K\right)=-\int_s^\infty q^2(x,\lambda) \, d\,x;
\end{equation}
and hence,
\begin{equation}
  \log\det\left(I-K\right)=-\int_s^\infty\left(\int_y^\infty q^2(x,\lambda)\,d\,x\right)\, d\,y =-\int_s^\infty (x-s) q^2(x,\lambda)\, d\,x.
\end{equation}
To summarize,  we have shown that $D_{2}(s,\lambda)$ has the Painlev\'e representation (\ref{D2})
where $q$  satisfies the Painlev\'e II equation (\ref{pII}) subject to the boundary condition (\ref{bc}).

\section{The Distribution of the $m^{th}$ Largest Eigenvalue in the GSE}

\sectionmark{The $m^{th}$ largest eigenvalue in the GSE}

\subsection{The distribution function as a Fredholm determinant} 

The GSE corresponds case corresponds to the specialization $\beta=4$ in \eqref{jointdensity} so that
\begin{equation}
  G_{4,N}(t,\lambda)=C_{4}^{(N)}\underset{x_{i}\in\mathbb{R}}{\int\cdots\int} \prod_{j<k}\left(x_{j}-x_{k}\right)^{4}\,\prod_{j}^{N}w(x_{j})\,\prod_{j}^{N}\left(1+f(x_{j})\right)\,d\,x_{1}\cdots d\,x_{N}
\end{equation}
where $w(x)=\exp\left(-x^{2}\right)$, $f(x)=-\lambda\,\rchi_{J}(x)$, and $C_{4}^{(N)}$ depends only on $N$. As in the GUE case, we will
absorb  into $C_{4}^{(N)}$ any constants depending only on $N$ that appear in the derivation. A simple argument  at the end will show that the final constant is $1$. These calculations  follow  \cite{Trac1}. By \eqref{vandthm}, $G_{4,N}(t,\lambda)$ is given by the integral
\begin{equation*}
  C_{4}^{(N)}\underset{x_{i}\in \mathbb R}{\int\cdots\int} \det\left(x_{k}^{j}\quad j\,x_{k}^{j-1}\right)_{\substack{j=0,\dots,2\,N-1 \\ k=1,\ldots,N}}\,\prod_{i=1}^{N}w(x_{i})\prod_{i=1}^{N}(1+f(x_{i}))\,d\,x_{1}\cdots d\,x_{N}
\end{equation*}
which, if we define $\varphi_{j}(x)=x^{j-1}\,w(x)\,(1+f(x))$ and $\psi_{j}(x)=(j-1)\,x^{j-2}$ and use the linearity of the determinant, becomes
\begin{equation*}
  G_{4,N}(t,\lambda)= C_{4}^{(N)}\underset{x_{i}\in \mathbb R}{\int\cdots\int} \det\left(\varphi_{j}(x_{k})\quad \psi(x_{k})\right)_{\substack{1\leq j \leq 2\,N \\
      1\leq k \leq N}}\,d\,x_{1}\cdots d\,x_{N}.
\end{equation*}
Now using \eqref{debruijn3}, we obtain
\begin{eqnarray*}
  G_{4,N}(t,\lambda) & = & C_{4}^{(N)}\,\pfaffian\left(\int \varphi_{j}(x)\psi_{k}{x}-\varphi_{k}(x)\psi_{j}(x)\,d\,x\right)_{_{1\leq j,k\leq 2\,N}}\\
  &=& C_{4}^{(N)}\,\pfaffian\left(\int (k-j)\,x^{j+k-3}\,w(x)\,(1+f(x))\,d\,x\right)_{_{1\leq j,k\leq 2\,N}} \\
  &=&  C_{4}^{(N)}\pfaffian\left(\int (k-j)\,x^{j+k-1}\,w(x)\,(1+f(x))\,d\,x\right)_{_{0\leq j,k\leq 2\,N-1}},
\end{eqnarray*}
where we let $k\to k+1$ and $j\to j+1$ in the last line. Remembering that the square of a Pfaffian is a determinant, we
obtain
\begin{equation*}
  G_{4,N}^{2}(t,\lambda)=C_{4}^{(N)}\, \det\left(\int (k-j)\,x^{j+k-1}\,w(x)\,(1+f(x))\,d\,x\right)_{_{0\leq j,k\leq 2\,N-1}}.
\end{equation*}
Row operations on the matrix do not change the determinant, so we can replace $\{x^{j}\}$ by an arbitrary sequence $\{p_{j}(x)\}$ of polynomials of degree $j$ obtained by adding rows to each other. Note that the general $(j,k)$ element in the matrix can be written as
\begin{equation*}
  \left[\left(\frac{d}{d\,x}\,x^{k}\right)\,x^{j} - \left(\frac{d}{d\,x}\,x^{j}\right)\,x^{k} \right]\,w(x) \,\left(1+f(x)\right).
\end{equation*}
Thus when we add rows to each other the polynomials we obtain will have the same general form (the derivatives factor). Therefore we can assume without loss of generality that $G_{4,N}^{2}(t,\lambda)$ equals
\begin{equation*}
  C_{4}^{(N)}\, \det\left(\int \left[ p_{j}(x)\,p_{k}^{\prime}(x) - \,p_{j}^{\prime}(x)\,p_{k}(x)\right]\,w(x)\,\left(1+f(x)\right)\,d\,x\right)_{_{0\leq j,k\leq 2\,N-1}},
\end{equation*}
where the sequence $\{p_{j}(x)\}$ of polynomials of degree $j$ is arbitrary. Let $\psi_{j}=p_{j}\,w^{1/2}$ so that $p_{j}=\psi_{j}\,w^{-1/2}$. Substituting this into the above formula and simplifying, we obtain
\begin{eqnarray*}
  G_{4,N}^{2}(t,\lambda)&=&C_{4}^{(N)}\, \det\left[ \int\left( \left( \psi_{j}(x)\,\psi_{k}^{\prime}(x) - \,\psi_{k}(x)\,\psi_{j}^{\prime}(x)\right)\left(1+f(x)\right)\right)\,d\,x \right]_{_{0\leq j,k\leq 2\,N-1}}\\
  &=&C_{4}^{(N)}\, \det\left[M + L\right] =C_{4}^{(N)}\, \det[M]\,\det[I+ M^{-1}\cdot L],
\end{eqnarray*}
where $M,L$ are matrices given by
\begin{equation*}
  M=\left(\int \left( \psi_{j}(x)\,\psi_{k}^{\prime}(x) - \,\psi_{k}(x)\,\psi_{j}^{\prime}(x)\right)\,d\,x\right)_{_{0\leq j,k\leq 2\,N-1}},
\end{equation*}
\begin{equation*}
  L=\left(\int \left( \psi_{j}(x)\,\psi_{k}^{\prime}(x) - \,\psi_{k}(x)\,\psi_{j}^{\prime}(x)\right)f(x)\,d\,x\right)_{_{0\leq j,k\leq 2\,N-1}}.
\end{equation*}
Note that $\det[M]$ is a constant which depends only on $N$ so we can absorb it into $C_{4}^{(N)}$. Also if we denote
\begin{equation*}
  M^{-1}=\left\{\mu_{jk}\right\}_{_{0\leq j,k\leq 2\,N-1}}\quad , \quad \eta_{j}=\sum_{k=0}^{2N-1}\mu_{jk}\psi_{k}(x),
\end{equation*} 
it follows that
\begin{equation*}
  M^{-1}\cdot N=\left\{\int \left( \eta_{j}(x)\,\psi_{k}^{\prime}(x) - \,\eta_{j}^{\prime}(x)\,\psi_{k}(x)\right)f(x)\,d\,x\right\}_{_{0\leq j,k\leq 2\,N-1}}.
\end{equation*}
Let $A:L^{2}(\mathbb R)\times L^{2}(\mathbb R)\to \mathbb C^{2N}$ be the operator defined by the $2N\times 2$ matrix
\begin{equation*}
  A(x)=\left(
    \begin{array}{cc} 
      \eta_{0}(x) & -\eta_{0}^{\prime}(x) \\
      \eta_{1}(x) & -\eta_{1}^{\prime}(x) \\
      \vdots & \vdots  \\
    \end{array}
  \right).
\end{equation*}
Thus if
\begin{equation*}
  g=
  \left(
    \begin{array}{c}
      g_{0}(x)  \\
      g_{1}(x) 
    \end{array}
  \right)
  \in L^{2}(\mathbb R)\times L^{2}(\mathbb R),
\end{equation*}
we have
\begin{equation*}
  A\,g=A(x)\,g= 
  \left(
    \begin{array}{c}
      \int\left(\eta_{0}g_{0}-\eta_{0}^{\prime}g_{1}\right)\,d\,x \\
      \int\left(\eta_{1}g_{0}-\eta_{1}^{\prime}g_{1}\right)\,d\,x \\
      \vdots  \\
    \end{array}
  \right)
  \in \mathbb C^{2N}.
\end{equation*}
Similarly we define $B:\mathbb C^{2n}\to L^{2}(\mathbb R)\times L^{2}(\mathbb R)$ given by the $2\times 2n$ matrix
\begin{equation*}
  B(x) = f\cdot
  \left(
    \begin{array}{ccc}
      \psi_{0}^{\prime}(x) & \psi_{1}^{\prime}(x) & \cdots \\
      \psi_{0}(x) & \psi_{1}(x) & \cdots
    \end{array}
  \right).
\end{equation*}
Explicitly if
\begin{equation*}
  \alpha=
  \left(
    \begin{array}{c}
      \alpha_{0}  \\
      \alpha_{1} \\
      \vdots\end{array}
  \right)
  \in \mathbb C^{2N},
\end{equation*}
then
\begin{equation*}
  B\alpha= B(x)\cdot\alpha= 
  \left(
    \begin{array}{c}
      f\displaystyle{\sum_{i=0}^{2N-1}\alpha_{i}\psi_{i}^{\prime}}\\
      f\displaystyle{\sum_{i=0}^{2N-1}\alpha_{i}\psi_{i}}
    \end{array}
  \right)\in L^{2}\times L^{2}.
\end{equation*}
Observe that $M^{-1}\cdot L=A B: \mathbb C^{2n}\to \mathbb C^{2n}$. Indeed
\begin{eqnarray*}
  A B \alpha &=& 
  \left(
    \begin{array}{c} 
      \displaystyle{\sum_{i=0}^{2N-1}\left[\int\left( \eta_{0}\psi_{i}^{\prime}-\eta_{0}^{\prime}\psi_{i}\right)f\, d\,x\right] \alpha_{i}} \\
      \displaystyle{\sum_{i=0}^{2N-1}\left[\int\left( \eta_{1}\psi_{i}^{\prime}-\eta_{1}^{\prime}\psi_{i}\right)f\, d\,x\right] \alpha_{i}} \\
      \vdots
    \end{array}
  \right) \\
  & = & \left\{\int\left(\eta_{j}\psi_{k}^{\prime}-\eta_{j}^{\prime}\psi_{k}\right)f\, d\,x \right\}\cdot 
  \left(
    \begin{array}{c}
      \alpha_{0}  \\
      \alpha_{1} \\ 
      \vdots
    \end{array}
  \right) 
  =\left( M^{-1}\cdot L\right) \alpha.
\end{eqnarray*}
Therefore, by \eqref{detthm}
\begin{equation*}G_{4,N}^{2}(t,\lambda) = C_{4}^{(N)}\det(I+M^{-1}\cdot L) =C_{4}^{(N)}\det(I+AB)=C_{4}^{(N)}\det(I+BA)\end{equation*}
where $BA:L^{2}\left(\mathbb R\right)\to L^{2}\left(\mathbb R\right)$. From our definition of $A$ and $B$ it follows that
\begin{eqnarray*}
  B\,A\,g & = & 
  \left(
    \begin{array}{c}
      f\displaystyle{\sum_{i=0}^{2n-1}\psi_{i}^{\prime} \left( \int\left(\eta_{i}g_{0}-\eta_{i}^{\prime}g_{1}\right)\,d\,x \right)} \\
      f \displaystyle{\sum_{i=0}^{2N-1}\psi_{i}^{\prime} \left(\int\left(\eta_{i}g_{0}-\eta_{i}^{\prime}g_{1}\right)\,d\,x \right)}
    \end{array}\right) \\
  & =&  f
  \left(
    \begin{array}{c}
      \displaystyle{\int \sum_{i=0}^{2N-1}\psi_{i}^{\prime}(x)\eta_{i}(y)g_{0}(y)\,d\,y - \int \sum_{i=0}^{2N-1}\psi_{i}^{\prime}(x)\eta_{i}^{\prime}(y)g_{1}(y)\,d\,y } \\
      \displaystyle{\int\sum_{i=0}^{2N-1}\psi_{i}(x)\eta_{i}(y)g_{0}(y)\,d\,y - \int\sum_{i=0}^{2N-1}\psi_{i}(x)\eta_{i}^{\prime}(y)g_{1}(y)\,d\,y}
    \end{array}
  \right)
  \\ 
  & = &f\, K_{4,N}\,g,
\end{eqnarray*} 
where $K_{4,N}$ is the integral operator with matrix kernel
\begin{equation*}
  K_{4,N}(x,y)=
  \left(
    \begin{array}{cc}
      \displaystyle{ \sum_{i=0}^{2N-1}\psi_{i}^{\prime}(x)\eta_{i}(y)} & \displaystyle{-\sum_{i=0}^{2N-1}\psi_{i}^{\prime}(x)\eta_{i}^{\prime}(y)} \\
      \displaystyle{ \sum_{i=0}^{2N-1}\psi_{i}(x)\eta_{i}(y)} & \displaystyle{- \sum_{i=0}^{2N-1}\psi_{i}(x)\eta_{i}^{\prime}(y)}
    \end{array}
  \right).
\end{equation*}
Recall that
$\displaystyle{\eta_{j}(x)=\sum_{k=0}^{2N-1}\mu_{jk}\psi_{k}(x)}$ so that
\begin{equation*}
  K_{4,N}(x,y)=
  \left(
    \begin{array}{cc}
      \displaystyle{ \sum_{j,k=0}^{2N-1}\psi_{j}^{\prime}(x)\mu_{jk}\psi_{k}(y)} & \displaystyle{-\sum_{j,k=0}^{2N-1}\psi_{j}^{\prime}(x)\mu_{jk}\psi_{k}^{\prime}(y)} \\
      \displaystyle{ \sum_{j,k=0}^{2N-1}\psi_{j}(x)\mu_{jk}\psi_{k}(y)}& \displaystyle{-\sum_{j,k=0}^{2N-1}\psi_{j}(x)\mu_{jk}\psi_{k}^{\prime}(y)}
    \end{array}
  \right).
\end{equation*}
Define $\epsilon$ to be the following integral operator
\begin{equation}
  \label{epsilonop}
\epsilon\doteq  \epsilon(x-y)=
  \begin{cases}
   \>\>\; \frac{1}{2} &\textrm{if $x>y$}, \\
    -\frac{1}{2} &\textrm{if $x<y$}.
  \end{cases}
\end{equation}
As before, let $D$ denote the operator that acts by differentiation with respect to $x$. The fundamental theorem of calculus implies that $D\,\epsilon~=~\epsilon\,D~=~I$. We also define
\begin{equation*}
  S_{N}(x,y)=\sum_{j,k=0}^{2N-1}\psi_{j}^{\prime}(x)\mu_{jk}\psi_{k}(y).
\end{equation*}
Since $M$ is antisymmetric,
\begin{equation*}
  S_{N}(y,x)=\sum_{j,k=0}^{2N-1}\psi_{j}^{\prime}(y)\mu_{jk}\psi_{k}(x)=-\sum_{j,k=0}^{2N-1}\psi_{j}^{\prime}(y)\mu_{kj}\psi_{k}(x)= -\sum_{j,k=0}^{2N-1}\psi_{j}(y)\mu_{kj}\psi_{k}^{\prime}(x),
\end{equation*}
after re-indexing. Note that
\begin{equation*}
  \epsilon\,S_{N}(x,y)=\sum_{j,k=0}^{2N-1}\epsilon\,D\,\psi_{j}(x)\mu_{jk}\psi_{k}(y) =\sum_{j,k=0}^{2N-1}\psi_{j}(x)\mu_{jk}\psi_{k}(y), 
\end{equation*}
and
\begin{equation*}
  -\frac{d}{d\,y}\,S_{N}(x,y)=\sum_{j,k=0}^{2N-1}\psi_{j}^{\prime}(x)\mu_{jk}\psi_{k}^{\prime}(y).
\end{equation*}
Thus we can now write succinctly
\begin{equation}
  \label{gsekernel}
  K_{N}(x,y)=
  \left(
    \begin{array}{cc}
      S_{N}(x,y) & -\frac{d}{d\,y}\,S_{N}(x,y) \\
      \epsilon\,S_{N}(x,y) & S_{N}(y,x)
    \end{array}
  \right).
\end{equation}
To summarize, we have shown that $G_{4,N}^{2}(t,\lambda)=C_{4}^{(N)}\det(I-K_{4,N}f)$. Setting $f\equiv 0$ on both sides (where the original definition of $G_{4,N}(t,\lambda)$ as an integral is used on the left) shows that $C_{4}^{(N)}=1$. Thus
\begin{equation}
  \label{eq:5}
  G_{4,N}(t,\lambda)=\sqrt{D_{4,N}(t,\lambda)},
\end{equation}
where we define
\begin{equation}
  D_{4,N}(t,\lambda)=\det(I+K_{4,N}f),
\end{equation}
and $K_{4,N}$ is the integral operator with matrix kernel \eqref{gsekernel}. 

\subsection{Gaussian specialization}

We would like to specialize the above results to the case of a Gaussian weight function
\begin{equation}
  \label{gseweight}
  w(x)=\exp\left(x^{2}\right)
\end{equation}
and indicator function
\begin{equation*}
  f(x)=-\lambda\,\rchi_{_{J}}\quad , \quad J=(t,\infty).
\end{equation*}
We want the matrix
\begin{equation*}
  M=\left\{\int \left(\psi_{j}(x)\,\psi_{k}^{\prime}(x) - \,\psi_{k}(x)\,\psi_{j}^{\prime}(x)\right)\,d\,x\right\}_{_{0\leq j,k\leq 2\,N-1}}
\end{equation*}
to be the direct sum of $N$ copies of
\begin{equation*}
  Z=\left(
    \begin{array}{cc}
      0 & 1 \\
      -1 & 0
    \end{array}
  \right),
\end{equation*}
so that the formulas are the simplest possible, since then $\mu_{jk}$ can only be $0$ or $\pm 1$.  In that case $M$ would be skew--symmetric so that $M^{-1}=-M$. In terms of the integrals defining the entries of $M$ this means that we would like to have
\begin{equation*}
  \int \left( \psi_{2j}(x)\,\frac{d}{d\,x}\,\psi_{2k+1}(x) - \,\psi_{2k+1}(x)\,\frac{d}{d\,x}\,\psi_{2j}(x)\right)\,d\,x = \delta_{j,k},
\end{equation*}
\begin{equation*}
  \int \left( \psi_{2j+1}(x)\,\frac{d}{d\,x}\,\psi_{2k}(x) - \,\psi_{2k}(x)\,\frac{d}{d\,x}\,\psi_{2j+1}(x)\right)\,d\,x = -\delta_{j,k}
\end{equation*}
and otherwise
\begin{equation*}
  \int \left( \psi_{j}(x)\,\frac{d}{d\,x}\,\psi_{k}(x) - \,\psi_{j}(x)\,\frac{d}{d\,x}\,\psi_{k}(x)\right)\,d\,x = 0.
\end{equation*}
It is easier to treat this last case if we replace it with three non-exclusive conditions
\begin{equation*}
  \int \left(\psi_{2j}(x)\,\frac{d}{d\,x}\,\psi_{2k}(x) -\,\psi_{2k}(x)\,\frac{d}{d\,x}\,\psi_{2j}(x)\right)\,d\,x = 0,
\end{equation*}
\begin{equation*}
  \int \left( \psi_{2j+1}(x)\,\frac{d}{d\,x}\,\psi_{2k+1}(x) - \,\psi_{2k+1}(x)\,\frac{d}{d\,x}\,\psi_{2j+1}(x)\right)\,d\,x = 0,
\end{equation*}
(so when the parity is the same for  $j,k$, which takes care of diagonal entries, among others) and
\begin{equation*}
  \int \left( \psi_{j}(x)\,\frac{d}{d\,x}\,\psi_{k}(x) - \,\psi_{j}(x)\,\frac{d}{d\,x}\,\psi_{k}(x)\right)\,d\,x = 0,
\end{equation*}
whenever $|j-k|>1$, which targets entries outside of the tridiagonal.
Define
\begin{equation}
  \label{gsevarphi}
  \varphi_{k}(x)=\frac{1}{c_{k}}H_{k}(x)\,\exp(-x^{2}/2)\quad\textrm{for}\quad c_{k}=\sqrt{2^{k}k!\sqrt{\pi}}
\end{equation}
where the $H_{k}$ are the usual Hermite polynomials defined by the orthogonality condition
\begin{equation*}
  \int_{\mathbb{R}} H_{j}(x)\,H_{k}(x)\,e^{-x^{2}}\,d\,x = c_{j}^{2}\,\delta_{j,k}.
\end{equation*}
Then it follows that
\begin{equation*}
  \int_{\mathbb{R}} \varphi_{j}(x)\,\varphi_{k}(x)\,d\,x = \delta_{j,k}.
\end{equation*}
Now let
\begin{equation*}
  \psi_{_{2j+1}}(x) = \frac{1}{\sqrt{2}}\,\varphi_{_{2j+1}}\left(x\right)\qquad \psi_{_{2j}}(x) = -\frac{1}{\sqrt{2}}\,\epsilon\,\varphi_{_{2j+1}}\left(x\right)
\end{equation*}
This definition satisfies our earlier requirement that $\psi_{j}=p_{j}\,w^{1/2}$ with $w$ defined in \eqref{gseweight}. In particular we have in this case
\begin{equation*}
  p_{2j+1}(x)=\frac{1}{c_{j}\sqrt{2}}\,H_{2j+1}\left(x\right).
\end{equation*}
Let $\epsilon$ as in \eqref{epsilonop}, and $D$ denote the operator that acts by differentiation with respect to $x$ as before, so that $D\,\epsilon=\epsilon\,D=I$. It follows that
\begin{equation*}
  \begin{aligned}
    \int_{\mathbb{R}} & \left[\psi_{_{2j}}(x)\frac{d}{d\,x}\,\psi_{_{2k+1}}(x)-\psi_{_{2k+1}}(x)\frac{d}{d\,x}\,\psi_{_{2j}}(x)\right]d\,x\\
    & = \frac{1}{2}\int_{\mathbb{R}}\left[-\epsilon\,\varphi_{_{2j+1}}\left(x\right)\frac{d}{d\,x}\,\varphi_{_{2k+1}}\left(x\right) + \varphi_{_{2k+1}}\left(x\right)\frac{d}{d\,x}\,\epsilon\,\varphi_{_{2j+1}}\left(x\right)\right]d\,x\\
    & = \frac{1}{2}\int_{\mathbb{R}}\left[ -\epsilon\,\varphi_{_{2j+1}}\left(x\right)\frac{d}{d\,x}\,\varphi_{_{2k+1}}\left(x\right) +\varphi_{_{2k+1}}\left(x\right)\varphi_{_{2j+1}}\left(x\right)\right]d\,x\\
  \end{aligned}
\end{equation*}
We integrate the first term by parts and use the fact that
\begin{equation*}
  \frac{d}{d\,x}\,\epsilon\,\varphi_{_{j}}\left(x\right) =\varphi_{_{j}}\left(x\right)
\end{equation*}
and also that $\varphi_{_{j}}$ vanishes at the boundary (i.e. $\varphi_{_{j}}\left(\pm\infty\right) =0$) to obtain
\begin{equation*}
  \begin{aligned} 
    \int_{\mathbb{R}} & \left[\psi_{_{2j}}(x)\frac{d}{d\,x}\,\psi_{_{2k+1}}(x)-\psi_{_{2k+1}}(x)\frac{d}{d\,x}\,\psi_{_{2j}}(x)\right]d\,x\\
    & = \frac{1}{2} \int_{\mathbb{R}}\left[-\epsilon\,\varphi_{_{2j+1}}\left(x\right)\frac{d}{d\,x}\,\varphi_{_{2k+1}}\left(x\right) +\varphi_{_{2k+1}}\left(x\right)\varphi_{_{2j+1}}\left(x\right)\right]d\,x\\
    & = \frac{1}{2}  \int_{\mathbb{R}}\left[ \varphi_{_{2j+1}}\left(x\right)\,\varphi_{_{2k+1}}\left(x\right)+\varphi_{_{2k+1}}\left(x\right)\varphi_{_{2j+1}}\left(x\right)\right]d\,x\\
    & = \frac{1}{2}  \int_{\mathbb{R}}\left[ \varphi_{_{2j+1}}\left(x\right)\,\varphi_{_{2k+1}}\left(x\right)+\varphi_{_{2k+1}}\left(x\right)\varphi_{_{2j+1}}\left(x\right)\right]d\,x\\
    &= \frac{1}{2} \left(\delta_{j,k}+ \delta_{j,k}\right)\\
    & =\delta_{j,k},
  \end{aligned}
\end{equation*}
as desired. Similarly
\begin{equation*}
  \begin{aligned}
    \int_{\mathbb{R}} & \left[\psi_{_{2j+1}}(x)\frac{d}{d\,x}\,\psi_{_{2k}}(x)-\psi_{_{2k}}(x)\frac{d}{d\,x}\,\psi_{_{2j+1}}(x)\right]d\,x\\
    & = \frac{1}{2}\int_{\mathbb{R}}\left[-\varphi_{_{2j+1}}\left(x\right)\frac{d}{d\,x}\,\epsilon\,\varphi_{_{2k+1}}\left(x\right) +\epsilon\,\varphi_{_{2k+1}}\left(x\right)\frac{d}{d\,x}\,\varphi_{_{2j+1}}\left(x\right)\right]d\,x \\
    & = -\delta_{j,k}.
  \end{aligned}
\end{equation*}
Moreover,
\begin{equation*}
  p_{2j+1}(x)=\frac{1}{c_{j}\sqrt{2}}\,H_{2j+1}\left(x\right)
\end{equation*}
is certainly an odd function, being the multiple of and odd Hermite polynomial. On the other hand, one easily checks that $\epsilon$ maps odd functions to even functions on $L^{2}(\mathbb{R})$. Therefore 
\begin{equation*}
  p_{2j}(x)=-\frac{1}{c_{j}\sqrt{2}}\,\epsilon\,H_{2j+1}\left(x\right)
\end{equation*}
is an even function, and it follows that
\begin{equation*}
  \begin{aligned}\int_{\mathbb{R}} & \left[ \psi_{_{2k}}(x)\frac{d}{d\,x}\,\psi_{_{2j}}(x)-\psi_{_{2j}}(x)\frac{d}{d\,x}\,\psi_{_{2k}}(x)\right]d\,x\\
    & = \int_{\mathbb{R}}\left[p_{2j}(x)\frac{d}{d\,x}\,p_{2k}(x)-p_{2k}(x)\frac{d}{d\,x}\,p_{2j}(x)\right]w(x)\,d\,x\\
    & = 0,
  \end{aligned}
\end{equation*}
since both terms in the integrand are odd functions, and the weight function is even. Similarly,
\begin{equation*}
  \begin{aligned}
    \int_{\mathbb{R}} & \left[\psi_{_{2k+1}}(x)\frac{d}{d\,x}\,\psi_{_{2j+1}}(x)-\psi_{_{2j+1}}(x)\frac{d}{d\,x}\,\psi_{_{2k+1}}(x)\right]d\,x\\
    & = \int_{\mathbb{R}}\left[p_{2j+1}(x)\frac{d}{d\,x}\,p_{2k+1}(x)-p_{2k+1}(x)\frac{d}{d\,x}\,p_{2j+1}(x)\right]w(x)\,d\,x\\
    & = 0.
  \end{aligned}
\end{equation*}
Finally it is easy to see that if $|j-k|>1$ then
\begin{equation*}
  \int_{\mathbb{R}} \left[\psi_{_{j}}(x)\frac{d}{d\,x}\,\psi_{_{k}}(x)-\psi_{_{j}}(x)\frac{d}{d\,x}\,\psi_{_{k}}(x)\right]d\,x=0.
\end{equation*}
Indeed both  differentiation and the action of $\epsilon$ can only ``shift'' the indices by $1$. Thus by orthogonality of the $\varphi_{j}$, this integral will always be $0$. Hence by choosing
\begin{equation*}
  \psi_{_{2j+1}}(x) = \frac{1}{\sqrt{2}}\,\varphi_{_{2j+1}}\left(x\right),\qquad \psi_{_{2j}}(x) = -\frac{1}{\sqrt{2}}\,\epsilon\,\varphi_{_{2j+1}}\left(x\right),
\end{equation*}
we force the matrix
\begin{equation*}
  M=\left\{\int_{\mathbb{R}} \left(\psi_{j}(x)\,\psi_{k}^{\prime}(x) -\,\psi_{k}(x)\,\psi_{j}^{\prime}(x)\right)\,d\,x\right\}_{_{0\leq
j,k\leq 2\,n-1}}
\end{equation*}
to be the direct sum of $N$ copies of
\begin{equation*}
  Z=\left(\begin{array}{cc} 0 & 1 \\ -1 & 0 \end{array}\right).
\end{equation*}
Hence $M^{-1}=-M$ where $M^{-1}=\{\mu_{j,k}\}_{j,k=0,2N-1}$. Moreover, with our above choice, $ \mu_{j,k}=0$ if $j,k$ have the same parity or $|j-k|>1$, and $\mu_{2k,2j+1}=\delta_{jk}=-\mu_{2j+1,2k}$ for $j,k=0.\ldots,N-1$.
Therefore
\begin{equation*}
  \begin{aligned}
    S_{N}(x,y) & =-\sum_{j,k=0}^{2N-1}\psi_{j}^{\prime}(x)\mu_{jk}\psi_{k}(y)\\
    & = -\sum_{j=0}^{N-1}\psi_{2j}^{\prime}(x)\,\psi_{2j+1}(y) + \sum_{j=0}^{N-1}\psi_{2j+1}^{\prime}(x)\,\psi_{2j}(y) \\
    & =   \frac{1}{2}\left[ \sum_{j=0}^{N-1}\varphi_{_{2j+1}}\left(x\right) \varphi_{_{2j+1}}\left(y\right)-\sum_{j=0}^{N-1}\varphi_{_{2j+1}}^{\prime}\left(x\right)\epsilon\,\varphi_{_{2j+1}}\left(y\right) \right].
  \end{aligned}
\end{equation*}
Recall that the $H_{j}$ satisfy  the differentiation formulas (see for example \cite{Andr1}, p. 280)
\begin{equation}
  \label{hermrec1}
  H_{j}^{\prime}(x) = 2\,x\,H_{j}(x)-H_{j-1}(x)\quad j=1,2,\ldots
\end{equation}
\begin{equation}
  \label{hermrec2}
  H_{j}^{\prime}(x) = 2\,j\,H_{j-1}(x)\quad j=1,2,\ldots
\end{equation}
Combining (\ref{gsevarphi}) and (\ref{hermrec1}) yields
\begin{equation}
  \label{eq:2}
  \varphi_{j}^{\prime}(x) = x\,\varphi_{j}(x) - \frac{c_{j+1}}{c_{j}}\,\varphi_{j+1}(x).
\end{equation}
Similarly, from (\ref{gsevarphi}) and (\ref{hermrec2}) we have
\begin{equation}
  \label{eq:3}
  \varphi_{j}^{\prime}(x) = -x\,\varphi_{j}(x) + 2j\,\frac{c_{j-1}}{c_{j}}\,\varphi_{j-1}(x).
\end{equation}
Combining \eqref{eq:2} and \eqref{eq:3}, we obtain
\begin{equation}
  \label{eq:4}
    \varphi_{j}^{\prime}(x) = \sqrt{\frac{j}{2}}\,\varphi_{j-1}(x) - \sqrt{\frac{j+1}{2}}\,\varphi_{j+1}(x).
\end{equation}
Let $\varphi=\left(\begin{array}{c}\varphi_{1} \\ \varphi_{2}\\ \vdots \end{array}\right)$ and $\varphi^{\prime}=\left(\begin{array}{c}\varphi_{1}^{\prime} \\ \varphi_{2}^{\prime}\\ \vdots \end{array}\right)$. Then we can rewrite \eqref{eq:4} as
\begin{equation*}
  \varphi^{\prime}=A\,\varphi
\end{equation*}
where $A=\{a_{j,k}\}$ is the infinite antisymmetric tridiagonal matrix with $a_{j,j-1}=\sqrt{\frac{j}{2}}$. Hence,
\begin{equation*}
  \varphi_{j}^{\prime}(x)=\sum_{k\geq 0}a_{jk}\varphi_{k}(x).
\end{equation*}
Moreover, using the fact that $D\,\epsilon=\epsilon\,D=I$ we also have
\begin{equation*}
  \begin{aligned}\varphi_{j}\left(x\right)
    & =\epsilon\,\varphi_{j}^{\prime}\left(x\right)=\epsilon\left(\sum_{k\geq 0}a_{jk}\varphi_{k}\left(x\right)\right)
    \\ & = \sum_{k\geq 0}a_{jk}\epsilon\,\varphi_{k}\left(x\right).
  \end{aligned}
\end{equation*}
Combining the above results, we have
\begin{equation*}
  \begin{aligned}
    \sum_{j=0}^{N-1}\varphi_{2j+1}^{\prime}\left(x\right)\epsilon\,\varphi_{2j+1}\left(y\right) & = \sum_{j=0}^{N-1}\sum_{k\geq 0} a_{2j+1,k}\varphi_{k}\left(x\right)\epsilon\,\varphi_{2j+1}\left(x\right)\\
    & = -\sum_{j=0}^{N-1}\sum_{k\geq 0} a_{k,2j+1}\varphi_{k}\left(x\right)\epsilon\,\varphi_{2j+1}\left(x\right).
  \end{aligned}
\end{equation*}
Note that $a_{k,2j+1}=0$ unless $|k-(2j+1)|=1$, that is unless $k$ is even. Thus we can rewrite the sum as 
\begin{equation*}
  \begin{aligned} 
    \sum_{j=0}^{N-1}\varphi_{2j+1}^{\prime}\left(x\right) \epsilon\,\varphi_{2j+1}\left(y\right)  & =-\sum_{\substack{k,j\geq 0\\ k\ \textrm{even} \\ k\leq 2N}} a_{k,j}\varphi_{k}\left(x\right)\epsilon\,\varphi_{j}\left(y\right) - a_{_{2N,2N+1}}\varphi_{_{2N}}\left(x\right)\epsilon\,\varphi_{_{2N+1}}\left(y\right)\\
    &  = -\sum_{\substack{k\geq 0\\ k\ \textrm{even} \\ k\leq 2N}} \varphi_{k}\left(x\right)\sum_{j\geq 0}a_{k,j}\epsilon\,\varphi_{j}\left(y\right) + a_{_{2N,2N+1}}\,\varphi_{_{2N}}\left(x\right)\epsilon\,\varphi_{_{2N+1}}\left(y\right)
  \end{aligned}
\end{equation*}
where the last term takes care of the fact that we are counting  an extra term in the sum that was not present before. The sum over $j$  on the right is just $\varphi_{k}(y)$, and  $a_{_{2N,2N+1}}=-\sqrt{\frac{2N+1}{2}}$. Therefore
\begin{equation*}
  \begin{aligned}
    \sum_{j=0}^{N-1}\varphi_{2j+1}^{\prime}\left(x\right)\epsilon\,\varphi_{2j+1}\left(y\right) & = \sum_{\substack{k\geq 0\\ k\ \textrm{even} \\ k\leq 2N}} \varphi_{k}\left(x\right)\varphi_{k}\left(y\right)-\sqrt{\frac{2N+1}{2}}\,\varphi_{_{2N}}\left(x\right)\epsilon\,\varphi_{_{2N+1}}\left(y\right)\\
    & = \sum_{j=0}^{N} \varphi_{2j}\left(x\right)\varphi_{2j}\left(y\right)-\sqrt{\frac{2N+1}{2}}\,\varphi_{_{2N}}\left(x\right)\epsilon\,\varphi_{_{2N+1}}\left(y\right).
  \end{aligned}
\end{equation*}
It follows that
\begin{equation*}
  \begin{aligned}
    S_{N}(x,y) =   \frac{1}{2}\left[ \sum_{j=0}^{2N}\varphi_{_{j}}\left(x\right) \varphi_{_{j}}\left(y\right)-\sqrt{\frac{2N+1}{2}}\, \varphi_{_{2N}}\left(x\right)\epsilon\,\varphi_{_{2N+1}} \left(y\right) \right].
  \end{aligned}
\end{equation*}
We redefine
\begin{equation}
  \label{newgseS}
  S_{N}(x,y) = \sum_{n=0}^{2N}\varphi_{_{n}}\left(x\right) \varphi_{_{n}}\left(y\right)=S_{N}(y,x)
\end{equation}
so that the top left entry of $K_{N}(x,y)$ is
\begin{equation*}
  S_{N}(x,y) +  \sqrt{\frac{2N+1}{2}}\, \varphi_{_{2N}}\left(x\right)\epsilon\,\varphi_{_{2N+1}}\left(y\right).
\end{equation*}
If $S_{N}$ is the operator with kernel $S_{N}(x,y)$ then integration by parts gives
\begin{equation*}
  S_{N}D f = \int_{\mathbb{R}} S(x,y) \frac{d}{d\,y}f(y)\,d\,y = \int_{\mathbb{R}}  \left(-\frac{d}{d\,y} S_{N}(x,y)\right) f(y)\,d\,y,
\end{equation*}
so that $-\frac{d}{d\,y} S_{N}(x,y)$ is in fact the kernel of $S_{N}D$. Therefore \eqref{eq:5} now holds with $K_{4,N}$ being the integral operator with matrix kernel $K_{4,N}(x,y)$ whose $(i,j)$--entry $K_{4,N}^{(i,j)}(x,y)$ is given by
\begin{equation*}
  \begin{aligned}
    K_{4,N}^{(1,1)}(x,y)& =\frac{1}{2}\left[ S_{N}(x,y) + \sqrt{\frac{2N+1}{2}}\, \varphi_{_{2N}}\left(x\right)\epsilon\,\varphi_{_{2N+1}} \left(y\right)\right], \\
    K_{4,N}^{(1,2)}(x,y) & = \frac{1}{2}\left[ SD_{N}(x,y) -\frac{d}{d\,y} \, \left(\sqrt{\frac{2N+1}{2}}\, \varphi_{_{2N}}\left(x\right)\epsilon\,\varphi_{_{2N+1}} \left(y\right) \right) \right], \\
    K_{4,N}^{(2,1)}(x,y) & = \frac{\epsilon}{2}\left[ S_{N}(x,y) +  \sqrt{\frac{2N+1}{2}}\, \varphi_{_{2N}}\left(x\right)\epsilon\,\varphi_{_{2N+1}}\left(y\right) \right],\\ 
    K_{4,N}^{(2,2)}(x,y) & = \frac{1}{2}\left[ S_{N}(x,y) + \sqrt{\frac{2N+1}{2}}\, \epsilon\,\varphi_{_{2N+1}}\left(x\right)\,\varphi_{_{2N}}\left(y\right)\right].
  \end{aligned}
\end{equation*}
We let $2N+1\to N$ so that $N$ is assumed to be odd from now on (this will not matter in the  end since we will take $N\to\infty$). Therefore the $K_{4,N}^{(i,j)}(x,y)$ are given by
\begin{equation*}
  \begin{aligned}
    K_{4,N}^{(1,1)}(x,y)& =\frac{1}{2}\left[ S_{N}(x,y) + \sqrt{\frac{N}{2}}\, \varphi_{_{N-1}}(x)\epsilon\,\varphi_{_{N}} (y)\right], \\
    K_{4,N}^{(1,2)}(x,y) & = \frac{1}{2}\left[ SD_{N}(x,y) -\sqrt{\frac{N}{2}}\, \varphi_{_{N-1}}(x)\,\varphi_{_{N}} (y) \right], \\
    K_{4,N}^{(2,1)}(x,y) & = \frac{\epsilon}{2}\left[ S_{N}(x,y) +  \sqrt{\frac{N}{2}}\, \varphi_{_{N-1}}(x)\epsilon\,\varphi_{_{N}}(y) \right], \\
    K_{4,N}^{(2,2)}(x,y) & = \frac{1}{2}\left[ S_{N}(x,y) + \sqrt{\frac{N}{2}}\, \epsilon\,\varphi_{_{N}} (x)\,\varphi_{_{N-1}}(y)\right],
  \end{aligned}
\end{equation*}
where
\begin{equation*}
S_{N}(x,y) = \sum_{n=0}^{N-1}\varphi_{_{n}}(x) \varphi_{_{n}}(y).
\end{equation*}
Define
\begin{equation*}
  \varphi(x)=\left(\frac{N}{2}\right)^{1/4}\varphi_{_{N}}(x)\qquad \psi(x)=\left(\frac{N}{2}\right)^{1/4}\varphi_{_{N-1}}(x),
\end{equation*}
so that
\begin{equation*}
  \begin{aligned}
    K_{4,N}^{(1,1)}(x,y)& =\frac{1}{2}\,\rchi(x)\,\left[ S_{N}(x,y) + \psi(x)\epsilon\,\varphi(y)\right] \,\rchi(y),\\ 
    K_{4,N}^{(1,2)}(x,y) & = \frac{1}{2}\,\rchi(x)\,\left[ SD_{N}(x,y) - \psi(x)\varphi(y)  \right]\,\rchi(y),\\
    K_{4,N}^{(2,1)}(x,y) & = \frac{1}{2}\,\rchi(x)\,\left[ \epsilon S_{N}(x,y) + \epsilon\,\psi(x)\epsilon\,\varphi(y) \right]\,\rchi(y),\\
    K_{4,N}^{(2,2)}(x,y) & = \frac{1}{2}\,\rchi(x)\,\left[ S_{N}(x,y) + \epsilon\,\varphi (x)\,\psi(y)\right]\,\rchi(y).
  \end{aligned}
\end{equation*}
Notice that
\begin{equation*}
  \begin{aligned}
    \frac{1}{2}\,\rchi\,\left( S + \psi\otimes\,\epsilon\,\varphi \right)\,\rchi & \doteq K_{4,N}^{(1,1)}(x,y),\\
    \frac{1}{2}\,\rchi\,\left(SD - \psi\otimes\,\varphi  \right)\,\rchi &\doteq K_{4,N}^{(1,2)}(x,y),\\
    \frac{1}{2}\,\rchi\,\left(\epsilon\,S + \epsilon\,\psi\otimes\,\epsilon\,\varphi \right)\,\rchi &\doteq K_{4,N}^{(2,1)}(x,y),\\
    \frac{1}{2}\,\rchi\,\left(S + \epsilon\,\varphi\otimes\,\epsilon\,\psi \right)\,\rchi &\doteq K_{4,N}^{(2,2)}(x,y).\\
  \end{aligned}
\end{equation*}
Therefore
\begin{equation}
  \label{eq:6}
  K_{4,N}=\frac{1}{2}\,\rchi\,\left(\begin{array}{cc} S + \psi\otimes\,\epsilon\,\varphi  & SD - \psi\otimes\,\varphi \\
      \epsilon\,S + \epsilon\,\psi\otimes\,\epsilon\,\varphi  &  S + \epsilon\,\varphi\otimes\,\psi  \end{array}\right)\,\rchi.
\end{equation}
Note that this is identical to the corresponding operator for $\beta=4$ obtained by Tracy and Widom in \cite{Trac2}, 
the only difference being that $\varphi$, $\psi$, and hence also $S$, are redefined to depend on $\lambda$. 
This will affect boundary conditions for the differential equations we will obtain later. 

\subsection{Edge scaling}

\subsubsection{Reduction of the determinant}

We want to compute the Fredholm determinant \eqref{eq:5} with $K_{4,N}$ given by \eqref{eq:6} and $f=\rchi_{(t,\infty)}$. This is the determinant of an operator on $L^{2}(J)~\times~L^{2}(J)$. Our first task will be to rewrite the determinant as that of an operator on $L^{2}(J)$. This part follows exactly the proof in \cite{Trac2}. To begin, note that 
\begin{equation}
 \label{eq:7}
  \com{S}{D}=\varphi\otimes\psi + \psi\otimes\varphi
\end{equation}
so that, using the fact that $D\,\epsilon=\epsilon\,D=I$,
\begin{eqnarray}
 \label{eq:8}
  \com{\epsilon}{S} &=&\epsilon\,S-S\,\epsilon\nonumber\\
  &=& \epsilon\,S\,D\,\epsilon-\epsilon\,D\,S\,\epsilon = \epsilon\,\com{S}{D}\,\epsilon\nonumber\\
  &=& \epsilon\,\varphi\otimes\psi\,\epsilon + \epsilon\,\psi\otimes\varphi\,\epsilon\nonumber\\
  &=&\epsilon\,\varphi\otimes\epsilon^{t}\psi + \epsilon\,\psi\otimes\epsilon^{t}\,\varphi\nonumber\\
  &=& - \epsilon\,\varphi\otimes\epsilon\,\psi -
  \epsilon\,\psi\otimes\epsilon\,\varphi,
\end{eqnarray}
where the last equality follows from the fact that $\epsilon^{t}=-\epsilon$. We thus have
\begin{eqnarray*}
  D\,\left(\epsilon\,S + \epsilon\,\psi\otimes\epsilon\,\varphi\right) & = & S + \psi\otimes\epsilon\,\varphi,\\
  D\,\left(\epsilon\,S\,D - \epsilon\,\psi\otimes\varphi\right) & =
  & S\,D - \psi\otimes\varphi.
\end{eqnarray*} 
The expressions on the right side are the top matrix entries in \eqref{eq:6}. Thus the first row of $K_{4,N}$ is, as a vector,
\begin{equation*}
  D\,\left(\epsilon\,S + \epsilon\,\psi\otimes\epsilon\,\varphi, \epsilon\,S\,D - \epsilon\,\psi\otimes\varphi\right).
\end{equation*} 
Now \eqref{eq:8} implies that
\begin{equation*}  
  \epsilon\,S + \epsilon\,\psi\otimes\epsilon\,\varphi = S\,\epsilon -\epsilon\,\varphi\otimes\epsilon\,\psi.
\end{equation*} 
Similarly \eqref{eq:7} gives
\begin{equation*} 
  \epsilon\,\com{S}{D} = \epsilon\,\varphi\otimes\psi + \epsilon\,\psi\otimes\varphi,
\end{equation*}
so that
\begin{equation*}
  \epsilon\,S\,D - \epsilon\,\psi\otimes\varphi = \epsilon\,D\,S + \epsilon\,\varphi\otimes\psi = S + \epsilon\,\varphi\otimes\psi.
\end{equation*}
Using these expressions we can rewrite the first row of $K_{4,N}$ as
\begin{equation*}
  D\,\left(S\,\epsilon - \epsilon\,\varphi\otimes\epsilon\,\psi, S + \epsilon\,\varphi\otimes\psi\right). 
\end{equation*}
Now use \eqref{eq:8} to show the second row of $K_{4,N}$ is
\begin{equation*}
  \left(S\,\epsilon - \epsilon\,\varphi\otimes\epsilon\,\psi, S + \epsilon\,\varphi\otimes\psi\right). 
\end{equation*} 
Therefore,
\begin{eqnarray*}
  K_{4,N} &=& \rchi\,\left(
    \begin{array}{cc}
      D\,\left(S\,\epsilon - \epsilon\,\varphi\otimes\epsilon\,\psi\right) & D\,\left(S + \epsilon\,\varphi\otimes\psi\right) \\
      S\,\epsilon - \epsilon\,\varphi\otimes\epsilon\,\psi &  S + \epsilon\,\varphi\otimes\psi
    \end{array}
  \right)\,\rchi \\
  & =& \left(
    \begin{array}{cc}
      \rchi\,D & 0 \\
      0 & \rchi
    \end{array}
  \right) 
  \left(
    \begin{array}{cc}
      \left(S\,\epsilon - \epsilon\,\varphi\otimes\epsilon\,\psi\right)\,\rchi & \left(S + \epsilon\,\varphi\otimes\psi\right)\,\rchi \\
      \left(S\,\epsilon - \epsilon\,\varphi\otimes\epsilon\,\psi\right)\rchi & \left(S + \epsilon\,\varphi\otimes\psi\right)\rchi
    \end{array}
  \right).
\end{eqnarray*}
Since $K_{4,N}$ is of the form $A\,B$, we can use \ref{detthm}  and deduce that $D_{4,N}(s,\lambda)$ is unchanged if instead we take $ K_{4,N}$ to be
\begin{eqnarray*} 
  K_{4,N} & = & \left(
    \begin{array}{cc}
      \left(S\,\epsilon - \epsilon\,\varphi\otimes\epsilon\,\psi\right)\,\rchi & \left(S + \epsilon\,\varphi\otimes\psi\right)\,\rchi \\
      \left(S\,\epsilon - \epsilon\,\varphi\otimes\epsilon\,\psi\right)\,\rchi & \left(S + \epsilon\,\varphi\otimes\psi\right)\,\rchi
    \end{array}
  \right) 
  \left(
    \begin{array}{cc}
      \rchi\,D & 0 \\
      0 & \rchi
    \end{array}
  \right)\\
  & = & \left(\begin{array}{cc}
      \left(S\,\epsilon - \epsilon\,\varphi\otimes\epsilon\,\psi\right)\,\rchi\,D & \left(S + \epsilon\,\varphi\otimes\psi\right)\,\rchi \\
      \left(S\,\epsilon - \epsilon\,\varphi\otimes\epsilon\,\psi\right)\rchi\,D & \left(S + \epsilon\,\varphi\otimes\psi\right)\rchi
    \end{array}\right).
\end{eqnarray*}
Therefore
\begin{equation}
  D_{4,N}(s,\lambda)=\det\left(
    \begin{array}{cc}
      I - \frac{1}{2}\,\left(S\,\epsilon - \epsilon\,\varphi\otimes\epsilon\,\psi\right)\,\lambda\,\rchi\,D & - \frac{1}{2}\,\left(S + \epsilon\,\varphi\otimes\psi\right)\,\lambda\,\rchi \\
      - \frac{1}{2}\,\left(S\,\epsilon - \epsilon\,\varphi\otimes\epsilon\,\psi\right)\,\lambda\,\rchi\,D & I - \frac{1}{2}\,\left(S + \epsilon\,\varphi\otimes\psi\right)\,\lambda\,\rchi
    \end{array}
  \right).
\end{equation}
Now we perform row and column operations on the matrix to simplify it, which do not change the Fredholm determinant. Justification of these operations is given in \cite{Trac2}. We start by subtracting row 1 from row 2 to get
\begin{equation*}  
  \left(\begin{array}{cc}
      I - \frac{1}{2}\,\left(S\,\epsilon - \epsilon\,\varphi\otimes\epsilon\,\psi\right)\,\lambda\,\rchi\,D & - \frac{1}{2}\,\left(S + \epsilon\,\varphi\otimes\psi\right)\,\lambda\,\rchi \\
      - I  & I \end{array}\right).
\end{equation*}
Next, adding column 2 to column 1 yields
\begin{equation*} 
  \left(\begin{array}{cc}
      I - \frac{1}{2}\,\left(S\,\epsilon - \epsilon\,\varphi\otimes\epsilon\,\psi\right)\,\lambda\,\rchi\,D - \frac{1}{2}\,\left(S + \epsilon\,\varphi\otimes\psi\right)\,\lambda\,\rchi & - \frac{1}{2}\,\left(S + \epsilon\,\varphi\otimes\psi\right)\,\lambda\,\rchi \\
      0  & I 
    \end{array} 
  \right). 
\end{equation*}

Thus the determinant we want equals the determinant of
\begin{equation}
  \label{eq:9}
  I - \frac{1}{2}\,\left(S\,\epsilon - \epsilon\,\varphi\otimes\epsilon\,\psi\right)\,\lambda\,\rchi\,D - \frac{1}{2}\,\left(S + \epsilon\,\varphi\otimes\psi\right)\,\lambda\,\rchi 
\end{equation}
So we have reduced the problem from the computation of the Fredholm determinant of an operator on $L^{2}(J)~\times~L^{2}(J)$, to that of an operator on $L^{2}(J)$.

\subsubsection{Differential equations}

Next we want to write the operator in \eqref{eq:9} in the form 
\begin{equation} 
  \left(I-  K_{2, N}\right)\left(I - \sum_{i=1}^{L}\alpha_{i}\otimes\beta_{i}\right), 
\end{equation}
where the $\alpha_{i}$  and $\beta_{i}$ are functions in $L^{2}(J)$. In other words, we want to rewrite the determinant for the GSE case as a finite dimensional perturbation of the corresponding GUE determinant. The Fredholm determinant of the product is then the product of the determinants. The limiting form for the GUE part is already known, and we can just focus on finding a limiting form for the determinant of the finite dimensional piece. It is here that the proof must be modified from that in \cite{Trac2}. A little rearrangement of \eqref{eq:9} yields (recall $\epsilon^{t}=-\epsilon$)
\begin{equation*}
  I - \frac{\lambda}{2}\,S\,\rchi  - \frac{\lambda}{2}\,S\,\epsilon\,\rchi\,D - \frac{\lambda}{2}\,\epsilon\,\varphi\,\otimes\,\rchi\,\psi - \frac{\lambda}{2}\,\epsilon\,\varphi\,\otimes\,\psi\,\epsilon\,\rchi\,D.
\end{equation*}
Writing $\epsilon\,\com{\rchi}{D}+\rchi$ for $\epsilon\,\rchi\,D$ and simplifying gives
\begin{equation*}
  I - \lambda\,S\,\rchi - \lambda\,\epsilon\,\varphi\,\otimes\,\psi\,\rchi - \frac{\lambda}{2}\,S\,\epsilon\,\com{\rchi}{D} - \frac{\lambda}{2}\,\epsilon\,\varphi\,\otimes\,\psi\,\epsilon\,\com{\rchi}{D}.
\end{equation*}
Let $\sqrt{\lambda}\,\varphi\to\varphi$, and $\sqrt{\lambda}\,\psi\to\psi$ so that $\lambda\,S\to S$ and \eqref{eq:9} goes to
\begin{equation*}
  I - S\,\rchi - \epsilon\,\varphi\,\otimes\,\psi\,\rchi - \frac{1}{2}\,S\,\epsilon\,\com{\rchi}{D} - \frac{1}{2}\,\epsilon\,\varphi\,\otimes\,\psi\,\epsilon\,\com{\rchi}{D}.
\end{equation*}
Now we define $R:=(I-S\,\rchi)^{-1}\,S\,\rchi=(I-S\,\rchi)^{-1}-I$ (the resolvent operator of $S\,\rchi$), whose kernel we denote by $R(x,y)$, and $Q_{\epsilon}:=(I-S\,\rchi)^{-1}\,\epsilon\,\varphi$. Then \eqref{eq:9} factors into
\begin{equation*}
  A =(I - S\,\rchi)\,B.
\end{equation*}
where $B$ is
\begin{equation*}
  I - Q_{\epsilon}\,\otimes\,\rchi\,\psi - \frac{1}{2}(I+R)\,S\,\epsilon\,\com{\rchi}{D} - \frac{1}{2}\,(Q_{\epsilon}\,\otimes\,\psi)\,\epsilon\,\com{\rchi}{D}
\end{equation*}
Hence
\begin{equation*}
  D_{4,N}(t,\lambda) = D_{2,N}(t,\lambda)\,\det(B).
\end{equation*}
In order to find $\det(B)$ we use the identity
\begin{equation}
  \label{eq:14}
  \epsilon\,\com{\rchi}{D}=\sum_{k=1}^{2m}(-1)^{k}\,\epsilon_{k}\otimes\delta_{k},
\end{equation}
where $\epsilon_{k}$ and $\delta_{k}$ are the functions $\epsilon(x-a_{k})$ and $\delta(x-a_{k})$ respectively, and the $a_{k}$ are the endpoints of the (disjoint) intervals considered, $J=\cup_{k=1}^{m}(a_{2\,k-1},a_{2\,k})$. In our case $m=1$ and $a_{1}=t$, $a_{2}=\infty$. We also make use of the fact that
\begin{equation}
  a\otimes b\cdot c\otimes d= \inprod{b}{c}\cdot a\otimes d
\end{equation}
where $\inprod{.}{.}$ is the usual $L^{2}$--inner product. Therefore
\begin{align*}
  (Q_{\epsilon}\otimes\psi)\,\epsilon\,\com{\rchi}{D} & = \sum_{k=1}^{2} (-1)^{k}Q_{\epsilon}\otimes\psi\cdot\epsilon_{k}\otimes\delta_{k} \\ 
  & =\sum_{k=1}^{2} (-1)^{k}\inprod{\psi}{\epsilon_{k}}\,Q_{\epsilon}\otimes\,\delta_{k}.
\end{align*}
It follows that
\begin{equation}
  \frac{D_{4,N}(t,\lambda)}{D_{2,N}(t,\lambda)} 
\end{equation}
is the determinant of
\begin{equation}
I - Q_{\epsilon}\otimes\rchi\psi - \frac{1}{2}\,\sum_{k=1}^{2} (-1)^{k}\left[(S+R\,S)\,\epsilon_{k} +\inprod{\psi}{\epsilon_{k}}\,Q_{\epsilon}\right]\otimes\delta_{k}.
\end{equation}
We now specialize to the case of one interval $J=(t,\infty)$, so $m=1$, $a_{1}=t$ and $a_{2}=\infty$. We write $\epsilon_{t}=\epsilon_{1}$, and $\epsilon_{\infty}=\epsilon_{2}$, and similarly for $\delta_{k}$. Writing out the terms in the summation and using the fact that
\begin{equation}
  \epsilon_{\infty}=-\frac{1}{2},
\end{equation}
yields
\begin{equation}
  I - Q_{\epsilon}\otimes\rchi\psi  + \frac{1}{2}\,\left[(S+R\,S)\,\epsilon_{t}+ \inprod{\psi}{\epsilon_{t}}\,Q_{\epsilon}\right]\otimes\delta_{t} + \frac{1}{4}\,\left[(S+R\,S)\,1 + \inprod{\psi}{1}\,Q_{\epsilon}\right]\otimes\delta_{\infty}
\end{equation}
Now we can use the formula
\begin{equation}
  \label{eq:16}
  \det\left(I-\sum_{i=1}^{L}\alpha_{i}\otimes\beta_{i}\right)=\det\left(\delta_{jk}-\inprod{\alpha_{j}}{\beta_{k}}\right)_{1\leq j,k\leq L}
\end{equation}
In order to simplify the notation in preparation for the computation of the various inner products, define
\begin{align}
  Q(x,\lambda,t)&:=(I-S\,\rchi)^{-1}\,\varphi, & P(x,\lambda,t)&:=(I-S\,\rchi)^{-1}\,\psi, \nonumber\\
  Q_{\epsilon}(x,\lambda,t)&:=(I-S\,\rchi)^{-1}\,\epsilon\,\varphi, & P_{\epsilon}(x,\lambda,t)&:=(I-S\,\rchi)^{-1}\,\epsilon\,\psi,
\end{align}
\begin{align}
  \label{eq:10}
  q_{_{N}}&:=Q(t,\lambda,t), & p_{_{N}}&:=P(t,\lambda,t),\nonumber\\ 
  q_{\epsilon}&:=Q_{\epsilon}(t,\lambda,t), & p_{\epsilon}&:=P_{\epsilon}(t,\lambda,t),\nonumber\\
  u_{\epsilon}&:=\inprod{Q}{\rchi\,\epsilon\,\varphi}=\inprod{Q_{\epsilon}}{\rchi\,\varphi}, & v_{\epsilon}&:=\inprod{Q}{\rchi\,\epsilon\,\psi}=\inprod{P_{\epsilon}}{\rchi\,\psi}, \nonumber\\
  \tilde{v}_{\epsilon}&:=\inprod{P}{\rchi\,\epsilon\,\varphi}=\inprod{Q_{\epsilon}}{\rchi\,\varphi}, & w_{\epsilon}&:=\inprod{P}{\rchi\,\epsilon\,\psi}=\inprod{P_{\epsilon}}{\rchi\,\psi},
\end{align}
\begin{equation}
  \label{eq:11}
  \mathcal{P}_{4} := \int_{\mathbb{R}}\,\epsilon_{t}(x)\,P(x,t)\,d\,x, \quad \mathcal{Q}_{4} := \int_{\mathbb{R}}\,\epsilon_{t}(x)\,Q(x,t)\,d\,x, \quad  \mathcal{R}_{4} := \int_{\mathbb{R}}\,\epsilon_{t}(x)\,R(x,t)\,d\,x,
\end{equation}
where we remind the reader that $\epsilon_{t}$ stands for the function $\epsilon(x-t)$. Note that all quantities in \eqref{eq:10} and \eqref{eq:11} are functions of $t$ and $\lambda$ alone. Furthermore, let
\begin{equation}
  c_{\varphi} = \epsilon\,\varphi(\infty)=\frac{1}{2}\int_{-\infty}^{\infty}\varphi(x)\,d\,x, \qquad c_{\psi} = \epsilon\,\psi(\infty)=\frac{1}{2}\int_{-\infty}^{\infty}\psi(x)\,d\,x.
\end{equation}
Recall from the previous section that when $\beta=4$ we take $N$ to be odd. It follows that $\varphi$ and $\psi$ are odd and even functions respectively. Thus when $\beta=4$, $c_{\varphi}=0$ while computation using known integrals for the Hermite polynomials gives
\begin{equation}
c_{\psi}=(\pi N)^{1/4} 2^{-3/4-N/2}\,\frac{(N!)^{1/2}}{(N/2)!}\,\sqrt{\lambda}.
\end{equation}
Hence computation yields
\begin{equation}
 \lim_{N\to\infty}c_{\psi}= \sqrt{\frac{\lambda}{2}}.
\end{equation}
At $t=\infty$,
\begin{equation}
  u_{\epsilon}(\infty)=0, \quad  q_{\epsilon}(\infty)=c_{\varphi}
\end{equation}
\begin{equation}
  \mathcal{P}_{4}(\infty) = -c_{\psi}, \quad \mathcal{Q}_{4}(\infty) = -c_{\varphi},\quad \mathcal{R}_{4}(\infty) = 0.
\end{equation}
In \eqref{eq:16}, $L=3$ and if we denote $a_{4}=\inprod{\psi}{\epsilon_{t}}$, then we have explicitly
\begin{equation*}
  \alpha_{1}=Q_{\epsilon}, \quad \alpha_{2} = -\frac{1}{2}\,\left[(S+R\,S)\,\epsilon_{t} + a_{4}\,Q_{\epsilon}\right],\quad\alpha_{3} = -\frac{1}{4}\,\left[(S+R\,S)\,1 + \inprod{\psi}{1}\,Q_{\epsilon}\right],
\end{equation*}
\begin{equation*}
  \beta_{1}=\rchi\psi, \quad  \beta_{2}=\delta_{t}, \quad  \beta_{3}=\delta_{\infty}.
\end{equation*}
However notice that
\begin{equation}
  \inprod{\left(S+R\,S\right)\,\epsilon_{t}}{\delta_{\infty}}=\inprod{\epsilon_{t}}{\delta_{\infty}}=0,\quad \inprod{\left(S+R\,S\right)\,1}{\delta_{\infty}}=\inprod{1}{R_{\infty}}=0
\end{equation}
and $\inprod{Q_{\epsilon}}{\delta_{\infty}}=c_{\varphi}=0$. Therefore the terms involving $\beta_{3}=\delta_{\infty}$ are all $0$ and we can discard them reducing our computation to that of a $2\times 2$ determinant instead with
\begin{equation}
  \alpha_{1}=Q_{\epsilon}, \quad \alpha_{2} = -\frac{1}{2}\,\left[(S+R\,S)\,\epsilon_{t} + a_{4}\,Q_{\epsilon}\right], \quad \beta_{1}=\rchi\psi, \quad  \beta_{2}=\delta_{t}.
\end{equation}
Hence
\begin{align}
  \inprod{\alpha_{1}}{\beta_{1}} & = \tilde{v}_{\epsilon}, \quad \inprod{\alpha_{1}}{\beta_{2}}=q_{\epsilon}, \\
  \inprod{\alpha_{2}}{\beta_{1}} & =-\frac{1}{2} \,\left(\mathcal{P}_{4} - a_{4}\, + a_{4}\,\tilde{v}_{\epsilon}\right), \\
  \inprod{\alpha_{2}}{\beta_{2}} & = -\frac{1}{2}\,\left(\mathcal{R}_{4} + a_{4}\,q_{\epsilon}\right)\label{eq:12}.
 \end{align}
We want the limit of the determinant
\begin{equation}
  \det\left(\delta_{jk}-\inprod{\alpha_{j}}{\beta_{k}}\right)_{1\leq j,k\leq 2},
\end{equation}
as $N\to \infty$. In order to get our hands on the limits of the individual terms involved in the determinant, we will find differential equations for them first as in \cite{Trac2}. Adding $a_{4}/2$ times row 1 to row 2 shows that ${a}_{4}$ falls out of the determinant, so we will not need to find differential equations for it. Thus our determinant is now
{\large
  \begin{equation}
    \det\left(
      \begin{array}{cc}
        1 - \tilde{v}_{\epsilon} & -q_{\epsilon} \\[6pt]
        \frac{1}{2} \,\mathcal{P}_{4} & 1 + \frac{1}{2}\,\mathcal{R}_{4} \\[6pt]
      \end{array}
    \right).
  \end{equation}
}
Proceeding as in \cite{Trac2} we find the following differential equations
\begin{align}
  \frac{d}{d\,t}\,u_{\epsilon} & = - q_{_{N}}\,q_{\epsilon}, & \frac{d}{d\,t}\,q_{\epsilon} & = q_{_{N}} -q_{_{N}}\, \tilde{v}_{\epsilon} - p_{_{N}}\,u_{\epsilon},\\
 \frac{d}{d\,t}\mathcal{Q}_{4} & = -q_{_{N}}\left(\mathcal{R}_{4} + 1\right), & \frac{d}{d\,t}\mathcal{P}_{4} & = -p_{_{N}}\left(\mathcal{R}_{4} + 1\right),\label{eq:13}\\
 \frac{d}{d\,t}\mathcal{R}_{4} & = -p_{_{N}}\,\mathcal{Q}_{4}-q_{_{N}}\,\mathcal{P}_{4}.   & 
\end{align}
Now we change variable from $t$ to $s$ where $t=\tau(s)= \sqrt{2\,N}+s/(\sqrt{2}\,N^{1/6})$ and take the limit $N\to \infty$, denoting the limits of $ q_{\epsilon}$, $\mathcal{P}_{4}$, $\mathcal{Q}_{4}$, $\mathcal{R}_{4},$ and the common limit of $u_{\epsilon}$ and $\tilde{v}_{\epsilon}$ respectively by $\overline{q}$, $\overline{\mathcal{P}}_{4}$, $\overline{\mathcal{Q}}_{4}$, $\overline{\mathcal{R}}_{4}$ and $\overline{u}$. Also $\overline{\mathcal{P}}_{4}$ and $\overline{\mathcal{Q}}_{4}$ differ by a constant, namely $\overline{\mathcal{Q}}_{4}=\overline{\mathcal{P}}_{4} + \sqrt{2}/{2}$. These limits hold uniformly for bounded $s$ so we can interchange $\lim_{N\to\infty}$ and $\frac{d}{d\,s}$. Also $\lim_{N\to\infty}N^{-1/6}q_{_{N}}=\lim_{N\to\infty}N^{-1/6}p_{_{N}}=q $ , where $q$ is as in \eqref{D2}. We obtain the systems
\begin{equation}
  \frac{d}{d\,s}\,\overline{u} = -\frac{1}{\sqrt{2}}\,q\,\overline{q},\qquad \frac{d}{d\,s}\,\overline{q} = \frac{1}{\sqrt{2}}\,q\,\left(1-2\,\overline{u}\right),
\end{equation}
\begin{equation}
  \frac{d}{d\,s}\overline{\mathcal{P}}_{4} = -\frac{1}{\sqrt{2}}\,q\,\left(\overline{\mathcal{R}}_{4} + 1\right), \qquad \frac{d}{d\,s}\overline{\mathcal{R}}_{4} = -\frac{1}{\sqrt{2}}\,q\,\left(2\,\overline{\mathcal{P}}_{4} + \sqrt{\frac{\lambda}{2}}\right),
\end{equation}
The change of variables $s\to\mu=\int_{s}^{\infty} q(x,\lambda)\,d\,x$ transforms these systems into constant coefficient ordinary differential equations
\begin{equation}
  \frac{d}{d\,\mu}\overline{u} = \frac{1}{\sqrt{2}}\,\overline{q}, \qquad \frac{d}{d\,\mu}\overline{q} = -\frac{1}{\sqrt{2}}\,\left(1-2\,\overline{u}\right),
\end{equation}
\begin{equation}
  \frac{d}{d\,\mu}\overline{\mathcal{P}}_{4} = \frac{1}{\sqrt{2}}\,\left(\overline{\mathcal{R}}_{4} + 1\right), \qquad \frac{d}{d\,\mu}\overline{\mathcal{R}}_{4} = \frac{1}{\sqrt{2}}\,\left(2\,\overline{\mathcal{P}}_{4}+ \sqrt{\frac{\lambda}{2}}\right).
\end{equation}
Since $\lim_{s\to \infty}\mu=0$, corresponding to the boundary values at $t=\infty$ which we found earlier for $\mathcal{P}_{4}, \mathcal{R}_{4}$, we now have initial values at $\mu=0$. Therefore
\begin{equation}
  \overline{u}(\mu=0)=\overline{q}(\mu=0)=0,
\end{equation}
\begin{equation}
  \overline{\mathcal{P}}_{4}(\mu=0)=-\sqrt{\frac{\lambda}{2}}, \qquad \overline{\mathcal{R}}_{4}(\mu=0)=0.
\end{equation}
We use this to solve the systems and get
\begin{align}
  \overline{q} & = \frac{1}{2\sqrt{2}}\,\left(e^{-\mu}-e^{\mu}\right)\\
 \overline{u} & = \frac{1}{2}\,\left(1 - \frac{1}{2}\,e^{\mu} - \frac{1}{2}\,e^{-\mu}\right)\\
 \overline{\mathcal{P}}_{4} & = \frac{1}{2\,\sqrt{2}}\,\left(\frac{2-\sqrt{\lambda}}{2}\,e^{\mu} - \frac{2+\sqrt{\lambda}}{2}\,e^{-\mu} - \sqrt{\lambda}\right) \label{eq:24}\\
 \overline{\mathcal{R}}_{4} & = \frac{2-\sqrt{\lambda}}{4}\,e^{\mu} + \frac{2+\sqrt{\lambda}}{4}\,e^{-\mu} - 1\label{eq:25} 
\end{align}
Substituting these expressions into the determinant gives
\eqref{gsedet}, namely
\begin{equation}
  \label{eq:26}
  D_{4}(s,\lambda)= D_{2}(s,\lambda)\,\cosh^{2}\left(\frac{\mu(s,\lambda)}{2}\right),
\end{equation} 
where $D_{\beta}=\lim_{N\to\infty}D_{\beta,N}$.
Note that even though there are $\lambda$--terms in \eqref{eq:24} and \eqref{eq:25}, these do not appear in the final result \eqref{eq:26}, making it similar to the GUE case where the main conceptual difference between the $m=1$ (largest eigenvalue) case and the general $m$ is the dependence of the function $q$ on $\lambda$. The right hand side of the above formula clearly reduces to the $\beta=4$ Tracy-Widom distribution when we set $\lambda=1$. 
Note that where we have $D_{4}(s,\lambda)$ above, Tracy and Widom (and hence many RMT references) write $D_{4}(s/\sqrt{2},\lambda)$ instead. Tracy and Widom applied the change of variable $s\to s/\sqrt{2}$ in their derivation in \cite{Trac2} so as to agree with Mehta's form of the $\beta=4$ joint eigenvalue density,\footnoterecall{norm} which has $-2x^{2}$ in the exponential in the weight function, instead of $-x^{2}$ in our case. To switch back to the other convention, one just needs to substitute in the argument $s/\sqrt{2}$ for $s$ everywhere in our results. At this point this is just a cosmetic discrepancy, and it does not change anything in our derivations since all the differentiations are done with respect to $\lambda$ anyway. It \underline{does} change conventions for rescaling data while doing numerical work though.

\section{The Distribution of the $m^{th}$ Largest Eigenvalue in the GOE}

\sectionmark{The $m^{th}$ largest eigenvalue in the GOE}

\subsection{The distribution function as a Fredholm determinant}
The GOE corresponds case corresponds to the specialization $\beta=1$ in \eqref{jointdensity} so that
\begin{equation}
\label{eq:17}
  G_{1,N}(t,\lambda)=C_{1}^{(N)}\underset{x_{i}\in\mathbb{R}}{\int\cdots\int} \prod_{j<k}\left|x_{j}-x_{k}\right|\,\prod_{j}^{N}w(x_{j})\,\prod_{j}^{N}\left(1+f(x_{j})\right)\,d\,x_{1}\cdots d\,x_{N}
\end{equation}
where $w(x)=\exp\left(-x^{2}\right)$, $f(x)=-\lambda\,\rchi_{J}(x)$, and $C_{1}^{(N)}$ depends only on $N$. As in the GSE case, we will lump into $C_{1}^{(N)}$ any constants depending only on $N$ that appear in the derivation. A simple argument  at the end will show that the final constant is $1$. These calculations more or less faithfully follow and expand on \cite{Trac1}. We want to use \eqref{debruijn2}, which requires an ordered space. Note that the above integrand is symmetric under permutations, so the integral is $n!$ times the same integral over ordered pairs $x_{1}\leq\ldots\leq x_{N}$. So we can rewrite \ref{eq:17} as
\begin{equation*} 
  (N!)\,\underset{x_{1}\leq\ldots\leq x_{N}\in \mathbb R\,\,}{\int\cdots\int} \prod_{j<k}(x_{k}-x_{j})\,\prod_{i=1}^{N}w(x_{k})\prod_{i=1}^{N}(1+f(x_{k}))\,d\,x_{1}\cdots d\,x_{N},
\end{equation*}
where we can remove the absolute values since the ordering insures that $(x_{j}-x_{i})\geq 0$ for $i<j$. Recall that the Vandermonde determinant is
\begin{equation*}
  \Delta_{N}(x)=\det(x^{j-1}_{k})_{_{1\leq j,k\leq N}}=(-1)^{\frac{N\,(N-1)}{2}}\prod_{j<k}(x_{j}-x_{k}).
\end{equation*}
Therefore what we have inside the integrand above is, up to sign
\begin{equation*}
  \det(x^{j-1}_{k}\,w(x_{k})\,(1+f(x_{k})))_{_{1\leq j,k\leq N}}.
\end{equation*}
Note that the sign depends only on $N$. Now we can use \eqref{debruijn2} with \[\varphi_{j}(x)=x^{j-1}\,w(x)\,(1+f(x)). \]
In using \eqref{debruijn2} we square both sides so that the right hand side is now a determinant instead of a Pfaffian. Therefore $G_{1,N}^{2}(t,\lambda)$ equals
\begin{equation*}
  C^{(N)}_{1}\,\det\left(\int\int\sgn(x-y)x^{j-1}\,y^{k-1}\,(1+f(x))\,w(x)\,w(y)\,d\,x\,d\,y \right)_{_{1\leq j,k\leq N}}.
\end{equation*}
Shifting indices, we can write it as
\begin{equation}
  C^{(N)}_{1}\,\det\left(\int\int\sgn(x-y)x^{j}\,y^{k}\,(1+f(x))\,w(x)\,w(y)\,d\,x\,d\,y \right)_{_{1\leq j,k\leq N-1}}
\end{equation}
where $C^{(N)}_{1}$ is a constant depending only on $N$, and is such that the right side is $1$ if $f\equiv 0$. Indeed this would correspond to the probability that $\lambda_{p}^{GOE(N)}<\infty$, or equivalently to the case where the excluded set $J$ is empty. We can replace $x^{j}$ and $y^{k}$ by any arbitrary polynomials $p_{j}(x)$ and $p_{k}(x)$, of degree $j$ and $k$ respectively, which are obtained by row operations on the matrix. Indeed such operations would not change the determinant. We also replace $\sgn(x-y)$ by $\epsilon(x-y)=\frac{1}{2}\,\sgn(x-y)$ which just produces a factor of $2$ that we absorb in $C^{(N)}_{1}$. Thus $G_{1,N}^{2}(t,\lambda)$ now equals
\begin{equation}
 C^{(N)}_{1}\,\det\left(\int\int\epsilon(x-y)\,p_{j}(x)\,p_{k}(y)\,(1 + f(x))\,(1 + f(y))\,w(x)\,w(y)\,d\,x\,d\,y \right)_{_{0\leq j,k\leq N-1}}.
\end{equation}
Let $\psi_{j}(x)=p_{j}(x)\,w(x)$  so the above integral becomes
\begin{equation}
 C^{(N)}_{1}\,\det\left(\int\int\epsilon(x-y)\,\psi_{j}(x)\,\psi_{k}(y)\,(1+f(x)+f(y)+f(x)\,f(y))\,d\,x\,d\,y \right)_{_{0\leq j,k\leq N-1}}.
\end{equation}
Partially multiplying out the term we obtain
\begin{align*}
  C^{(N)}_{1}\,&\det\left(\int\int\epsilon(x-y)\,\psi_{j}(x)\,\psi_{k}(y)\,d\,x\,d\,y \right. \\ 
  & + \left. \int\int\epsilon(x-y)\,\psi_{j}(x)\,\psi_{k}(y)\,(f(x)+f(y)+f(x)\,f(y))\,d\,x\,d\,y \right)_{_{0\leq j,k\leq N-1}}.
\end{align*}
Define
\begin{equation}
  \label{eq:18}
  M=\left(\int\int\epsilon(x-y)\,\psi_{j}(x)\,\psi_{k}(y)\,d\,x\,d\,y \right)_{_{0\leq j,k\leq N-1}},
\end{equation}
so that $G_{1,N}^{2}(t,\lambda)$ is now
\begin{equation*}
  C^{(N)}_{1}\,\det\left(M + \int\int\epsilon(x-y)\,\psi_{j}(x)\,\psi_{k}(y)\,(f(x)+f(y)+f(x)\,f(y))\,d\,x\,d\,y \right)_{_{0\leq j,k\leq N-1}}.
\end{equation*}
Let $\epsilon$ be the operator defined in \eqref{epsilonop}. We can use operator notation to simplify the expression for $G_{1,N}^{2}(t,\lambda)$ a great deal by rewriting the double integrals as single integrals. Indeed
\begin{eqnarray*}
  \int\int\epsilon(x-y)\,\psi_{j}(x)\,\psi_{k}(y)\,f(x)\,d\,x\,d\,y & = & \int f(x)\,\psi_{j}(x)\,\int\epsilon(x-y)\,\psi_{k}(y)\,d\,y\,d\,x\\ 
  & = & \int f\,\psi_{j}\,\epsilon\,\psi_{k}\,d\,x.
\end{eqnarray*}
Similarly,
\begin{eqnarray*}
  \int\int\epsilon(x-y)\,\psi_{j}(x)\,\psi_{k}(y)\,f(y)\,d\,x\,d\,y & = & - \int\int\epsilon(y-x)\,\psi_{j}(x)\,\psi_{k}(y)\,f(y)\,d\,x\,d\,y \\
  & = & - \int f(y)\,\psi_{k}(y)\,\int\epsilon(y-x)\,\psi_{j}(x)\,d\,x\,d\,y \\
  & = &  - \int f(x)\,\psi_{k}(x)\,\int\epsilon(x-y)\,\psi_{j}(y)\,d\,y\,d\,x \\
  & = &  - \int f\,\psi_{k}\,\epsilon\,\psi_{j}\,d\,x.
\end{eqnarray*}
Finally,
\begin{equation}
  \begin{aligned}
    \int\int \epsilon(x-y)\,\psi_{j}(x)\,&\psi_{k}(y)\,f(x)\,f(y)\,d\,x\,d\,y  \\
    &  = - \int\int\epsilon(y-x)\,\psi_{j}(x)\,\psi_{k}(y)\,f(x)\,f(y)\,d\,x\,d\,y \\
    &  = - \int f(y)\,\psi_{k}(y)\,\int\epsilon(y-x)\,f(x)\,\psi_{j}(x)\,d\,x\,d\,y \\
    &  = - \int f(x)\,\psi_{k}(x)\,\int\epsilon(x-y)\,f(y)\,\psi_{j}(y)\,d\,y\,d\,x \\
    &  = - \int f\,\psi_{k}\,\epsilon\,(f\,\psi_{j})\,d\,x.
  \end{aligned}
\end{equation}
It follows that
\begin{equation}
  G_{1,N}^{2}(t,\lambda)=C^{(N)}_{1}\,\det\left(M + \int\left[ f\,\psi_{j}\,\epsilon\,\psi_{k} - f\,\psi_{k}\,\epsilon\,\psi_{j} -  f\,\psi_{k}\,\epsilon\,(f\,\psi_{j}) \right]\,d\,x \right)_{_{0\leq j,k\leq N-1}}.
\end{equation}
If we let $M^{-1}=\left(\mu_{j\,k} \right)_{_{0\leq j,k\leq N-1}}$, and factor $\det(M)$ out, then $G_{1,N}^{2}(t,\lambda) $ equals
\begin{align}
  C^{(N)}_{1}\,& \det(M)\,\det\left(\vphantom{M^{-1}\cdot\left(\int\left[ f\,\psi_{j}\,\epsilon\,\psi_{k} - f\,\psi_{k}\,\epsilon\,\psi_{j} -  f\,\psi_{k}\,\epsilon\,(f\,\psi_{j}) \right]\,d\,x\right)_{_{0\leq j,k\leq N-1}}} I + \right. \nonumber\\ & \left. M^{-1}\cdot\left(\int\left[ f\,\psi_{j}\,\epsilon\,\psi_{k} - f\,\psi_{k}\,\epsilon\,\psi_{j} -  f\,\psi_{k}\,\epsilon\,(f\,\psi_{j}) \right]\,d\,x\right)_{_{0\leq j,k\leq N-1}} \right)_{_{0\leq j,k\leq N-1}}
\end{align}
where the dot denotes matrix multiplication of $M^{-1}$ and the matrix with the integral as its $(j, k)$--entry. define $\eta_{j}=\sum_{j}\mu_{j\,k}\,\psi_{k}$ and use it to simplify the result of carrying out the matrix multiplication. From \eqref{eq:18} it follows that $\det(M)$ depends only on $N$ we lump it into $C^{(N)}_{1}$. Thus $G_{1,N}^{2}(t,\lambda)$ equals
\begin{equation}
  C^{(N)}_{1}\,\det\left(I + \left(\int\left[ f\,\eta_{j}\,\epsilon\,\psi_{k} - f\,\psi_{k}\,\epsilon\,\eta_{j} -  f\,\psi_{k}\,\epsilon\,(f\,\eta_{j}) \right]\,d\,x\right)_{_{0\leq j,k\leq N-1}} \right)_{_{0\leq j,k\leq N-1}}.
\end{equation}
Recall our remark at the very beginning of the section that if $f\equiv 0$ then the integral we started with evaluates to $1$ so that
\begin{equation}
  C^{(N)}_{1}\det(I) = C^{(N)}_{1},
\end{equation}
which implies that $C^{(N)}_{1}=1$. Now $G_{1,N}^{2}(t,\lambda)$ is of the form $\det(I+AB)$ where $A:L^{2}(J)\times L^{2}(J) \to \mathbb C^{N}$ is a $N\times 2$ matrix
\begin{equation*}
  A=\left(
    \begin{array}{c}
      A_{1} \\
      A_{2} \\
      \vdots \\
      A_{N}
    \end{array}
  \right),
\end{equation*}
whose $j^{th}$ row is given by
\begin{equation*}
  A_{j}= A_{j}(x)=\left(-f\,\epsilon\,\eta_{j}-f\,\epsilon\,\left(f\,\eta_{j}\right)\qquad f\,\eta_{j}\right).
\end{equation*}
Therefore, if
\begin{equation*}
  g=\left(
    \begin{array}{c}
      g_{1} \\
      g_{2}
    \end{array}
  \right)\in L^{2}(J)\times L^{2}(J),
\end{equation*}
then $A\,g$ is a column vector whose $j^{th}$ row is $\inprod{A_{j}}{g}_{_{L^{2}\times L^{2}}}$
\begin{equation*}
  \left(A\,g\right)_{j}=  \int\left[-f\,\epsilon\,\eta_{j}-f\,\epsilon\,\left(f\,\eta_{j}\right)\right]\,g_{1}\,d\,x + \int f\,\eta_{j}\,g_{2}\,d\,x.
\end{equation*}
Similarly, $B:\mathbb C^{N} \to L^{2}(J)\times L^{2}(J)$ is a $2\times N$ matrix
\begin{equation*}
  B=\left(
    \begin{array}{cccc}
      B_{1} & B_{2} & \ldots & B_{N}
    \end{array}\right),
\end{equation*}
whose $j^{th}$ column is given by
\begin{equation*}
  B_{j}=B_{j}(x)=\left(\begin{array}{c} \psi_{j} \\ \epsilon\,\psi_{j}
    \end{array}
  \right).
\end{equation*}
Thus if 
\begin{equation*}
  h=\left(
    \begin{array}{c}
      h_{1}  \\ 
      \vdots \\
      h_{N}
    \end{array}
  \right)\in \mathbb C^{N},
\end{equation*}
then $B\,h$ is the column vector of $L^{2}(J)\times L^{2}(J)$ given by
\begin{equation*}
  B\,h=\left(
    \begin{array}{c}
      \sum_{j}h_{i}\,\psi_{j}  \\
      \sum_{j}h_{i}\,\epsilon\,\psi_{j}
    \end{array}
  \right).
\end{equation*}
Clearly $A\,B:\mathbb C^{N} \to \mathbb C^{N}$ and $B\,A: L^{2}(J)\times L^{2}(J) \to L^{2}(J)\times L^{2}(J) $ with kernel
\begin{equation*}
  \left(
    \begin{array}{cc}
      -\sum_{j}\psi_{j}\otimes f\,\epsilon\,\eta_{j} - \sum_{j}\psi_{j}\otimes f\,\epsilon\,(f\,\eta_{j}) & \sum_{j}\psi_{j}\otimes f\,\eta_{j} \\
      -\sum_{j}\epsilon\,\psi_{j}\otimes f\,\epsilon\,\eta_{j} - \sum_{j}\epsilon\,\psi_{j}\otimes f\,\epsilon\,(f\,\eta_{j}) &
      \sum_{j}\epsilon\,\psi_{j}\otimes f\,\eta_{j}
    \end{array}
  \right).
\end{equation*}
Hence $I+BA$ has kernel
\begin{equation*}
  \left(
    \begin{array}{cc}
      I -\sum_{j}\psi_{j}\otimes f\,\epsilon\,\eta_{j} - \sum_{j}\psi_{j}\otimes f\,\epsilon\,(f\,\eta_{j}) &  \sum_{j}\psi_{j}\otimes f\,\eta_{j} \\
      -\sum_{j}\epsilon\,\psi_{j}\otimes f\,\epsilon\,\eta_{j} - \sum_{j}\epsilon\,\psi_{j}\otimes f\,\epsilon\,(f\,\eta_{j}) & I + \sum_{j}\epsilon\,\psi_{j}\otimes f\,\eta_{j}
    \end{array}
  \right),
\end{equation*}
which can be written as
\begin{equation*}
  \left(
    \begin{array}{cc}
      I -\sum_{j}\psi_{j}\otimes f\,\epsilon\,\eta_{j} &   \sum_{j}\psi_{j}\otimes f\,\eta_{j} \\
      -\sum_{j}\epsilon\,\psi_{j}\otimes f\,\epsilon\,\eta_{j} - \epsilon\,f & I + \sum_{j}\epsilon\,\psi_{j}\otimes f\,\eta_{j}
    \end{array}
  \right)\cdot \left(
    \begin{array}{cc}
      I & 0 \\
      \epsilon\,f & I
    \end{array}
  \right).
\end{equation*}
Since we are taking the determinant of this operator expression, and the determinant of the second term is just 1, we can drop it. Therefore
\begin{align*}
  G_{1,N}^{2}(t,\lambda) & =\det\left(
    \begin{array}{cc}
      I -\sum_{j}\psi_{j}\otimes f\,\epsilon\,\eta_{j} &  \sum_{j}\psi_{j}\otimes f\,\eta_{j} \\
      -\sum_{j}\epsilon\,\psi_{j}\otimes f\,\epsilon\,\eta_{j} - \epsilon\,f & I + \sum_{j}\epsilon\,\psi_{j}\otimes f\,\eta_{j}
    \end{array}
  \right) \\
  &= \det(I+\,K_{1,N}\,f),
\end{align*}
where
\begin{align*}
  K_{1,N}& = \left(
    \begin{array}{cc}
      -\sum_{j}\psi_{j}\otimes \,\epsilon\,\eta_{j}  &   \sum_{j}\psi_{j}\otimes \,\eta_{j} \\
      -\sum_{j}\epsilon\,\psi_{j}\otimes \,\epsilon\,\eta_{j} - \epsilon & \sum_{j}\epsilon\,\psi_{j}\otimes \,\eta_{j}
    \end{array}
  \right)\\
  &= \left(
    \begin{array}{cc}
      -\sum_{j,k}\psi_{j}\otimes \,\mu_{jk}\,\epsilon\,\psi_{k}  &   \sum_{j,k}\psi_{j}\otimes \,\mu_{jk}\,\psi_{k} \\
      -\sum_{j,k}\epsilon\,\psi_{j}\otimes \,\mu_{jk}\,\epsilon\,\psi_{k} - \epsilon & \sum_{j,k}\epsilon\,\psi_{j}\otimes \,\mu_{jk}\,\psi_{k}
    \end{array}
  \right)
\end{align*}
and $K_{1,N}$ has matrix kernel
\begin{align*}
  K_{1,N}(x,y)&= \left(
    \begin{array}{cc}
      -\sum_{j,k}\psi_{j}(x)\,\mu_{jk}\,\epsilon\,\psi_{k}(y) &   \sum_{j,k}\psi_{j}(x)\,\mu_{jk}\,\psi_{k}(y) \\
      -\sum_{j,k}\epsilon\,\psi_{j}(x) \,\mu_{jk}\,\epsilon\,\psi_{k}(y) - \epsilon(x-y) & \sum_{j,k}\epsilon\,\psi_{j}(x)\,\mu_{jk}\,\psi_{k}(y)
    \end{array}
  \right).
\end{align*}
We define
\begin{equation*}
  S_{N}(x,y)=- \sum_{j,k=0}^{N-1}\psi_{j}(x)\mu_{jk}\epsilon\,\psi_{k}(y).
\end{equation*}
Since $M$ is antisymmetric,
\begin{equation*}
  S_{N}(y,x)=- \sum_{j,k=0}^{N-1}\psi_{j}(y)\mu_{jk}\epsilon\,\psi_{k}(x)=\sum_{j,k=0}^{N-1}\psi_{j}(y)\mu_{kj}\epsilon\,\psi_{k}(x)= \sum_{j,k=0}^{N-1}\epsilon\,\psi_{j}(x)\mu_{jk}\psi_{k}(y).
\end{equation*}
Note that
\begin{equation*}
  \epsilon\,S_{N}(x,y)= \sum_{j,k=0}^{N-1}\epsilon\,\psi_{j}(x)\,\mu_{jk}\,\epsilon\,\psi_{k}(y),
\end{equation*}
whereas
\begin{equation*}
  -\frac{d}{d\,y}\,S_{N}(x,y)=\sum_{j,k=0}^{N-1}\psi_{j}(x)\mu_{jk}\psi_{k}(y).
\end{equation*}
So we can now write succinctly
\begin{equation}
  \label{eq:19}
  K_{1,N}(x,y)=\left(
    \begin{array}{cc}
      S_{N}(x,y) & -\frac{d}{d\,y}\,S_{N}(x,y) \\
      \epsilon\,S_{N}(x,y) - \epsilon & S_{N}(y,x)
    \end{array}
  \right)
\end{equation}
So we have shown that
\begin{equation}
\label{eq:21}
   G_{1,N}(t,\lambda)=\sqrt{D_{1,N}\left(t,\lambda\right)}
\end{equation}
where
\begin{equation*}
   D_{1,N}(t,\lambda)=\det\left(I+K_{1,N}\,f\right)
\end{equation*}
where $K_{1,N}$ is the integral operator with matrix kernel $K_{1,N}(x,y)$ given in \eqref{eq:19}.

\subsection{Gaussian specialization}
\label{sec:gauss-spec}

We specialize the results above to the case of a Gaussian weight function
\begin{equation}
  \label{eq:36}
  w(x)=\exp\left(-x^{2}/2\right)
\end{equation}
and indicator function
\begin{equation*}
  f(x)=-\lambda\,\rchi_{_{J}}\qquad J=(t,\infty)
\end{equation*}
Note that this does not agree with the weight function in \eqref{jointdensity}. However it is a necessary choice if we want the technical convenience of working with exactly the same orthogonal polynomials (the Hermite functions) as in the $\beta=2,4$ cases. In turn the Painlev\'e function in the limiting distribution will be unchanged. The discrepancy is resolved by the choice of standard deviation. Namely here the standard deviation on the diagonal matrix elements is taken to be $1$, corresponding to the weight function \eqref{eq:36}. In the $\beta=2,4$ cases the standard deviation on the diagonal matrix elements is $1/\sqrt{2}$, giving the weight function \eqref{gseweight}.
Now we again want the matrix
\begin{equation*}
  M=\left(\int\int\epsilon(x-y)\,\psi_{j}(x)\,\psi_{k}(y)\,d\,x\,d\,y \right)_{_{0\leq j,k\leq N-1}} = \left(\int\,\psi_{j}(x)\,\epsilon\,\psi_{k}(x)\,d\,x \right)_{_{0\leq j,k\leq N-1}}
\end{equation*}
to be the direct sum of $\frac{N}{2}$ copies of
\begin{equation*}
  Z=\left(\begin{array}{cc} 0 & 1 \\ -1 & 0 \end{array}\right)
\end{equation*}
so that the formulas are the simplest possible, since then $\mu_{jk}$ can only be $0$ or $\pm 1$.  In that case $M$ would be skew--symmetric so that $M^{-1}=-M$. In terms of the integrals defining the entries of $M$ this means that we would like to have
\begin{equation*}
  \int\psi_{2m}(x)\,\epsilon\,\psi_{2n+1}(x)\,d\,x = \delta_{m,n},
\end{equation*}
\begin{equation*}
  \int \psi_{2m+1}(x)\,\epsilon\,\psi_{2n}(x)\,d\,x = -\delta_{m,n},
\end{equation*}
and otherwise
\begin{equation*}
  \int \psi_{j}(x)\,\frac{d}{d\,x}\,\psi_{k}(x) \,d\,x = 0.
\end{equation*}
It is easier to treat this last case if we replace it with three non-exclusive conditions
\begin{equation*}
  \int \psi_{2m}(x)\,\epsilon\,\psi_{2n}(x)d\,x = 0,
\end{equation*}
\begin{equation*}
  \int \psi_{2m+1}(x)\,\epsilon\,\psi_{2n+1}(x)\,d\,x = 0
\end{equation*}
(so when the parity is the same for  $j,k$, which takes care of diagonal entries, among others), and
\begin{equation*}
  \int \psi_{j}(x)\,\epsilon\,\psi_{k}(x)\,d\,x = 0.
\end{equation*}
whenever $|j-k|>1$, which targets entries outside of the tridiagonal. Define
\begin{equation*}
  \varphi_{n}(x)=\frac{1}{c_{n}}H_{n}(x)\,\exp(-x^{2}/2)\quad\textrm{for}\quad c_{n}=\sqrt{2^{n}n!\sqrt{\pi}}
\end{equation*}
where the $H_{n}$ are the usual Hermite polynomials defined by the orthogonality condition
\begin{equation*}
  \int H_{j}(x)\,H_{k}(x)\,e^{-x^{2}}\,d\,x = c_{j}^{2}\,\delta_{j,k}.
\end{equation*}
It follows that
\begin{equation*}
  \int \varphi_{j}(x)\,\varphi_{k}(x)\,d\,x = \delta_{j,k}.
\end{equation*}
Now let
\begin{equation}
  \label{eq:20}
  \psi_{_{2n+1}}(x) = \frac{d}{d\,x}\,\varphi_{_{2n}}\left(\,x\right)\qquad
  \psi_{_{2n}}(x) = \varphi_{_{2n}}\left(\,x\right).
\end{equation}
This definition satisfies our earlier requirement that $\psi_{j}=p_{j}\,w$ for 
\begin{equation*}
  w(x)=\exp\left(-x^{2}/2\right).
\end{equation*}
In this case for example
\begin{equation*}
  p_{2n}(x)=\frac{1}{c_{n}}\,H_{2n}\left(x\right).
\end{equation*}
With $\epsilon$ defined as in \eqref{epsilonop}, and recalling that, if $D$ denote the operator that acts by differentiation with respect to $x$, then $D\,\epsilon=\epsilon\,D=I$, it follows that
\begin{equation*}
  \begin{aligned}
    \int \psi_{_{2m}}(x)\epsilon\,\psi_{_{2n+1}}(x)d\,x & =\int \varphi_{_{2m}}\left(x\right)\epsilon\,\frac{d}{d\,x}\, \varphi_{_{2n+1}}\left(x\right)d\,x \\
    &=\int \varphi_{_{2m}}\left(x\right) \varphi_{_{2n+1}}\left(x\right)d\,x \\
    &= \int \varphi_{_{2m}}\left(x\right) \varphi_{_{2n+1}}\left(x\right)d\left(x\right) \\ 
    &=\delta_{m,n},
  \end{aligned}
\end{equation*}
as desired. Similarly, integration by parts gives
\begin{equation*}
  \begin{aligned}
    \int \psi_{_{2m+1}}(x)\epsilon\,\psi_{_{2n}}(x)d\,x & =\int \frac{d}{d\,x}\,\varphi_{_{2m}}\left(x\right)\,\epsilon\,\varphi_{_{2n}} \left(x\right)d\,x \\
    &= -\int \,\varphi_{_{2m}}\left(x\right)\,\varphi_{_{2n}} \left(x\right)d\,x \\
    &=-\int \varphi_{_{2m}}\left(x\right) \varphi_{_{2n+1}}\left(x\right)d\left(x\right) \\& = -\delta_{m,n}.
  \end{aligned}
\end{equation*}
Also $\psi_{_{2n}}$ is even since $H_{2n}$ and $\varphi_{_{2n}}$ are. Similarly, $\psi_{_{2n+1}}$ is odd.  It follows  that $\epsilon\,\psi_{_{2n}}$, and $\epsilon\,\psi_{_{2n+1}}$, are respectively odd and even functions. From these observations, we obtain
\begin{equation*}
  \begin{aligned}\int \psi_{_{2n}}(x)\,\epsilon\,\psi_{_{2m}}(x)d\,x  = 0,
  \end{aligned}
\end{equation*}
since the integrand is a product of an odd and an even function. Similarly
\begin{equation*}
  \begin{aligned}
    \int \psi_{_{2n+1}}(x)\,\epsilon\,\psi_{_{2m+1}}(x)d\,x= 0.
  \end{aligned}
\end{equation*}
Finally it is easy to see that if $|j-k|>1$, then
\begin{equation*}
  \int\psi_{_{j}}(x)\,\epsilon\,\psi_{_{k}}(x)d\,x=0.
\end{equation*}
Indeed both  differentiation and the action of $\epsilon$ can only ``shift'' the indices by $1$. Thus by orthogonality of the $\varphi_{j}$, this integral will always be $0$. Thus by our choice in \eqref{eq:20}, we force the matrix
\begin{equation*}
M= \left(\int\,\psi_{j}(x)\,\epsilon\,\psi_{k}(x)\,d\,x \right)_{_{0\leq j,k\leq N-1}}
\end{equation*}
to be the direct sum of $\frac{N}{2}$ copies of
\begin{equation*}
  Z=\left(
    \begin{array}{cc}
      0 & 1 \\
      -1 & 0 
    \end{array}
  \right)
\end{equation*}
This means $M^{-1}=-M$ where $M^{-1}=\{\mu_{j,k}\}$. Moreover, $\mu_{j,k}=0$ if $j,k$ have the same parity or $|j-k|>1$, and $\mu_{2j,2k+1}=\delta_{jk}=-\mu_{2k+1,2j}$ for $j,k=0.\ldots,\frac{N}{2}-1$. Therefore
\begin{equation*}
  \begin{aligned}S_{N}(x,y) & =- \sum_{j,k=0}^{N-1}\psi_{j}(x)\mu_{jk}\epsilon\,\psi_{k}(y)\\
    & = -\sum_{j=0}^{N/2-1}\,\psi_{2j}(x)\,\epsilon\,\psi_{2j+1}(y) + \sum_{j=0}^{N/2-1}\,\psi_{2j+1}(x)\,\epsilon\,\psi_{2j}(y) \\
    & =  \left[ \sum_{j=0}^{N/2-1}\varphi_{_{2j}}\left(\frac{x}{\,\sigma}\right) \,\varphi_{_{2j}}\left(\frac{y}{\,\sigma}\right)-\sum_{j=0}^{N/2-1}\frac{d}{d\,x}\, \varphi_{_{2j}}\left(\frac{x}{\,\sigma}\right)\,\epsilon\,\varphi_{_{2j}} \left(\frac{y}{\,\sigma}\right) \right].
  \end{aligned}
\end{equation*}
Manipulations similar to those in the $\beta=4$ case (see \eqref{hermrec1} through \eqref{newgseS}) yield
\begin{equation*}
  \begin{aligned}
    S_{N}(x,y) =  \left[ \sum_{j=0}^{N-1}\varphi_{_{j}}\left(x\right) \varphi_{_{j}}\left(y\right)-\sqrt{\frac{N}{2}}\, \varphi_{_{N-1}}\left(x\right)\left(\epsilon\,\varphi_{_{N}}\right) \left(y\right) \right].
  \end{aligned}
\end{equation*}
We redefine
\begin{equation*}
  S_{N}(x,y) = \sum_{j=0}^{N-1}\varphi_{_{j}}\left(x\right) \varphi_{_{j}}\left(y\right)=S_{N}(y,x),
\end{equation*}
so that the top left entry of $K_{1,N}(x,y)$ is
\begin{equation*}
  S_{N}(x,y) +  \sqrt{\frac{N}{2}}\, \varphi_{_{N-1}}\left(x\right)\left(\epsilon\,\varphi_{_{N}}\right) \left(y\right).
\end{equation*}
If $S_{N}$ is the operator with kernel $S_{N}(x,y)$ then integration by parts gives
\begin{equation*}
  S_{N}D f = \int S(x,y) \frac{d}{d\,y}f(y)\,d\,y = \int  \left(-\frac{d}{d\,y} S_{N}(x,y)\right) f(y)\,d\,y,
\end{equation*}
so that $-\frac{d}{d\,y} S_{N}(x,y)$ is in fact the kernel of $S_{N}D$. Therefore \eqref{eq:21} now holds with $K_{1,N}$ being the integral operator with matrix kernel $K_{1,N}(x,y)$ whose $(i,j)$--entry $K_{1,N}^{(i,j)}(x,y)$ is given by
\begin{equation*}
  \begin{aligned} 
    K _{1,N}^{(1,1)}(x,y) & =\left[ S_{N}(x,y) + \sqrt{\frac{N}{2}}\, \varphi_{_{N-1}}\left(x\right)\left(\epsilon\,\varphi_{_{N}}\right) \left(y\right)\right], \\
    K_{1,N}^{(1,2)}(x,y) & =\left[ SD_{N}(x,y) -\frac{d}{d\,y} \, \left(\sqrt{\frac{N}{2}}\, \varphi_{_{N-1}}\left(x\right)\left(\epsilon\,\varphi_{_{N}}\right) \left(y\right) \right) \right], \\
    K_{1,N}^{(2,1)}(x,y) & = \epsilon\,\left[ S_{N}(x,y) +  \sqrt{\frac{N}{2}}\, \varphi_{_{N-1}}\left(x\right)\left(\epsilon\,\varphi_{_{N}}\right)\left(y\right)-1 \right], \\
    K_{1,N}^{(2,2)}(x,y) & =\left[ S_{N}(x,y) + \sqrt{\frac{N}{2}}\,\left(\epsilon\,\varphi_{_{N}}\right)\left(x\right)\, \varphi_{_{N-1}} \left(y\right)\right].
  \end{aligned}
\end{equation*}
Define
\begin{equation*}
  \varphi(x)=\left(\frac{N}{2}\right)^{1/4}\varphi_{_{N}}(x),\qquad \psi(x)=\left(\frac{N}{2}\right)^{1/4}\varphi_{_{N-1}}(x),
\end{equation*}
so that
\begin{equation*}
  \begin{aligned}
    K_{1,N}^{(1,1)}(x,y)& =\rchi(x)\,\left[ S_{N}(x,y) + \psi(x)\,\epsilon\,\varphi(y)\right] \,\rchi(y),\\
    K_{1,N}^{(1,2)}(x,y) & = \rchi(x)\,\left[ SD_{N}(x,y) - \psi(x)\,\varphi(y)  \right]\,\rchi(y), \\
    K_{1,N}^{(2,1)}(x,y) & = \rchi(x)\,\left[ \epsilon S_{N}(x,y) + \epsilon\,\psi(x)\,\epsilon\,\varphi(y) - \epsilon(x-y) \right]\,\rchi(y), \\
    K_{1,N}^{(2,2)}(x,y) & = \rchi(x)\,\left[ S_{N}(x,y) + \epsilon\,\varphi(x)\,\psi(y)\right]\,\rchi(y).
  \end{aligned}
\end{equation*}
Note that
\begin{equation*}
  \begin{aligned}
    \rchi\,\left( S + \psi\otimes\,\epsilon\,\varphi \right)\,\rchi &\doteq K_{1,N}^{(1,1)}(x,y),   \\
    \rchi\,\left( SD - \psi\otimes\,\varphi  \right)\,\rchi &\doteq K_{1,N}^{(1,2)}(x,y),  \\
    \rchi\,\left( \epsilon\,S + \epsilon\,\psi\otimes\,\epsilon\,\varphi - \epsilon \right)\,\rchi &\doteq K_{1,N}^{(2,1)}(x,y) , \\
    \rchi\,\left( S + \epsilon\,\varphi\otimes\,\epsilon\,\psi \right)\,\rchi &\doteq   K_{1,N}^{(2,2)}(x,y). \\
  \end{aligned}
\end{equation*}
Hence
\begin{equation*}
  K_{1,N}=\rchi\,\left(
    \begin{array}{cc}
      S + \psi\otimes\,\epsilon\,\varphi  & SD - \psi\otimes\,\varphi \\
      \epsilon\,S + \epsilon\,\psi\otimes\,\epsilon\,\varphi  - \epsilon &  S + \epsilon\,\varphi\otimes\,\epsilon\,\psi 
    \end{array}
  \right)\,\rchi.
\end{equation*}

Note that this is identical to the corresponding operator for $\beta=1$ obtained by Tracy and Widom in \cite{Trac2}, the only difference being that $\varphi$, $\psi$, and hence also $S$, are redefined to depend on $\lambda$. 

\subsection{Edge scaling}

\subsubsection{Reduction of the determinant}

The above determinant is that of an operator on $L^{2}(J)\oplus L^{2}(J)$. Our first task will be to rewrite these determinants as those of operators on $L^{2}(J)$. This part follows exactly the proof in \cite{Trac2}. To begin, note that \begin{equation}\com{S}{D}=\varphi\otimes\psi + \psi\otimes\varphi\label{sdcom}\end{equation} so that (using the fact that $D\,\epsilon=\epsilon\,D=I$ )
\begin{eqnarray}
  \com{\epsilon}{S} &=&\epsilon\,S-S\,\epsilon\nonumber\\
  &=& \epsilon\,S\,D\,\epsilon-\epsilon\,D\,S\,\epsilon = \epsilon\,\com{S}{D}\,\epsilon\nonumber\\
  &=& \epsilon\,\varphi\otimes\psi\,\epsilon + \epsilon\,\psi\otimes\varphi\,\epsilon\nonumber\\
  &=&\epsilon\,\varphi\otimes\epsilon^{t}\psi + \epsilon\,\psi\otimes\epsilon^{t}\,\varphi\nonumber\\
  &=& - \epsilon\,\varphi\otimes\epsilon\,\psi -
  \epsilon\,\psi\otimes\epsilon\,\varphi,
  \label{escom}\end{eqnarray} 
where the last equality follows from the fact that $\epsilon^{t}=-\epsilon$. We thus have

\begin{eqnarray*}
  D\,\left(\epsilon\,S + \epsilon\,\psi\otimes\epsilon\,\varphi\right) & = & S + \psi\otimes\epsilon\,\varphi,\\
  D\,\left(\epsilon\,S\,D - \epsilon\,\psi\otimes\varphi\right) & =
  & S\,D - \psi\otimes\varphi.
\end{eqnarray*}
The expressions on the right side are the top entries of $K_{1,N}$. Thus the first row of $K_{1,N}$ is, as a vector, 
\begin{equation*}
  D\,\left(\epsilon\,S + \epsilon\,\psi\otimes\epsilon\,\varphi, \epsilon\,S\,D - \epsilon\,\psi\otimes\varphi\right).
\end{equation*} 
Now \eqref{escom} implies that
\begin{equation*} 
\epsilon\,S + \epsilon\,\psi\otimes\epsilon\,\varphi = S\,\epsilon -\epsilon\,\varphi\otimes\epsilon\,\psi.
\end{equation*}
Similarly \eqref{sdcom} gives
\begin{equation*}
  \epsilon\,\com{S}{D} = \epsilon\,\varphi\otimes\psi + \epsilon\,\psi\otimes\varphi,
\end{equation*} 
so that
\begin{equation*}
  \epsilon\,S\,D - \epsilon\,\psi\otimes\varphi = \epsilon\,D\,S + \epsilon\,\varphi\otimes\psi = S + \epsilon\,\varphi\otimes\psi.
\end{equation*}
Using these expressions we can rewrite the first row of $K_{1,N}$ as
\begin{equation*}
D\,\left(S\,\epsilon - \epsilon\,\varphi\otimes\epsilon\,\psi, S + \epsilon\,\varphi\otimes\psi\right).
\end{equation*}
Applying $\epsilon$ to this expression shows the second row of $K_{1,N}$ is given by
\begin{equation*}
\left(\epsilon\,S - \epsilon + \epsilon\,\psi\otimes\epsilon\,\varphi, S + \epsilon\,\varphi\otimes\psi\right) 
\end{equation*}
Now use
\eqref{escom} to show the second row of $K_{1,N}$ is
\begin{equation*}
\left(S\,\epsilon - \epsilon + \epsilon\,\varphi\otimes\epsilon\,\psi, S + \epsilon\,\varphi\otimes\psi\right).
\end{equation*} 
Therefore,
\begin{eqnarray*}
  K_{1,N} &=& \rchi\,
  \left(
    \begin{array}{cc} D\,\left(S\,\epsilon - \epsilon\,\varphi\otimes\epsilon\,\psi\right) & D\,\left(S + \epsilon\,\varphi\otimes\psi\right) \\
      S\,\epsilon - \epsilon + \epsilon\,\varphi\otimes\epsilon\,\psi &  S + \epsilon\,\varphi\otimes\psi
    \end{array}
  \right)
  \,\rchi \\
  & =& 
  \left(
    \begin{array}{cc}
      \rchi\,D & 0 \\
      0 & \rchi
    \end{array}
  \right)
  \left(
    \begin{array}{cc}
      \left(S\,\epsilon - \epsilon\,\varphi\otimes\epsilon\,\psi\right)\,\rchi & \left(S + \epsilon\,\varphi\otimes\psi\right)\,\rchi \\
      \left(S\,\epsilon - \epsilon + \epsilon\,\varphi\otimes\epsilon\,\psi\right)\rchi & \left(S + \epsilon\,\varphi\otimes\psi\right)\rchi
    \end{array}
  \right).
\end{eqnarray*}
Since $K_{1,N}$ is of the form $A\,B$, we can use the fact that $\det(I-A\,B)=~\det(I-~B\,A)$ and deduce that $D_{1,N}(s,\lambda)$ is unchanged if instead we take $K_{1,N}$ to be
\begin{eqnarray*}
  K_{1,N} & = & 
  \left(
    \begin{array}{cc}
      \left(S\,\epsilon - \epsilon\,\varphi\otimes\epsilon\,\psi\right)\,\rchi & \left(S + \epsilon\,\varphi\otimes\psi\right)\,\rchi \\
      \left(S\,\epsilon - \epsilon + \epsilon\,\varphi\otimes\epsilon\,\psi\right)\,\rchi & \left(S + \epsilon\,\varphi\otimes\psi\right)\,\rchi
    \end{array}
  \right)
  \left(
    \begin{array}{cc}
      \rchi\,D & 0 \\
      0 & \rchi
    \end{array}
  \right)\\
  & = & 
  \left(
    \begin{array}{cc}
      \left(S\,\epsilon - \epsilon\,\varphi\otimes\epsilon\,\psi\right)\,\rchi\,D & \left(S + \epsilon\,\varphi\otimes\psi\right)\,\rchi \\
      \left(S\,\epsilon - \epsilon + \epsilon\,\varphi\otimes\epsilon\,\psi\right)\rchi\,D & \left(S + \epsilon\,\varphi\otimes\psi\right)\rchi
    \end{array}\right).
\end{eqnarray*}
Therefore
\begin{equation}
  D_{1,N}(s,\lambda)=\det
  \left(
    \begin{array}{cc}
      I - \left(S\,\epsilon - \epsilon\,\varphi\otimes\epsilon\,\psi\right)\,\lambda\,\rchi\,D & - \left(S + \epsilon\,\varphi\otimes\psi\right)\,\lambda\,\rchi \\
      - \left(S\,\epsilon - \epsilon + \epsilon\,\varphi\otimes\epsilon\,\psi\right)\,\lambda\,\rchi\,D & I - \left(S + \epsilon\,\varphi\otimes\psi\right)\,\lambda\,\rchi
    \end{array}
  \right).
\end{equation}
Now we perform row and column operations on the matrix to simplify it, which do not change the Fredholm determinant. Justification of these operations is given in \cite{Trac2}. We start by subtracting row 1 from row 2 to get

\begin{equation*}  \left(\begin{array}{cc}
I - \left(S\,\epsilon - \epsilon\,\varphi\otimes\epsilon\,\psi\right)\,\lambda\,\rchi\,D & - \left(S + \epsilon\,\varphi\otimes\psi\right)\,\lambda\,\rchi \\
  - I + \epsilon\,\lambda\,\rchi\,D  & I \end{array}\right). \end{equation*}
Next, adding column 2 to column 1 yields
\begin{equation*} 
  \left(
    \begin{array}{cc}
      I - \left(S\,\epsilon - \epsilon\,\varphi\otimes\epsilon\,\psi\right)\,\lambda\,\rchi\,D - \left(S + \epsilon\,\varphi\otimes\psi\right)\,\lambda\,\rchi & - \left(S + \epsilon\,\varphi\otimes\psi\right)\,\lambda\,\rchi \\
      \epsilon\,\lambda\,\rchi\,D  & I 
    \end{array} 
  \right). 
\end{equation*}
Then right-multiply column 2 by $-\epsilon\,\lambda\,\rchi\,D$ and add it to column 1, and multiply row 2 by $S + \epsilon\,\varphi\otimes\psi$ and add it to row 1 to arrive at
\begin{equation*} 
  \det
  \left(
    \begin{array}{cc} 
      I - \left(S\,\epsilon - \epsilon\,\varphi\otimes\epsilon\,\psi\right)\,\lambda\,\rchi\,D +  \left(S + \epsilon\,\varphi\otimes\psi\right)\,\lambda\,\rchi\,\left(\epsilon\,\lambda\,\rchi\,D - I\right) & 0 \\
      0  & I 
    \end{array} 
  \right). 
\end{equation*}
Thus the determinant we want equals the determinant of
\begin{equation}
  I - \left(S\,\epsilon - \epsilon\,\varphi\otimes\epsilon\,\psi\right)\,\lambda\,\rchi\,D + \left(S +   \epsilon\,\varphi\otimes\psi\right)\,\lambda\,\rchi\,\left(\epsilon\,\lambda\,\rchi\,D - I\right) 
  \label{operator}.
\end{equation}
So we have reduced the problem from the computation of the Fredholm determinant of an operator on $L^{2}(J)~\times~L^{2}(J)$, to that of an operator on $L^{2}(J)$.

\subsubsection{Differential equations}

Next we want to write the operator in \eqref{operator} in the form
\begin{equation} \left(I- K_{2, N}\right)\left(I -
\sum_{i=1}^{L}\alpha_{i}\otimes\beta_{i}\right), \end{equation}
where the $\alpha_{i}$  and $\beta_{i}$ are functions in
$L^{2}(J)$. In other words, we want to rewrite the determinant for
the GOE case as a finite dimensional perturbation of the
corresponding GUE determinant. The Fredholm determinant of the
product is then the product of the determinants. The limiting form
for the GUE part is already known, and we can just focus on
finding a limiting form for the determinant of the finite
dimensional piece. It is here that the proof must be modified from
that in \cite{Trac2}. A little simplification of \eqref{operator}
yields
\begin{equation*}
I -
\lambda\,S\,\rchi-\lambda\,S\,\left(1-\lambda\,\rchi\right)\,\epsilon\,\rchi\,D
- \lambda\,\left(\epsilon\,\varphi\,\otimes\,\rchi\,\psi\right) -
\lambda\,\left(\epsilon\,\varphi\,\otimes\,\psi\right)\left(1-\lambda\,\rchi\right)\,\epsilon\,\rchi\,D.
\end{equation*}
Writing $\epsilon\,\com{\rchi}{D}+\rchi$ for $\epsilon\,\rchi\,D$
and simplifying $\left(1-\lambda\,\rchi\right)\,\rchi$ to
$\left(1-\lambda\right)\,\rchi$ gives

 \begin{align*}& I -
\lambda\,S\,\rchi - \lambda\,\left(1-\lambda\right)\,S\,\rchi -
\lambda\,\left(\epsilon\,\varphi\,\otimes\,\rchi\,\psi\right)
-\lambda\,\left(1-\lambda\right)\,\left(\epsilon\,\varphi\,\otimes\,\rchi\,\psi\right)
\\ & \qquad
-\lambda\,S\,\left(1-\lambda\,\rchi\right)\,\epsilon\,\com{\rchi}{D}
-
\lambda\,\left(\epsilon\,\varphi\,\otimes\,\psi\right)\,\left(1-\lambda\,\rchi\right)\,\epsilon\,\com{\rchi}{D}
 \\ & = I -  (2\lambda-\lambda^{2})\,S\,\rchi - (2\lambda-\lambda^{2})\,(\epsilon\,\varphi\,\otimes\,\rchi\,\psi)
  - \lambda\,S\,(1-\lambda\,\rchi)\,\epsilon\,\com{\rchi}{D} \\ &
  \qquad
  -  \lambda\,(\epsilon\,\varphi\,\otimes\,\psi)\,(1-\lambda\,\rchi)\,\epsilon\,\com{\rchi}{D}.\end{align*}
Define $\tilde{\lambda}=2\,\lambda-\lambda^{2}$ and let
$\sqrt{\tilde{\lambda}}\,\varphi\to\varphi$, and
$\sqrt{\tilde{\lambda}}\,\psi\to\psi$ so that
$\tilde{\lambda}\,S\to S$ and \eqref{operator} goes to
\begin{align*} I -  & S\,\rchi - (\epsilon\,\varphi\,\otimes\,\rchi\,\psi)
  - \frac{\lambda}{\tilde{\lambda}}\,S\,(1-\lambda\,\rchi)\,\epsilon\,\com{\rchi}{D} \\ & -
 \frac{\lambda}{\tilde{\lambda}}\,(\epsilon\,\varphi\,\otimes\,\psi)\,(1-\lambda\,\rchi)\,\epsilon\,\com{\rchi}{D}.\end{align*}
Now we define $R:=(I-S\,\rchi)^{-1}\,S\,\rchi=(I-S\,\rchi)^{-1}-I$
(the resolvent operator of $S\,\rchi$), whose kernel we denote by
$R(x,y)$, and
$Q_{\epsilon}:=(I-S\,\rchi)^{-1}\,\epsilon\,\varphi$. Then
\eqref{operator} factors into

\begin{equation*}A =(I - S\,\rchi)\,B.\end{equation*}
where $B$ is

\begin{align*}I - & (Q_{\epsilon}\,\otimes\,\rchi\,\psi)
  - \frac{\lambda}{\tilde{\lambda}}\,(I+R)\,S\,(1-\lambda\,\rchi)\,\epsilon\,\com{\rchi}{D}\\ & -
  \frac{\lambda}{\tilde{\lambda}}\,(Q_{\epsilon}\,\otimes\,\psi)\,(1-\lambda\,\rchi)\,\epsilon\,\com{\rchi}{D},\qquad \lambda\neq 1.\end{align*} Hence
\begin{equation*}D_{1,N}(s,\lambda)
=D_{2,N}(s,\tilde{\lambda})\,\det(B).\end{equation*} 
Note that because of the change of variable $\tilde{\lambda}\,S\to S$, we are in effect factoring $I-(2\lambda-\lambda^{2})\,S$, rather that $I-\lambda\,S$ as we did in the $\beta=4$ case. The fact that we factored $I -  (2\lambda-\lambda^{2})\,S\,\rchi$ as opposed to $I - \lambda\,S\,\rchi$ is crucial here for it is what makes $B$ finite rank. If we had factored $I - \lambda\,S\,\rchi$ instead, $B$ would have been 
\[\begin{split} B =  I & -
\lambda\,\sum_{k=1}^{2}(-1)^{k}\left(S + R\,S\right)\,\left(I-\lambda\,\rchi\right)\epsilon_{k}\otimes\delta_{k} - \lambda\,\left(I + R\right)\,\epsilon\,\varphi\,\otimes\,\rchi\,\psi \\
& - \lambda\,\sum_{k=1}^{2}(-1)^{k}\inprod{\psi}{\left(I -\lambda\,\rchi\right)\,\epsilon_{k}}\,\left(\left(I + R\right)\,\epsilon\,\varphi\right)\otimes\delta_{k} \\ & - \lambda\,(1-\lambda)\,\left(S + R\,S\right)\,\rchi - \lambda\,(1-\lambda)\,\left(\left(I + R\right)\,\epsilon\,\varphi\right)\otimes\rchi\,\psi
\end{split}\]
The first term on the last line is not finite rank, and the methods we have used previously in the $\beta=4$ case would not work here. It is also interesting to note that these complications disappear when we are dealing with the case of the largest eigenvalue; then is no differentiation with respect to $\lambda$, and we just set $\lambda=1$ in all these formulae. All the new troublesome terms vanish!
\par\noindent In order to find $\det(B)$ we use the identity
\begin{equation}\epsilon\,\com{\rchi}{D}=\sum_{k=1}^{2m}
(-1)^{k}\,\epsilon_{k}\otimes\delta_{k},\end{equation} where
$\epsilon_{k}$ and $\delta_{k}$ are the functions
$\epsilon(x-a_{k})$ and $\delta(x-a_{k})$ respectively, and the
$a_{k}$ are the endpoints of the (disjoint) intervals considered,
$J=\cup_{k=1}^{m}(a_{2\,k-1},a_{2\,k})$. We also make use of the
fact that
\begin{equation}a\otimes b\cdot c\otimes d= \inprod{b}{c}\cdot a\otimes d\end{equation}
where $\inprod{.}{.}$ is the usual $L^{2}$--inner product.
Therefore

\begin{align*}(Q_{\epsilon}\otimes\psi)\,(1-\lambda\,\rchi)\,\epsilon\,\com{\rchi}{D} &=
\sum_{k=1}^{2m} (-1)^{k}Q_{\epsilon}\otimes\psi\cdot
(1-\lambda\,\rchi)\,\epsilon_{k}\otimes\delta_{k}
\\ &=\sum_{k=1}^{2m}
(-1)^{k}\inprod{\psi}{(1-\lambda\,\rchi)\,\epsilon_{k}}\,Q_{\epsilon}\otimes\,\delta_{k}.
\end{align*}
It follows that

\begin{equation*}\frac{D_{1,N}(s,\lambda)}{D_{2,N}(s,\tilde{\lambda})}\end{equation*}
equals the determinant of
\begin{equation} 
I  - Q_{\epsilon}\otimes\rchi\psi - \frac{\lambda}{\tilde{\lambda}}\,\sum_{k=1}^{2m} (-1)^{k}\left[(S+R\,S)\,(1-\lambda\,\rchi)\,\epsilon_{k} +\inprod{\psi}{(1-\lambda\,\rchi)\,\epsilon_{k}}\,Q_{\epsilon}\right]\otimes\delta_{k}.
\end{equation}
We now specialize to the case of one interval $J=(t,\infty)$, so
$m=1$, $a_{1}=t$ and $a_{2}=\infty$. We write
$\epsilon_{t}=\epsilon_{1}$, and $\epsilon_{\infty}=\epsilon_{2}$,
and similarly for $\delta_{k}$. Writing the terms in the summation
and using the facts that
\begin{equation}\epsilon_{\infty}=-\frac{1}{2},\end{equation}
 and
 \begin{equation}(1-\lambda\,\rchi)\,\epsilon_{t}=-\frac{1}{2}\,(1-\lambda\,\rchi)+(1-\lambda\,\rchi)\,\rchi,\end{equation}
 then yields
\begin{align*}
I - Q_{\epsilon}\otimes\rchi\psi  -
\frac{\lambda}{2\tilde{\lambda}}\,
\left[(S+R\,S)\,(1-\lambda\,\rchi)+
\inprod{\psi}{(1-\lambda\,\rchi)}\,Q_{\epsilon}\right]\otimes(\delta_{t}-\delta_{\infty})\\
\qquad   +\frac{\lambda}{\tilde{\lambda}}
\left[(S+R\,S)\,(1-\lambda\,\rchi)\,\rchi+
\inprod{\psi}{(1-\lambda\,\rchi)\,\rchi}\,Q_{\epsilon}\right]\otimes\delta_{t}
\end{align*}
which, to simplify notation, we write as
\begin{align*}
 I - Q_{\epsilon}\otimes\rchi\psi  - \frac{\lambda}{2\tilde{\lambda}}\,
\left[(S+R\,S)\,(1-\lambda\,\rchi)+
a_{1,\lambda}\,Q_{\epsilon}\right]\otimes(\delta_{t}-\delta_{\infty}) \\
\qquad   +\frac{\lambda}{\tilde{\lambda}}
\left[(S+R\,S)\,(1-\lambda\,\rchi)\,\rchi+
\tilde{a}_{1,\lambda}\,Q_{\epsilon}\right]\otimes\delta_{t},
\end{align*}
where
\begin{equation}a_{1,\lambda} = \inprod{\psi}{(1-\lambda\,\rchi)},\qquad \tilde{a}_{1,\lambda}= \inprod{\psi}{(1-\lambda\,\rchi)\,\rchi}.\end{equation}
 Now we can use the formula:
\begin{equation}\det\left(I-\sum_{i=1}^{L}\alpha_{i}\otimes\beta_{i}\right)=\det\left(\delta_{jk}-\inprod{\alpha_{j}}{\beta_{k}}\right)_{1\leq j,k\leq L}\end{equation}
In this case, $L=3$, and
\begin{align}\alpha_{1}=Q_{\epsilon}, &\qquad
 \alpha_{2} =\frac{\lambda}{\tilde{\lambda}}\,
\left[(S+R\,S)\,(1-\lambda\,\rchi)+ a_{1,\lambda}\,Q_{\epsilon}\right]\nonumber, \\
\alpha_{3}= & -\frac{\lambda}{\tilde{\lambda}}
\left[(S+R\,S)\,(1-\lambda\,\rchi)\,\rchi+
\tilde{a}_{1,\lambda}\,Q_{\epsilon}\right],\nonumber \\ &
\beta_{1}=\rchi\psi, \qquad  \beta_{2}=\delta_{t}-\delta_{\infty},
\qquad \beta_{3}=\delta_{t}.\end{align} In order to simplify the
notation, define

\begin{align}
Q(x,\lambda,t)&:=(I-S\,\rchi)^{-1}\,\varphi, &
P(x,\lambda,t)&:=(I-S\,\rchi)^{-1}\,\psi, \nonumber\\
Q_{\epsilon}(x,\lambda,t)&:=(I-S\,\rchi)^{-1}\,\epsilon\,\varphi,
& P_{\epsilon}(x,\lambda,t)&:=(I-S\,\rchi)^{-1}\,\epsilon\,\psi,
\end{align}

\begin{align}
q_{_{N}}&:=Q(t,\lambda,t), & p_{_{N}}&:=P(t,\lambda,t),\nonumber\\
q_{\epsilon}&:=Q_{\epsilon}(t,\lambda,t), & p_{\epsilon}&:=P_{\epsilon}(t,\lambda,t),\nonumber\\
u_{\epsilon}&:=\inprod{Q}{\rchi\,\epsilon\,\varphi}=\inprod{Q_{\epsilon}}{\rchi\,\varphi},
&
v_{\epsilon}&:=\inprod{Q}{\rchi\,\epsilon\,\psi}=\inprod{P_{\epsilon}}{\rchi\,\psi}, \nonumber\\
\tilde{v}_{\epsilon}&:=\inprod{P}{\rchi\,\epsilon\,\varphi}=\inprod{Q_{\epsilon}}{\rchi\,\varphi},
&
w_{\epsilon}&:=\inprod{P}{\rchi\,\epsilon\,\psi}=\inprod{P_{\epsilon}}{\rchi\,\psi},
\label{inproddef1}\end{align}

\begin{align}
\mathcal{P}_{1,\lambda}&:= \int(1-\lambda\,\rchi)\,P\,d\,x, &
\tilde{\mathcal{P}}_{1,\lambda} &:=
\int(1-\lambda\,\rchi)\,\rchi\,P\,d\,x, \nonumber\\
\mathcal{Q}_{1,\lambda} &:= \int(1-\lambda\,\rchi)\,Q\,d\,x, &
\tilde{\mathcal{Q}}_{1,\lambda} &:=
\int(1-\lambda\,\rchi)\,\rchi\,Q\,d\,x,\nonumber\\
\mathcal{R}_{1,\lambda}&:= \int(1-\lambda\,\rchi)\,R(x,t)\,d\,x, &
\tilde{\mathcal{R}}_{1,\lambda} &:=
\int(1-\lambda\,\rchi)\,\rchi\,R(x,t)\,d\,x.
\label{inproddef2}\end{align} Note that all quantities in
\eqref{inproddef1} and \eqref{inproddef2} are functions of $t$
alone. Furthermore, let
\begin{equation} c_{\varphi} =
\epsilon\,\varphi(\infty)=\frac{1}{2}\int_{-\infty}^{\infty}\varphi(x)\,d\,x,
\qquad c_{\psi} =
\epsilon\,\psi(\infty)=\frac{1}{2}\int_{-\infty}^{\infty}\psi(x)\,d\,x.\end{equation}
Recall from the previous section that when $\beta=1$ we take $N$ to be even. It follows that $\varphi$ and $\psi$ are even and odd functions respectively. Thus $c_{\psi}=0$ for $\beta=1$, and computation gives
\begin{equation}
c_{\varphi}=(\pi N)^{1/4} 2^{-3/4-N/2}\,\frac{(N!)^{1/2}}{(N/2)!}\,\sqrt{\lambda}.
\end{equation}
Hence computation yields
\begin{equation}
 \lim_{N\to\infty}c_{\varphi}= \sqrt{\frac{\lambda}{2}},
\end{equation}
and at $t=\infty$ we have
\begin{equation*}
  u_{\epsilon}(\infty)=0, \quad  q_{\epsilon}(\infty)=c_{\varphi}
\end{equation*}
\begin{equation*}
\mathcal{P}_{1,\lambda}(\infty) = 2\,c_{\psi}, \quad 
\mathcal{Q}_{1,\lambda}(\infty) =2\,c_{\varphi},\quad
\mathcal{R}_{1,\lambda}(\infty) = 0 ,
\end{equation*}
\begin{equation*}
\tilde{\mathcal{P}}_{1,\lambda}(\infty) = 
\tilde{\mathcal{Q}}_{1,\lambda}(\infty) =\tilde{
\mathcal{R}}_{1,\lambda}(\infty) = 0.
\end{equation*} 
Hence
\begin{align}
\inprod{\alpha_{1}}{\beta_{1}} & = \tilde{v}_{\epsilon}, \quad
\inprod{\alpha_{1}}{\beta_{2}}=q_{\epsilon}-c_{\varphi}, \quad
\inprod{\alpha_{1}}{\beta_{3}}=q_{\epsilon},
\\
\inprod{\alpha_{2}}{\beta_{1}} &
=\frac{\lambda}{2\,\tilde{\lambda}}
\,\left[\mathcal{P}_{1,\lambda}-a_{1,\lambda}\,(1-\tilde{v}_{\epsilon})\right], \\
\inprod{\alpha_{2}}{\beta_{2}} & =
\frac{\lambda}{2\,\tilde{\lambda}}\,\left[\mathcal{R}_{1,\lambda}
+ a_{1,\lambda}\,(q_{\epsilon}-c_{\varphi})\right]\label{example}, \\
\inprod{\alpha_{2}}{\beta_{3}} & =
\frac{\lambda}{2\,\tilde{\lambda}}\,\left[\mathcal{R}_{1,\lambda}
+ a_{1,\lambda}\,q_{\epsilon}\right], \\
\inprod{\alpha_{3}}{\beta_{1}} & =
-\frac{\lambda}{\tilde{\lambda}}\,\left[\tilde{\mathcal{P}}_{1,\lambda}
- \tilde{a}_{1,\lambda}\,(1-\tilde{v}_{\epsilon})\right], \\
\inprod{\alpha_{3}}{\beta_{2}} & =
-\frac{\lambda}{\tilde{\lambda}}\,\left[\tilde{\mathcal{R}}_{1,\lambda}+\tilde{a}_{1,\lambda}\,(q_{\epsilon}
- c_{\varphi})\right], \\
 \inprod{\alpha_{3}}{\beta_{3}} & =
-\frac{\lambda}{\tilde{\lambda}}\,
\left[\tilde{\mathcal{R}}_{1,\lambda}+\tilde{a}_{1,\lambda}\,q_{\epsilon}\right].
\end{align}
As an illustration, let us do the computation that led to
\eqref{example} in detail. As in \cite{Trac2}, we use the facts
that $S^{t}=S$, and $(S+S\,R^{t})\,\rchi=R$ which can be easily
seen by writing $R=\sum_{k=1}^{\infty}(S\,\rchi)^{k}$. Furthermore
we write $R(x,a_{k})$ to mean
\begin{equation*} \lim_{\substack{y\to a_{k}\\y\in J}}R(x,y).\end{equation*}
In general, since all evaluations are done by taking the limits
from within $J$, we can use the identity
$\rchi\,\delta_{k}=\delta_{k}$ inside the inner products. Thus
\begin{align*}
\inprod{\alpha_{2}}{\beta_{2}} & =
\frac{\lambda}{\tilde{\lambda}}\,\left[
\inprod{(S+R\,S)\,(1-\lambda\,\rchi)}{\delta_{t}-\delta_{\infty}}
+
a_{1,\lambda}\,\inprod{Q_{\epsilon}}{\delta_{t}-\delta_{\infty}}\right]\\
&=\frac{\lambda}{\tilde{\lambda}}\,\left[\inprod{(1-\lambda\,\rchi)}{(S+R^{t}\,S)\,\left(\delta_{t}-\delta_{\infty}\right)}
+ a_{1,\lambda}\left(Q_{\epsilon}(t)-Q_{\epsilon}(\infty)\right)\right]\\
&
=\frac{\lambda}{\tilde{\lambda}}\,\left[\inprod{(1-\lambda\,\rchi)}{(S+R^{t}\,S)\,\rchi\,\left(\delta_{t}-\delta_{\infty}\right)}
+ a_{1,\lambda}\left(q_{\epsilon}-c_{\varphi}\right)\right]\\
&
=\frac{\lambda}{\tilde{\lambda}}\,\left[\inprod{(1-\lambda\,\rchi)}{R(x,t)-R(x,\infty)}
+ a_{1,\lambda}\left(q_{\epsilon}-c_{\varphi}\right)\right]\\
&
=\frac{\lambda}{\tilde{\lambda}}\,\left[\mathcal{R}_{1,\lambda}(t)-\mathcal{R}_{1,\lambda}(\infty)
+ a_{1,\lambda}\left(q_{\epsilon}-c_{\varphi}\right)\right]\\
&
=\frac{\lambda}{\tilde{\lambda}}\,\left[\mathcal{R}_{1,\lambda}(t)
+ a_{1,\lambda}\left(q_{\epsilon}-c_{\varphi}\right)\right].
\end{align*}
We want the limit of the determinant
\begin{equation}\det\left(\delta_{jk}-\inprod{\alpha_{j}}{\beta_{k}}\right)_{1\leq
j,k\leq 3},\end{equation} as $N\to \infty$. In order to get our
hands on the limits of the individual terms involved in the
determinant, we will find differential equations for them first as
in \cite{Trac2}. Row operation on the matrix show that $a_{1,\lambda}$ and
$\tilde{a}_{1,\lambda}$ fall out of the determinant; to see this
add $\lambda\,a_{1,\lambda}/(2\,\tilde{\lambda})$ times row 1 to
row 2 and $\lambda\,\tilde{a}_{1,\lambda}/\tilde{\lambda}$ times
row 1 to row 3. So we will not need to find differential equations
for them. Our determinant is
{\large
  \begin{equation}
    \det\left(
      \begin{array}{ccc}
        1 - \tilde{v}_{\epsilon} & -(q_{\epsilon}-c_{\varphi}) & -q_{\epsilon} \\[6pt]
        -\frac{\lambda\,\mathcal{P}_{1,\lambda}}{2\,\tilde{\lambda}} & 1 - \frac{\lambda\,\mathcal{R}_{1,\lambda}}{2\,\tilde{\lambda}}  & -\frac{\lambda\,\mathcal{R}_{1,\lambda}}{2\,\tilde{\lambda}} \\[6pt] 
        \frac{\lambda\,\tilde{\mathcal{P}}_{1,\lambda}}{\tilde{\lambda}} & \frac{\lambda\,\tilde{\mathcal{R}}_{1,\lambda}}{\tilde{\lambda}} & 1 + \frac{\lambda\,\tilde{\mathcal{R}}_{1,\lambda}}{\tilde{\lambda}}
      \end{array}
    \right).
  \end{equation}
}
Proceeding as in \cite{Trac2} we find the following differential
equations
\begin{align} \frac{d}{d\,t}\,u_{\epsilon} & = q_{_{N}}\,q_{\epsilon}, & \frac{d}{d\,t}\,q_{\epsilon} & = q_{_{N}} -q_{_{N}}\, \tilde{v}_{\epsilon} - p_{_{N}}\,u_{\epsilon},\\
 \frac{d}{d\,t}\mathcal{Q}_{1,\lambda} & =
q_{_{N}}\left(\lambda-\mathcal{R}_{1,\lambda}\right), &
\frac{d}{d\,t}\mathcal{P}_{1,\lambda} & =
p_{_{N}}\left(\lambda-\mathcal{R}_{1,\lambda}\right),\label{ex2}\\
 \frac{d}{d\,t}\mathcal{R}_{1,\lambda} & =
-p_{_{N}}\,\mathcal{Q}_{1,\lambda}-q_{_{N}}\,\mathcal{P}_{1,\lambda},
  & \frac{d}{d\,t}\tilde{\mathcal{R}}_{1,\lambda} & =
-p_{_{N}}\,\tilde{\mathcal{Q}}_{1,\lambda}-q_{_{N}}\,\tilde{\mathcal{P}}_{1,\lambda},
\\  \frac{d}{d\,t}\tilde{\mathcal{Q}}_{1,\lambda} & =
q_{_{N}}\left(\lambda-1-\tilde{\mathcal{R}}_{1,\lambda}\right),
 & \frac{d}{d\,t}\tilde{\mathcal{P}}_{1,\lambda}  & =
p_{_{N}}\left(\lambda-1-\tilde{\mathcal{R}}_{1,\lambda}\right).\end{align}
Let us derive the first equation in \eqref{ex2} for example. From
\cite{Trac3} (equation $2.17$), we have
\begin{equation*}
\frac{\partial Q}{\partial t}=-R(x,t)\,q_{_{N}}.
\end{equation*}
Therefore
\begin{align*}
\frac{\partial \mathcal{Q}_{1,\lambda}}{\partial t} & =
\frac{d}{d\,t}\left[\int_{-\infty}^{t}Q(x,t)\,d\,x-(1-\lambda)\,\int_{\infty}^{t}Q(x,t)\,d\,x\right]\\
& = q_{_{N}} + \int_{-\infty}^{t}\frac{\partial Q}{\partial
t}\,d\,x -
(1-\lambda)\left[q_{_{N}}+\int_{\infty}^{t}\frac{\partial
Q}{\partial t}\,d\,x\right]
\\
& = q_{_{N}} - q_{_{N}}\int_{-\infty}^{t}R(x,t)\,d\,x -
(1-\lambda)\,q_{_{N}}+(1-\lambda)\,q_{_{N}}\,\int_{\infty}^{t}R(x,t)\,d\,x \\
& = \lambda\,q_{_{N}}-q_{_{N}}\,\int_{-\infty}^{\infty}(1-\lambda)\,R(x,t)\,d\,x\\
& = \lambda\,q_{_{N}}-q_{_{N}}\,\mathcal{R}_{1,\lambda} =
q_{_{N}}\left(\lambda-\mathcal{R}_{1,\lambda}\right).
\end{align*}
 Now we change variable from $t$ to $s$ where $t=\tau(s)=
2\,\sigma\,\sqrt{N}+\frac{\sigma\,s}{N^{1/6}}$. Then we take the
limit $N\to \infty$, denoting the limits of $ q_{\epsilon},
\mathcal{P}_{1,\lambda},
\mathcal{Q}_{1,\lambda},\mathcal{R}_{1,\lambda},
\tilde{\mathcal{P}}_{1,\lambda} ,\tilde{\mathcal{Q}}_{1,\lambda} ,
\tilde{\mathcal{R}}_{1,\lambda}$ and the common limit of
$u_{\epsilon}$ and $\tilde{v}_{\epsilon}$ respectively by
$\overline{q}, \overline{\mathcal{P}}_{1,\lambda},
\overline{\mathcal{Q}}_{1,\lambda},
\overline{\mathcal{R}}_{1,\lambda},
\overline{\overline{\mathcal{P}}}_{1,\lambda} ,
\overline{\overline{\mathcal{Q}}}_{1,\lambda} ,
\overline{\overline{\mathcal{R}}}_{1,\lambda}$ and $\overline{u}$.
We eliminate $ \overline{\mathcal{Q}}_{1,\lambda}$ and
$\overline{\overline{\mathcal{Q}}}_{1,\lambda}$ by using the facts
that
$\overline{\mathcal{Q}}_{1,\lambda}=\overline{\mathcal{P}}_{1,\lambda}+\lambda\,\sqrt{2}$
and $\overline{\overline{\mathcal{Q}}}_{1,\lambda}=
\overline{\overline{\mathcal{P}}}_{1,\lambda}$. These limits hold
uniformly for bounded $s$ so we can interchange $\lim$ and
$\frac{d}{d\,s}$. Also
$\lim_{N\to\infty}N^{-1/6}q_{_{N}}=\lim_{N\to\infty}N^{-1/6}p_{_{N}}=q
$ , where $q$ is as in \eqref{D2}. We obtain the systems
\begin{equation}\frac{d}{d\,s}\,\overline{u} = -\frac{1}{\sqrt{2}}\,q\,\overline{q}
,
 \qquad
\frac{d}{d\,s}\,\overline{q} =
\frac{1}{\sqrt{2}}\,q\,\left(1-2\,\overline{u}\right),
\end{equation}
\begin{equation}\frac{d}{d\,s}\overline{\mathcal{P}}_{1,\lambda} =
-\frac{1}{\sqrt{2}}\,q\,\left(\overline{\mathcal{R}}_{1,\lambda}-\lambda\right),
\qquad \frac{d}{d\,s}\overline{\mathcal{R}}_{1,\lambda} =
-\frac{1}{\sqrt{2}}\,q\,\left(2\,\overline{\mathcal{P}}_{1,\lambda}+\sqrt{2\tilde{\lambda}}\right),
\end{equation}
\begin{equation}\frac{d}{d\,s}\overline{\overline{\mathcal{P}}}_{1,\lambda} =
\frac{1}{\sqrt{2}}\,q\,\left(1-\lambda-\overline{\overline{\mathcal{R}}}_{1,\lambda}\right),
\qquad
 \frac{d}{d\,s}\overline{\overline{\mathcal{R}}}_{1,\lambda} =
 -q\,\sqrt{2}\,\overline{\overline{\mathcal{P}}}_{1,\lambda}.
\end{equation}
 The change of variables $s\to\mu=\int_{s}^{\infty} q(x, \lambda)\,d\,x$
transforms these systems into constant coefficient ordinary differential equations

\begin{equation}\frac{d}{d\,\mu}\overline{u} =
\frac{1}{\sqrt{2}}\,\overline{q}, \qquad
\frac{d}{d\,\mu}\overline{q} =
-\frac{1}{\sqrt{2}}\,\left(1-2\,\overline{u}\right),
\end{equation}
\begin{equation}\frac{d}{d\,\mu}\overline{\mathcal{P}}_{1,\lambda} =
\frac{1}{\sqrt{2}}\,\left(\overline{\mathcal{R}}_{1,\lambda}-\lambda\right),
\qquad \frac{d}{d\,\mu}\overline{\mathcal{R}}_{1,\lambda} =
\frac{1}{\sqrt{2}}\,\left(2\,\overline{\mathcal{P}}_{1,\lambda}+\sqrt{2\tilde{\lambda}}\right),
\end{equation}
\begin{equation}\frac{d}{d\,\mu}\overline{\overline{\mathcal{P}}}_{1,\lambda} =
-\frac{1}{\sqrt{2}}\,\left(1-\lambda-\overline{\overline{\mathcal{R}}}_{1,\lambda}\right),
\qquad
 \frac{d}{d\,\mu}\overline{\overline{\mathcal{R}}}_{1,\lambda} =
 \sqrt{2}\,\overline{\overline{\mathcal{P}}}_{1,\lambda}.
\end{equation} Since $\lim_{s\to \infty}\mu=0$, corresponding to the boundary
values at $t=\infty$ which we found earlier for
$\mathcal{P}_{1,\lambda}, \mathcal{R}_{1,\lambda},
\tilde{\mathcal{P}}_{1,\lambda} ,
\tilde{\mathcal{R}}_{1,\lambda}$, we now have initial values at
$\mu=0$. Therefore
\begin{equation} \overline{\mathcal{P}}_{1,\lambda}(0) = \overline{\mathcal{R}}_{1,\lambda}(0)=
\overline{\overline{\mathcal{P}}}_{1,\lambda}(0)=
\overline{\overline{\mathcal{R}}}_{1,\lambda}(0)=0.\end{equation}
We use this to solve the systems and get
\begin{align}\overline{q} & = \frac{\sqrt{\tilde{\lambda}}-1}{2\,\sqrt{2}}\,e^{\mu} + \frac{\sqrt{\tilde{\lambda}} + 1}{2\,\sqrt{2}}\,e^{-\mu},\\
 \overline{u} & =
\frac{\sqrt{\tilde{\lambda}}-1}{4}\,e^{\mu} -
\frac{\sqrt{\tilde{\lambda}}+1}{4}\,e^{-\mu}+\frac{1}{2},\\
\overline{\mathcal{P}}_{1,\lambda} & =
\frac{\sqrt{\tilde{\lambda}}-\lambda}{2\,\sqrt{2}}\,e^{\mu} +
\frac{\sqrt{\tilde{\lambda}} + \lambda}{2\,\sqrt{2}}\,e^{-\mu} -
\sqrt{\frac{\tilde{\lambda}}{2}},\\
 \overline{\mathcal{R}}_{1,\lambda} & =
\frac{\sqrt{\tilde{\lambda}}-\lambda}{2}\,e^{\mu} -
\frac{\sqrt{\tilde{\lambda}} + \lambda}{2}\,e^{-\mu} +
\lambda, \\
\overline{\overline{\mathcal{P}}}_{1,\lambda} & =
\frac{1-\lambda}{2\,\sqrt{2}}\,(e^{\mu}-e^{-\mu}), \qquad
\overline{\overline{\mathcal{R}}}_{1,\lambda}=\frac{1-\lambda}{2}\,(e^{\mu}+e^{-\mu}-2).
\end{align}
Substituting these expressions into the determinant gives
\eqref{goedet}, namely
\begin{equation} D_{1}(s,\lambda)= D_{2}(s,\tilde{\lambda})\,\frac{\lambda - 1
- \cosh{\mu(s,\tilde{\lambda})} +
\sqrt{\tilde{\lambda}}\,\sinh{\mu(s,\tilde{\lambda})}}{\lambda -
2},\end{equation} where $D_{\beta}=\lim_{N\to\infty}D_{\beta,N}$.
As mentioned in Section \ref{sec:edgeDistr}, the functional form of the $\beta=1$ limiting determinant is very different from what one would expect, unlike in the $\beta=4$ case. Also noteworthy is the dependence on $\tilde{\lambda}=2\lambda-\lambda^{2}$ instead of just $\lambda$. However one should also note that when $\lambda$ is set equal to $1$, then $\tilde{\lambda}=\lambda=1$. Hence in the largest eigenvalue case, where there is no prior differentiation with respect to $\lambda$, and $\lambda$ is just set to $1$, a great deal of simplification occurs. The above formula then nicely reduces to the $\beta=1$ Tracy-Widom distribution. 

\section{An Interlacing Property}

The following series of lemmas establish
Corollary~\eqref{interlacingcor}:
\begin{lemma}
Define
\begin{equation}a_{j}=\frac{d^{j}}{d\,\lambda^{j}}\,\sqrt{\frac{\lambda}{2-\lambda}}\,\,\bigg{\vert}_{\lambda=1}.\label{ajdef}\end{equation}
Then $a_{j}$ satisfies the following recursion
\begin{equation} a_{j} = \begin{cases}
          \quad 1 &\text{if} \quad j=0, \\
           \quad (j-1)\,a_{j-1} &\text{for $j\geq 1$, $j$ even,}  \\
           \quad j\,a_{j-1} &\text{for $j\geq 1$, $j$ odd.} \\
         \end{cases}\end{equation}
\label{ajlemma}\end{lemma}
\begin{proof}
Consider the expansion of the generating function $f(\lambda)=\sqrt{\frac{\lambda}{2-\lambda}}$ around $\lambda=1$
\begin{equation*}
f(\lambda)= \sum_{j\geq 0}\frac{a_{j}}{j!}\,(\lambda-1)^{j}=\sum_{j\geq 0} b_{j}\,(\lambda-1)^{j}.
\end{equation*}
Since $a_{j}=j!\,b_{j}$, the statement of the lemma reduces to proving the following recurrence for the $b_{j}$ 
\begin{equation}
  b_{j} =
  \begin{cases}
    \quad 1 &\text{if} \quad j=0, \\
    \quad \frac{j-1}{j}\,b_{j-1} &\text{for $j\geq 1$, $j$ even,}  \\
    \quad b_{j-1} &\text{for $j\geq 1$, $j$ odd.} \\
  \end{cases}\label{bjrecurrence}
\end{equation}
Let
\begin{equation*}
  f^{even}(\lambda)= \frac{1}{2}\left(\sqrt{\frac{\lambda}{2-\lambda}}+\sqrt{\frac{2-\lambda}{\lambda}}\right), \qquad f^{odd}(\lambda)= \frac{1}{2}\left(\sqrt{\frac{\lambda}{2-\lambda}}-\sqrt{\frac{2-\lambda}{\lambda}}\right).
\end{equation*}
These are the even and odd parts of $f$ relative to the reflection $\lambda-1\to-(\lambda-1)$ or $\lambda\to 2-\lambda$. Recurrence \eqref{bjrecurrence} is equivalent to
\begin{equation*}
  \frac{d}{d\,\lambda}\,f^{even}(\lambda)=(\lambda-1)\,\frac{d}{d\,\lambda}\,f^{odd}(\lambda),
\end{equation*}
which is easily shown to be true.
\end{proof}
\begin{lemma}
  \label{flemma}
  Define
  \begin{equation}
    f(s,\lambda)=1-\sqrt{\frac{\lambda}{2-\lambda}}\,\,\tanh{\frac{\mu(s,\tilde{\lambda})}{2}},
  \end{equation}
  for $\tilde{\lambda}=2\,\lambda-\lambda^{2}$. Then
  \begin{equation}\frac{\partial^{2\,n}}{\partial\,\lambda^{2\,n}}\,f(s,\lambda)\,\,\bigg{\vert}_{\lambda=1}-
    \frac{1}{2\,n+1}\,\frac{\partial^{2\,n+1}}{\partial\,\lambda^{2\,n+1}}\,f(s,\lambda)\,\,\bigg{\vert}_{\lambda=1}=
    \begin{cases}
      \quad 1 &\text{if $n=0$,}\\
      \quad0 &\text{if $n\geq 1$.}  \\
    \end{cases} \label{flemmaeq}\end{equation}
\end{lemma}
\begin{proof}
The case $n=0$ is readily checked. The main ingredient for the
general case is Fa\'a di Bruno's formula
\begin{equation}
\frac{d^{n}}{d t^{n}}g(h(t))=\sum\frac{n!}{k_{1}!\cdots k_{n}!}
\left(\frac{d^{k}g}{dh^{k}}(h(t))\right)\left(\frac{1}{1!}\frac{dh}{d
t}\right)^{k_{1}}\cdots\left(\frac{1}{n!}\frac{d^{n} h}{d
t^{n}}\right)^{k_{n}}, \label{faa}\end{equation} where
$k=\sum_{i=1}^{n}k_{i}$ and the above sum is over all partitions
of $n$, that is all values of $k_{1},\ldots, k_{n}$ such that
$\sum_{i=1}^{n} i\,k_{i}=n$. We apply Fa\'a di Bruno's formula to
derivatives of the function
$\tanh{\frac{\mu(s,\tilde{\lambda})}{2}}$, which we treat as some
function $g(\tilde{\lambda}(\lambda))$. Notice that for $j\geq 1$,
$\frac{d^{j}\tilde{\lambda}}{d\,\lambda^{j}}\,\,\big{\vert}_{\lambda=1}$
is nonzero only when $j=2$, in which case it equals $-2$. Hence,
in \eqref{faa}, the only term that survives is the one
corresponding to the partition all of whose parts equal $2$. Thus
we have
\begin{align*}\frac{\partial^{2n-k}}{\partial\,\lambda^{2n-k}}\,&\tanh{\frac{\mu(s,\tilde{\lambda})}{2}}\,\,\bigg{\vert}_{\lambda=1}&&\\
 &= \begin{cases}
           0 &\text{if $k=2j+1$, $j\geq 0$}\\
          \frac{(-1)^{n-j}\,(2\,n-k)!}{(n-j)!} \frac{\partial^{n-j}}{\partial\,\tilde{\lambda}^{n-j}}\,\tanh{\frac{\mu(s,\tilde{\lambda})}{2}}\,\,\bigg{\vert}_{\tilde{\lambda}=1} &\text{for $k=2j$, $j\geq 0$}  \\
         \end{cases} \end{align*}

 \begin{align*}\frac{\partial^{2n-k+1}}{\partial\,\lambda^{2n+1-k}}\,&\tanh{\frac{\mu(s,\tilde{\lambda})}{2}}\,\,\bigg{\vert}_{\lambda=1}
 \\
 &= \begin{cases}
           0 &\text{if $k=2j$, $j\geq 0$}\\
         \frac{(-1)^{n-j}\,(2\,n+1-k)!}{(n-j)!} \frac{\partial^{n-j}}{\partial\,\tilde{\lambda}^{n-j}}\,\tanh{\frac{\mu(s,\tilde{\lambda})}{2}}\,\,\bigg{\vert}_{\tilde{\lambda}=1} &\text{for $k=2j+1$, $j\geq 0$}  \\
         \end{cases} \end{align*}
Therefore, recalling the definition of $a_{j}$ in \eqref{ajdef}
and setting $k=2\,j$, we obtain
\begin{eqnarray*}
\frac{\partial^{2\,n}}{\partial\,\lambda^{2\,n}}\,f(s,\lambda)\,\,\bigg{\vert}_{\lambda=1}&=&
\sum_{k=0}^{2\,n}\binom{2\,n}{k}\,\frac{\partial^{k}}{\partial\,\lambda^{k}}\,\sqrt{\frac{\lambda}{2-\lambda}}\,\frac{\partial^{2\,n-k}}{\partial\,\lambda^{2\,n-k}}\,\tanh{\frac{\mu(s,\tilde{\lambda})}{2}}\,\,\bigg{\vert}_{\lambda=1}\\
&=&
\sum_{j=0}^{n}\frac{(2\,n)!\,(-1)^{n-j}}{(2\,j)!\,(n-j)!}\,a_{2\,j}\,\frac{\partial^{n-j}}{\partial\,\tilde{\lambda}^{n-j}}\,\tanh{\frac{\mu(s,\tilde{\lambda})}{2}}\,\,\bigg{\vert}_{\tilde{\lambda}=1}.
\end{eqnarray*}
Similarly, using $k=2\,j+1$ instead yields
\begin{eqnarray*}
  \frac{\partial^{2\,n+1}}{\partial\,\lambda^{2\,n+1}} \,f(s,\lambda)\,\bigg{\vert}_{\lambda=1} &=& \sum_{k=0}^{2\,n+1}\binom{2\,n+1}{k}\,\frac{\partial^{k}}{\partial\,\lambda^{k}}\,\sqrt{\frac{\lambda}{2-\lambda}}\,\frac{\partial^{2\,n+1-k}}{\partial\,\lambda^{2\,n+1-k}}\,\tanh{\frac{\mu(s,\tilde{\lambda})}{2}}\,\,\bigg{\vert}_{\lambda=1}\\
  &=& (2\,n+1)\,\sum_{j=0}^{n}\frac{(2\,n)!\,(-1)^{n-j}}{(2\,j)!\,(n-j)!}\,\frac{a_{2\,j+1}}{2\,j+1}\,\frac{\partial^{n-j}}{\partial\,\tilde{\lambda}^{n-j}}\,\tanh{\frac{\mu(s,\tilde{\lambda})}{2}}\,\,\bigg{\vert}_{\tilde{\lambda}=1}\\
  &=& (2\,n+1)\,\frac{\partial^{2\,n}}{\partial\,\lambda^{2\,n}}\,f(s,\lambda)\,\,\bigg{\vert}_{\lambda=1},
\end{eqnarray*}
since $a_{_{2\,j+1}}/(2\,j+1)=a_{_{2j}}$. Rearranging this last
equality leads to \eqref{flemmaeq}.
\end{proof}

\begin{lemma}
Let $D_{1}(s,\lambda)$ and $D_{4}(s,\tilde{\lambda})$ be as in
\eqref{goedet} and \eqref{gsedet}. Then
\begin{equation}D_{1}(s,\lambda)=D_{4}(s,\tilde{\lambda})\,
\left(1-\sqrt{\frac{\lambda}{2-\lambda}}\,\tanh{\frac{\mu(s,\tilde{\lambda})}{2}}\right)^{2}.\end{equation}
\end{lemma}
\begin{proof}
Using the facts that $-1-\cosh{x}=-2\,\cosh^{2}\frac{x}{2}$,
$1=\cosh^{2}{x}-\sinh^{2}{x}$ and
$\sinh{x}=2\,\sinh\frac{x}{2}\,\cosh\frac{x}{2}$ we get

\begin{eqnarray*}
  D_{1}(s,\lambda) &=& \frac{-1-\cosh{\mu(s,\tilde{\lambda})}}{\lambda-2}\,D_{2}(s,\tilde{\lambda}) + D_{2}(s,\tilde{\lambda})\,\frac{\lambda + \sqrt{\tilde{\lambda}}\,\sinh{\mu(s,\tilde{\lambda})}}{\lambda-2}\\
  &=&\frac{-2}{\lambda-2}\,D_{4}(s,\tilde{\lambda}) \\
  & & \qquad + D_{2}(s,\tilde{\lambda})\,\frac{\lambda\left(\cosh^{2}{\frac{\mu(s,\tilde{\lambda})}{2}} - \sinh^{2}{\frac{\mu(s,\tilde{\lambda})}{2}}\right) + \sqrt{\tilde{\lambda}}\,\sinh{\mu(s,\tilde{\lambda})}}{\lambda-2}\\
  &=& D_{4}(s,\tilde{\lambda}) + \frac{D_{4}(s,\tilde{\lambda})}{\cosh^{2}\left(\frac{\mu(s,\lambda)}{2}\right)}\,\frac{\lambda\,\sinh^{2}{\frac{\mu(s,\tilde{\lambda})}{2}} - \sqrt{\tilde{\lambda}}\,\sinh{\mu(s,\tilde{\lambda})}}{2 - \lambda}\\
  &=& D_{4}(s,\tilde{\lambda})\,\left(1 + \frac{\lambda\,\sinh^{2}{\frac{\mu(s,\tilde{\lambda})}{2}} - 2\,\sqrt{\tilde{\lambda}}\,\sinh\left(\frac{\mu(s,\lambda)}{2}\right)\,\cosh\left(\frac{\mu(s,\lambda)}{2}\right)}{(2-\lambda)\,\cosh^{2}\left(\frac{\mu(s,\lambda)}{2}\right)}\right)\\
  &=& D_{4}(s,\tilde{\lambda})\,\left(1 - 2\,\sqrt{\frac{\lambda}{2-\lambda}}\,\tanh{\frac{\mu(s,\tilde{\lambda})}{2}} + \frac{\lambda}{2-\lambda}\,\tanh^{2}{\frac{\mu(s,\tilde{\lambda})}{2}} \right)\\
  &=& D_{4}(s,\tilde{\lambda})\,\left(1-\sqrt{\frac{\lambda}{2-\lambda}}\,\,\tanh{\frac{\mu(s,\tilde{\lambda})}{2}}\right)^{2}.
\end{eqnarray*}
\end{proof}
For notational convenience, define $d_{1}(s,\lambda)=D_{1}^{1/2}(s,\lambda)$, and $d_{4}(s,\lambda)=D_{4}^{1/2}(s,\lambda)$. Then

\begin{lemma}For $n\geq 0$,
\begin{equation*}\left[-\frac{1}{(2\,n+1)!}\,\frac{\partial^{2\,n+1}}{\partial\,\lambda^{2\,n+1}}
+\frac{1}{(2\,n)!}\,\frac{\partial^{2\,n}}{\partial\,\lambda^{2\,n}}\right]\,d_{1}(s,\lambda)\,\,\bigg{\vert}_{\lambda=1}=\frac{(-1)^{n}}{n!}\,\frac{\partial^{n}}{\partial\,\lambda^{n}}\,d_{4}(s,\lambda)\,\,\bigg{\vert}_{\lambda=1}.\end{equation*}
\label{dlemma}\end{lemma}
\begin{proof}
Let
\begin{equation*}f(s,\lambda)=1-\sqrt{\frac{\lambda}{2-\lambda}}\,\,\tanh{\frac{\mu(s,\tilde{\lambda})}{2}}\end{equation*}
by the previous lemma, we need to show that
\begin{equation}
  \label{eq:38}
  \left[-\frac{1}{(2\,n+1)!}\,\frac{\partial^{2\,n+1}}{\partial\,\lambda^{2\,n+1}} + \frac{1}{(2\,n)!}\,\frac{\partial^{2\,n}}{\partial\,\lambda^{2\,n}}\right]\,d_{4}(s,\tilde{\lambda})\,f(s,\lambda)\,\,\bigg{\vert}_{\lambda=1}
\end{equation}
equals
\begin{equation*}
\frac{(-1)^{n}}{n!}\,\frac{\partial^{n}}{\partial\,\tilde{\lambda}^{n}}\,d_{4}(s,\tilde{\lambda})\,\,\bigg{\vert}_{\lambda=1}.
\end{equation*} Now formula \eqref{faa} applied to $d_{4}(s,\tilde{\lambda})$
gives
\begin{equation*}\frac{\partial^{k}}{\partial\,\lambda^{k}}\,d_{4}(s,\tilde{\lambda})\,\,\bigg{\vert}_{\lambda=1}=
\begin{cases}
          \quad 0 &\text{if $k=2j+1$, $j\geq 0$,}\\
           \quad\frac{(-1)^{j}\,k!}{j!} \frac{\partial^{j}}{\partial\,\tilde{\lambda}^{j}}\,d_{4}(s,\tilde{\lambda}) &\text{if $k=2j$, $j\geq 0$.}  \\
         \end{cases} \end{equation*}
Therefore
\begin{align*}
-\frac{1}{(2\,n+1)!}\,\frac{\partial^{2\,n+1}}{\partial\,\lambda^{2\,n+1}}&\,d_{4}(s,\tilde{\lambda})\,f(s,\lambda)\,\,\bigg{\vert}_{\lambda=1}\\
&=-\frac{1}{(2\,n+1)!}\,\sum_{k=0}^{2\,n+1}\binom{2\,n+1}{k}\,\frac{\partial^{k}}{\partial\,\lambda^{k}}\,d_{4\,}
\frac{\partial^{2\,n+1-k}}{\partial\,\lambda^{2\,n+1-k}}\,f
\,\,\bigg{\vert}_{\lambda=1}\\
&=
-\sum_{j=0}^{n}\frac{(-1)^{j}}{(2\,n-2\,j+1)!\,j!}\,\frac{\partial^{j}}{\partial\,\tilde{\lambda}^{j}}\,d_{4\,}
\frac{\partial^{2\,n-2\,j+1}}{\partial\,\lambda^{2\,n-2\,j+1}}\,f
\,\,\bigg{\vert}_{\lambda=1}
\end{align*}
Similarly,
\begin{align*}
  \frac{1}{(2\,n)!}\,\frac{\partial^{2\,n}}{\partial\,\lambda^{2\,n}}\,d_{4}(s,\tilde{\lambda})\,f(s,\lambda)\,\,\bigg{\vert}_{\lambda=1}
  &=\frac{1}{(2\,n)!}\,\sum_{k=0}^{2\,n}\binom{2\,n}{k}\,\frac{\partial^{k}}{\partial\,\lambda^{k}}\,d_{4\,}\frac{\partial^{2\,n-k}}{\partial\,\lambda^{2\,n-k}}\,f\,\,\bigg{\vert}_{\lambda=1}\\
  &=\sum_{j=0}^{n}\frac{(-1)^{j}}{(2\,n-2\,j)!\,j!}\,\frac{\partial^{j}}{\partial\,\tilde{\lambda}^{j}}\,d_{4\,}\frac{\partial^{2\,n-2\,j}}{\partial\,\lambda^{2\,n-2\,j}}\,f\,\,\bigg{\vert}_{\lambda=1}.
\end{align*}
Therefore the expression in \eqref{eq:38} equals
\begin{align*}
  \sum_{j=0}^{n}\frac{(-1)^{j}}{(2\,n-2\,j)!\,j!}\,\frac{\partial^{j}}{\partial\,\tilde{\lambda}^{j}}\,d_{4\,}(s,\tilde{\lambda}) \left[\frac{\partial^{2\,n-2\,j}}{\partial\,\lambda^{2\,n-2\,j}}\,f- \frac{1}{2\,n-2\,j+1}\,\frac{\partial^{2\,n-2\,j+1}}{\partial\,\lambda^{2\,n-2\,j+1}}\,f\right]\,\,\bigg{\vert}_{\lambda=1}.
\end{align*}

Now Lemma~\ref{flemma} shows that the square bracket inside the
summation is zero unless $j=n$, in which case it is $1$. The
result follows.
\end{proof}
In an inductive proof of Corollary~\ref{interlacingcor}, the base case $F_{4}(s,2)=F_{1}(s,1)$ is easily checked by direct calculation. Lemma~\ref{dlemma} establishes the inductive step in the proof since, with the assumption $F_{4}(s,n)=F_{1}(s,2\,n)$, it is equivalent to the statement
\begin{equation*}
  F_{4}(s,n+1)=F_{1}(s,2\,n+2).
\end{equation*}

\section{Numerics}\label{sec:numerics}

\subsection{Partial derivatives of $q(x,\lambda)$}

Let \begin{equation}q_{n}(x)=\frac{\partial^{n}}{\partial\lambda^{n}}\,q(x,\lambda)\bigg\vert_{\lambda=1},\end{equation}
so that $q_{0}$ equals $q$ from \eqref{pII}.  In order to compute
$F_{\beta}(s,m)$ it is crucial to know $q_{n}$ with $0\leq n\leq m$ accurately.
Asymptotic expansions for $q_{n}$ at $-\infty$ are given in
\cite{Trac3}.  In particular, we know that, as $t\to+\infty$, $q_{0}(-t/2)$ is given by
\begin{equation*}
  \frac{1}{2}\sqrt{t}\left(1-\frac{1}{t^{3}}-\frac{73}{2t^{6}}-\frac{10657}{2t^{9}}-\frac{13912277}{8t^{12}}+\textrm{O}\left(\frac{1}{t^{15}}\right)\right),
\end{equation*}
whereas $q_{1}(-t/2)$ can be expanded as
\begin{equation*}
  \frac{\exp{(\frac{1}{3}t^{3/2})}}{2\sqrt{2\pi}\,t^{1/4}}\left(1+\frac{17}{24t^{3/2}}+\frac{1513}{2^{7}3^{2}t^{3}}+\frac{850193}{2^{10}3^{4}t^{9/2}}-\frac{407117521}{2^{15}3^{5}t^{6}}+\textrm{O}\left(\frac{1}{t^{15/2}}\right)\right).
\end{equation*}
These expansions are used in the algorithms below.
\subsection{Algorithms}
Quantities needed to compute $F_{\beta}(s,m), m=1,2,$ are not only
$q_{0}$ and $q_{1}$ but also integrals involving $q_{0}$, such as
\begin{equation}
I_{0}=\int_{s}^{\infty}(x-s)\,q_{0}^{2}(x)\,d\,x,\quad
J_{0}=\int_{s}^{\infty}q_{0}(x)\,d\,x.
\end{equation}
Instead of computing these integrals afterward, it is better to
include them as variables in a system together with $q_{0}$, as
suggested in \cite{Pers1}. Therefore all quantities needed are
computed in one step, greatly reducing errors, and taking full
advantage of the powerful numerical tools in MATLAB\texttrademark\,. Since
\begin{equation}
I_{0}'=-\int_{s}^{\infty}q_{0}^{2}(x)\,d\,x,\quad
I_{0}''=q_{0}^{2},\qquad J_{0}'=-q_{0},
\end{equation}
the system closes, and can be concisely written
\begin{equation}
\frac{d}{ds}\left(\begin{array}{c}q_{0} \\ q_{0}' \\ I_{0} \\ I_{0}' \\
J_{0}
\end{array}\right) =\left(\begin{array}{c}q_{0}' \\ s\,q_{0}+2q_{0}^3 \\ I_{0}' \\ q_{0}^{2} \\ -q_{0}
\end{array}\right).
\label{q0sys}\end{equation} We first use the MATLAB\texttrademark\, built--in
Runge--Kutta--based ODE solver \texttt{ode45} to obtain a first
approximation to  the solution of \eqref{q0sys} between $x=6$, and
$x=-8$, with an initial values obtained using the Airy function on
the right hand side. Note that it is not possible to extend the
range to the left due to the high instability of the solution a
little after $-8$. (This is where the transition region between
the three different regimes in the so--called ``connection
problem'' lies.  We circumvent this limitation by patching up our
solution with the asymptotic expansion to the left of $x=-8$.)
The approximation obtained is then used as a trial solution in the
MATLAB\texttrademark\, boundary value problem solver \texttt{bvp4c}, resulting in
an accurate solution vector between $x=6$ and $x=-10$.
Similarly, if we define
\begin{equation}
I_{1}=\int_{s}^{\infty}(x-s)\,q_{0}(x)\,q_{1}(x)\,d\,x,\quad
J_{1}=\int_{s}^{\infty}q_{0}(x)\,q_{1}(x)\,d\,x,
\end{equation}
then we have the first--order system
\begin{equation}
\frac{d}{ds}\left(\begin{array}{c}q_{1} \\ q_{1}' \\ I_{1} \\ I_{1}' \\
J_{1}
\end{array}\right) =\left(\begin{array}{c}q_{1}' \\ s\,q_{1}+6q_{0}^2\,q_{1} \\ I_{1}' \\ q_{0}\,q_{1} \\ -q_{0}\,q_{1}
\end{array}\right),
\end{equation}
which can be implemented using \texttt{bvp4c} together with a
``seed'' solution obtained in the same way as for $q_{0}$.

The MATLAB\texttrademark\, code is freely available, and may be obtained by contacting the first author.
\clearpage
\subsection{Tables}
\begin{center}
{\large
\begin{table}[!h]
\label{tab:2}
\newcommand\T{\rule{0pt}{2.6ex}}
\newcommand\B{\rule[-1.2ex]{0pt}{0pt}}
\begin{tabular}[!b]{c|c|c|c|c|}
\cline{2-5}
\T \B & \multicolumn{4}{c|}{\textbf{Statistic}} \\
\hline
\multicolumn{1}{|c|}{\textbf{Eigenvalue}}  \T \B & $\mu$ & $\sigma$ & $\gamma_{1}$ & $\gamma_{2}$ \\
\hline
\multicolumn{1}{|c|}{$\lambda_{1}$} \T & $-1.771087$ & $0.901773$ & $0.224084$ & $0.093448$ \\[6pt]
\multicolumn{1}{|c|}{$\lambda_{2}$} \T \B & $-3.675440$ & $0.735214$ & $0.125000$ & $0.021650$ \\
\hline
\end{tabular}
\caption{Mean, standard deviation, skewness and kurtosis data for first two edge--scaled eigenvalues in the $\beta=2$ Gaussian ensemble. Compare to Table~$1$ in \cite{Trac4}.}
\end{table}
}
\end{center}

\begin{center}
{\large
\begin{table}[!h]
\label{tab:1}
\newcommand\T{\rule{0pt}{2.6ex}}
\newcommand\B{\rule[-1.2ex]{0pt}{0pt}}
\begin{tabular}[!t]{c|c|c|c|c|}
\cline{2-5}
\T \B & \multicolumn{4}{c|}{\textbf{Statistic}} \\
\hline
\multicolumn{1}{|c|}{\textbf{Eigenvalue}}  \T \B & $\mu$ & $\sigma$ & $\gamma_{1}$ & $\gamma_{2}$ \\
\hline
\multicolumn{1}{|c|}{$\lambda_{1}$} \T & $-1.206548$ & $1.267941$ & $0.293115$ & $0.163186$ \\[6pt]
\multicolumn{1}{|c|}{$\lambda_{2}$} \T & $-3.262424$ & $1.017574$ & $0.165531$ & $0.049262$ \\[6pt]
\multicolumn{1}{|c|}{$\lambda_{3}$} \T & $-4.821636$ & $0.906849$ & $0.117557$ & $0.019506$ \\[6pt]
\multicolumn{1}{|c|}{$\lambda_{4}$} \T \B & $-6.162036$ & $0.838537$ & $0.092305$ & $0.007802$ \\
\hline
\end{tabular}
\caption{Mean, standard deviation, skewness and kurtosis data for first four edge--scaled eigenvalues in the $\beta=1$ Gaussian ensemble. Corollary~\ref{interlacingcor} implies that rows $2$ and $4$ give corresponding data for the two largest eigenvalues in the $\beta=4$ Gaussian ensemble. Compare to Table~$1$ in \cite{Trac4}, keeping in mind that the discrepancy in the $\beta=4$ data is caused by the different normalization in the 
definition of $F_{4}(s,1)$.}
\end{table}
}
\end{center}
\begin{center}
{\large
\begin{table}[!h]
\newcommand\T{\rule{0pt}{2.6ex}}
\newcommand\B{\rule[-1.2ex]{0pt}{0pt}}
\begin{tabular}[!b]{c|ccc|ccc|}
\cline{2-7}
\T \B & \multicolumn{3}{c|}{$\mathbf{100\times 100}$} & \multicolumn{3}{c|}{$\mathbf{100\times 400}$} \\
\hline
\multicolumn{1}{|c|}{\textbf{$F_{1}$-Percentile}} \T \B & $\lambda_{1}$ & $\lambda_{2}$ & $\lambda_{3}$ & $\lambda_{1}$ & $\lambda_{2}$ & $\lambda_{3}$ \\
\hline
\multicolumn{1}{|c|}{$0.01$} \T & $0.008$ & $0.005$  & $0.004$  & $0.008$ & $0.006$ & $0.004$ \\[6pt]
\multicolumn{1}{|c|}{$0.05$} & $0.042$ & $0.033$  & $0.025$  & $0.042$ & $0.037$ & $0.032$ \\[6pt]
\multicolumn{1}{|c|}{$0.10$} & $0.090$ & $0.073$  & $0.059$  & $0.088$ & $0.081$ & $0.066$ \\[6pt]
\multicolumn{1}{|c|}{$0.30$} & $0.294$ & $0.268$  & $0.235$  & $0.283$ & $0.267$ & $0.254$ \\[6pt]
\multicolumn{1}{|c|}{$0.50$} & $0.497$ & $0.477$  & $0.440$  & $0.485$ & $0.471$ & $0.455$ \\[6pt]
\multicolumn{1}{|c|}{$0.70$} & $0.699$ & $0.690$  & $0.659$  & $0.685$ & $0.679$ & $0.669$ \\[6pt]
\multicolumn{1}{|c|}{$\mathbf{0.90}$} & $\mathbf{0.902}$  & $\mathbf{0.891}$  & $\mathbf{0.901}$ & $\mathbf{0.898}$ & $\mathbf{0.894}$ & $\mathbf{0.884}$ \\[6pt]
\multicolumn{1}{|c|}{$\mathbf{0.95}$} & $\mathbf{0.951}$  & $\mathbf{0.948}$  & $\mathbf{0.950}$ & $\mathbf{0.947}$ & $\mathbf{0.950}$ & $\mathbf{0.941}$ \\[6pt]
\multicolumn{1}{|c|}{$\mathbf{0.99}$} \B & $\mathbf{0.992}$  & $\mathbf{0.991}$  & $\mathbf{0.991}$ & $\mathbf{0.989}$ & $\mathbf{0.991}$ & $\mathbf{0.989}$ \\
\hline
\end{tabular}
\caption{Percentile comparison of $F_{1}$ vs. empirical
distributions for $100\times 100$ and $100\times 400$ Wishart
matrices with identity covariance.}\label{table}
\end{table}
}
\end{center}
\clearpage
\noindent{\textbf{\large Acknowledgments:}}  The authors wish to thank Harold Widom without whom none of
this would have been possible.  We thank John Harnad for the invitation to participate in the
Program on \textit{Random Matrices, Random Processes and Integrable Systems} at the Centre de recherches
math\'ematiques on the campus of the Universit\'e de Montr\'eal. This work was supported by the National Science
Foundation under grant DMS--0304414.


\begin{thebibliography}{10}
\bibitem{Pear1}
K.~Pearson,
\newblock{On the criterion that a given system of deviations from the probable in the case
of a correlated system of variables is such that it can be reasonably supposed to have
arisen from random sampling},
\newblock{\em Phil.  Mag.}, V, 50:157, 1900.
\end{thebibliography}
\end{document}